\documentclass[12pt,a4paper,reqno]{article}
\usepackage[margin=2.2cm]{geometry}
\usepackage{latexsym}
\usepackage{mathrsfs}
\usepackage{array,tabularx}
\usepackage{tocloft}
\cftsetindents{section}{0em}{2em}
\cftsetindents{subsection}{0em}{2em}
\usepackage{amsmath,amssymb,booktabs,tabularx,array}
\AtBeginDocument{
  \setlength{\abovedisplayskip}{6pt}
  \setlength{\belowdisplayskip}{6pt}
  \setlength{\abovedisplayshortskip}{3pt}
  \setlength{\belowdisplayshortskip}{3pt}
}
\usepackage{etoolbox}
\apptocmd{\thebibliography}{%
    \setlength{\itemsep}{0pt}%
    \setlength{\parsep}{0pt}%
}{}{}
\usepackage{amsthm}
\usepackage{tikz}
\usepackage{tikz-cd}
\usepackage{hyperref}
\usepackage{enumerate}
\usepackage{enumitem}
\newcommand{\RN}[1]{%
  \textup{\uppercase\expandafter{\romannumeral#1}}%
}
\newcommand{\rn}[1]{%
  \textup{\lowercase\expandafter{\romannumeral#1}}%
}
\newtheorem{theorem}{Theorem}[section]
\newtheorem{proposition}[theorem]{Proposition}

\newtheorem{lemma}[theorem]{Lemma}

\newtheorem{corollary}[theorem]{Corollary}
\newtheorem{question}[theorem]{Question}

\newtheorem{definition}[theorem]{Definition}

\newtheorem{remark}[theorem]{Remark}

\newtheorem{example}[theorem]{Example}
\usepackage{hyperref}
\usepackage{enumitem}
\usepackage[title]{appendix}
\hypersetup{
colorlinks=true,linkcolor=blue,citecolor=blue,urlcolor=blue,citebordercolor={0 0 1},urlbordercolor={0 0 1},linkbordercolor={0 0 1}}
 
\newcommand{\ZZ}{\mathbb{Z}}

\newcommand{\CC}{\mathbb{C}}
\newcommand{\RR}{\mathbb{R}}
\newcommand{\PP}{\mathbb{P}}

\newcommand{\CP}{\mathbb{CP}}

\title{Geography of Klein bottle Luttinger surgery}
\author{Xinle Dai, Boxi Hao, and Shengzhen Ning}

\date{\today}

\begin{document}

\maketitle

\begin{abstract}
We propose and study a modification of Auroux's uniqueness question on Luttinger surgery by replacing Lagrangian tori with Klein bottles. In Kodaira dimension $-\infty$, we introduce two simple graphs whose connected components give bounds and, in several cases, exact counts of Klein bottle Luttinger equivalence classes. In Kodaira dimension zero, we provide examples of numerical invariants supporting infinitely many equivalence classes and prove a homological non-triviality theorem for Lagrangian Klein bottles in all symplectic Calabi--Yau surfaces with positive Euler characteristic. In positive Kodaira dimensions, we construct examples in which Klein bottle Luttinger surgery changes the symplectic Kodaira dimension.
\end{abstract}
\tableofcontents
\section{Introduction}
\subsection{Auroux's question}\label{sec:Auroux-question}
Luttinger surgery \cite{Luttinger1995,EP94,AurouxDonaldsonKatzarkov2003} was originally introduced as a symplectic cut-and-paste operation along an embedded Lagrangian torus $L$ in a symplectic four-manifold. It removes a Weinstein neighborhood $D^*L$ and reglues it by a contactomorphism of $S^*L\cong T^3$ with its canonical contact structure. In the Lagrangian framing, this map fixes the longitude classes and sends the meridian class $[\mu]$ to $[\mu]+k[\gamma']$, where $k\in\mathbb{Z}$ and $[\gamma']$ is the class of the Lagrangian push-off of an essential simple closed curve $\gamma\subset L$. The resulting symplectic form agrees with the original outside the surgery region, while the global topology may change. It has become a powerful tool for constructing symplectic $4$-manifolds, including exotic examples, and revealing the diversity of their topology and smooth structures.

For a closed symplectic $4$-manifold $(X,\omega)$, we define its \textbf{Auroux data} by
\[\mathcal{A}(X,\omega)=\bigl(c_1^2,c_2,c_1\cdot[\omega],[\omega]^2\bigr).
\]
Here $c_1$ is the first Chern class of any $\omega$-compatible almost-complex structure, $c_2$ is the Euler characteristic of $X$, and “integral” means that $[\omega]$ is the image of a class in $H^2(X;\ZZ)$. Auroux posed the following question in \cite[Question~2.6]{Auroux2005OpenQuestions}.

\begin{question}[Auroux]\label{question:Auroux}
Let $(X_1,\omega_1)$ and $(X_2,\omega_2)$ be closed integral symplectic $4$-manifolds with the same Auroux data. Is $(X_2,\omega_2)$ symplectomorphic to a symplectic manifold obtained from $(X_1,\omega_1)$ by a finite sequence of Luttinger surgeries?
\end{question}

The formulation of the above question in McDuff--Salamon  \cite[Problem~11]{McDuffSalamon2017} imposes the additional requirement that the surgeries be performed along pairwise disjoint Lagrangian tori. Surgeries along a disjoint collection can be performed successively, so the {\bf disjoint version} is more restrictive than the {\bf sequential version}. In Appendix \ref{section:two notions}, we show that a hyperelliptic surface admits two integral symplectic forms with the same Auroux data that are related by sequential Luttinger surgeries, but cannot be related by disjoint Luttinger surgeries. In the rest of the paper, we will only consider the sequential version. The related question  in the smooth category has an affirmative answer under the simply-connected assumption by Baykur--Sunukjian \cite{BaykurSunukjian}.

Recently, Question \ref{question:Auroux} was answered negatively in \cite{iida2026}. More precisely, Iida constructs  integral symplectic forms with the same Auroux data on \(\CP^2\#10\overline{\CP}^2\) and \(S\#2\overline{\CP}^2,\) where $S$ is the Craighero--Gattazzo numerical Godeaux surface. These manifolds are homeomorphic, but their symplectic Kodaira dimensions are $-\infty$ and $2$, respectively. They therefore cannot be related by Luttinger surgeries, since Luttinger surgeries preserve symplectic Kodaira dimension \cite[Theorem~1.1]{HoLi2012}.

There are also counterexamples for which the symplectic Kodaira dimensions agree. We give some  examples addressing the additional conditions proposed in \cite[Question~4.1]{iida2026}.

\begin{example}[Obstruction from minimality]\label{example:Auroux-minimality}
Let $m>n>0$ be integers. Equip $X_0=S^2\times S^2$ with the product symplectic form $\omega_0$ with \([\omega_0]=mA+2nB,\) where $A$ and $B$ denote the Poincar\'e duals of the two $S^2$-factors. On $X_1=\CP^2\#\overline{\CP}^2$, choose a blowup symplectic form $\omega_1$ with \([\omega_1]=(m+n)H-(m-n)E,\) where $H$ and $E$ denote the Poincar\'e duals of the line and exceptional classes. Both forms are integral and have the same Auroux data:
\[\mathcal{A}(X_0,\omega_0)=\mathcal{A}(X_1,\omega_1)=(8,4,2m+4n,4mn).\]
However, $(X_0,\omega_0)$ is symplectically minimal, whereas $(X_1,\omega_1)$ contains a symplectic exceptional sphere of area $m-n>0$. Since Luttinger surgery preserves symplectic minimality \cite{Welschinger} \cite[Proposition~3.1]{HoLi2012}, these two integral symplectic manifolds have the same Auroux data but cannot be related by a sequence of Luttinger surgeries.
\end{example}

\begin{example}[Minimal irrational ruled surfaces]
\label{example:Auroux-minimal-ruled}
Fix a closed Riemann surface $\Sigma_g$ of genus $g>0$, and let $\mathcal{L}$ be a holomorphic line bundle of degree $-1$ over $\Sigma_g$. Consider \(X_0=\Sigma_g\times\CP^1,\,X_1=\mathbb P(\mathcal O\oplus \mathcal{L}).\)
These are the trivial and nontrivial $S^2$-bundles over $\Sigma_g$. For $i=0,1$, let $S_i$ be a section and $F$ a fiber, chosen so that \(S_i^2=-i,S_i\cdot F=1,F^2=0.\) For integers $m>n>0$, we may choose integral symplectic forms with \([\omega_0]=2nS_0+mF,\,[\omega_1]=2nS_1+(m+n)F.\) Since \(c_1(X_i)=2S_i+(2-2g+i)F,\) a direct calculation gives
\[\mathcal{A}(X_0,\omega_0)=\mathcal{A}(X_1,\omega_1)=(8-8g,\;4-4g,\;4n(1-g)+2m,\;4mn).\]
Both manifolds are symplectically minimal. However, by \cite[Theorem~1.2]{HoLi2012}, every Luttinger surgery on a symplectic four-manifold with Kodaira dimension $-\infty$ leaves the manifold unchanged. Consequently, imposing minimality does not suffice to make Question~\ref{question:Auroux} true.
\end{example}

\begin{example}[Two forms on the same rational surface]
\label{example:Auroux-two-blowups}
Let \(X=\CP^2\#2\overline{\CP}^2\) and choose blowup symplectic forms $\omega_0,\omega_1$ with \([\omega_0]=5H-E_1-E_2,\,[\omega_1]=6H-3E_1-2E_2.\) Their Auroux data agree:
\[\mathcal{A}(X,\omega_0)=\mathcal{A}(X,\omega_1)=(7,5,13,23).\]
Nevertheless, these forms are not symplectomorphic by \cite[Lemma~1.10]{KK17}. On the other hand, every Luttinger surgery on a symplectic four-manifold with Kodaira dimension $-\infty$ leaves its symplectomorphism class unchanged \cite[Theorem~1.2]{HoLi2012}. Thus $(X,\omega_0)$ and $(X,\omega_1)$ cannot be related by Luttinger surgeries, despite having the same underlying smooth manifold.
\end{example}

We next give a counterexample in which both forms are defined on the \emph{same minimal} smooth four-manifold. The obstruction comes from the restrictions on the lattice automorphisms.

\begin{example}[Enriques surface]
\label{example:Auroux-Enriques}
Let $X$ be the underlying smooth four-manifold of an Enriques surface. Consider the Enriques lattice \(\Lambda=H^2(X;\mathbb Z)/\operatorname{Tor}\cong U\oplus E_8(-1)\) and its automorphism group  $O(\Lambda).$ Choose $D_0,D_1\in \Lambda$ such that \(D_0^2=D_1^2=d>0,\) but $D_0$ and $D_1$ lie in distinct orbits of $O(\Lambda)$. Such pairs exist already in the hyperbolic summand $U\subset\Lambda$. Indeed, choose a basis $s_1,s_2$ with $s_1^2=s_2^2=0$ and $s_1\cdot s_2=1$, and take \(D_0=6s_1+s_2,D_1=3s_1+2s_2.\) so that $D_0^2=D_1^2=12$. But by \cite[Proposition~3.4]{Ning2025}, the smallest positive pairings of $D_0$ and $D_1$ with integral square-zero classes are $1$ and $2$, respectively. Since lattice isometries preserve these pairings, $D_0$ and $D_1$ lie in distinct $O(\Lambda)$-orbits.

Note that every positive class in $H^2(X;\RR)$ admits a K\"ahler representative since a generic complex structure on $X$ contains no $(-2)$-curves; see \cite{Dolgachev2016}. Thus we may choose two symplectic forms $\omega_0,\omega_1$ representing $D_0,D_1$, respectively. They have the same Auroux data \[\mathcal{A}(X,\omega_0)=\mathcal{A}(X,\omega_1)=(0,12,0,12).\] Suppose that $(X,\omega_0)$ and $(X,\omega_1)$ were related by a sequence of Luttinger surgeries. By \cite[Theorems~1.1 and~1.3, Proposition~3.1]{HoLi2012}, every intermediate manifold is minimal of symplectic Kodaira dimension zero and has the same integral homology as $X$. In particular, every intermediate manifold has $b_2^+=1$. Corollary~\ref{cor:torus-period-lattice} therefore applies at each step. Composing the resulting lattice isometries gives an automorphism of $\Lambda$ which maps $\operatorname{PD}([\omega_0])$ to $\operatorname{PD}([\omega_1])$, contrary to our choice of $D_0,D_1$. Hence $(X,\omega_0)$ and $(X,\omega_1)$ cannot be related by Luttinger surgeries. This gives a negative answer to \cite[Question~4.1]{iida2026}, even for symplectic forms on the same manifold which is also minimal.
\end{example}

\subsection{Lagrangian Klein bottle Luttinger surgery}\label{sec:intro-KB-Luttinger}

The preceding examples show that Auroux data do not determine Luttinger surgery equivalence in many cases, even after imposing minimality and fixing the underlying smooth manifold. Nevertheless, Luttinger surgery can also be performed along Lagrangian Klein bottles \cite{Nemirovski2009,Shevchishin2009}. In contrast to surgery along Lagrangian tori, the two manifolds in Examples~\ref{example:Auroux-minimality}, \ref{example:Auroux-minimal-ruled}, and \ref{example:Auroux-two-blowups} are all related by a single Klein bottle (KB) Luttinger surgery. We refer to Adaloglou--Evans \cite[\S 1.1, \S 3.1, \S 3.2]{AdaloglouEvans2026} for Lagrangian Klein bottles in $S^2\times S^2$ whose surgery produces $\CC\PP^2\#\overline{\CP}^2$, and Sections \ref{sec:rational surfaces} and \ref{sec:irrational-ruled} for the other two examples.

All the counterexamples to Auroux's question in the previous section, including Iida's one \cite{iida2026}, rely on various properties for torus Luttinger surgery established by Ho--Li \cite{HoLi2012}, which in turn build on Welschinger's result \cite{Welschinger} using SFT compactness \cite{SFTcompactness}. This leads us to investigate which of these properties extend to Klein bottle Luttinger surgery and what new phenomena could arise. Comparisons between torus and Klein bottle Luttinger surgery are summarized in Table~\ref{tab:luttinger-comparison}, where all the assertions can be found in \cite{HoLi2012} or the subsequent sections of this paper. Remarkably, in Section~\ref{sec:KB-positive-Kodaira}, we use symplectic cyclic branched covers to show the following.

\begin{theorem}\label{thm:Kod-change}
    Klein bottle Luttinger surgery can change symplectic Kodaira dimensions.
\end{theorem}

\begin{table}[htbp]
\centering
\footnotesize
\renewcommand{\arraystretch}{1.15}
\setlength{\tabcolsep}{4pt}
\renewcommand{\tabularxcolumn}[1]{m{#1}}
\begin{tabularx}{0.95\linewidth}{
    |>{\centering\arraybackslash}m{0.34\linewidth}
    |>{\centering\arraybackslash}X
    |>{\centering\arraybackslash}X|}
\hline
\textbf{Property}
& \textbf{Torus Luttinger surgery}
& \textbf{KB Luttinger surgery}
\\ \hline
Symplectic minimality
& Preserved
& Can change
\\ \hline
Symplectic Kodaira dimension $\kappa$
& Preserved
& Can change
\\ \hline
$\kappa=-\infty$ symplectomorphism type
& Preserved
& Can change
\\ \hline
Symplectic Calabi--Yau (Torsion $c_1$)
& Preserved
& Preserved
\\ \hline
Auroux data $\mathcal{A}$
& Preserved
& Preserved
\\ \hline
$b_1,\ b_2,\ b_2^\pm$
& Can change
& Preserved
\\ \hline
Rational intersection form
& Can change
& Preserved
\\ \hline
Integral intersection form
& Can change
& Can change
\\ \hline
\end{tabularx}
\caption{Torus versus Klein bottle Luttinger surgery.}
\label{tab:luttinger-comparison}
\end{table}

For a closed integral symplectic four-manifold $(X,\omega)$, we define its \textbf{extended Auroux data} by
\[\widetilde{\mathcal{A}}(X,\omega)\:=\bigl(c_1^2,c_2,c_1\cdot[\omega],[\omega]^2,b_1\bigr),\]
where $b_1=\dim_{\mathbb{R}}H_1(X;\mathbb{R})$. This is motivated by the fact that KB Luttinger surgery preserves the first Betti number; see Lemma \ref{lemma:surgery-betti-change}. 
A $5$-tuple in $\ZZ^5$ is called {\bf realizable} if it can be realized as the extended Auroux data of a symplectic four-manifold. Realizable tuples have to satisfy several well-known conditions, such as $c_1^2+c_2\equiv 0\pmod {12},c_1\cdot[\omega]\equiv [\omega]^2\pmod 2$. Given a realizable integral $5$-tuple $\widetilde{\mathcal A}$, define the count
\[\mathcal N_{\widetilde{\mathcal A}} \:=\#\left(\left\{ (X,\omega)\middle|\widetilde{\mathcal A}(X,\omega)=\widetilde{\mathcal A}\right\}\big/\sim_{KB}\right),\]
where $(X,\omega)$ ranges over closed integral symplectic four-manifolds, and $\sim_{KB}$ is the equivalence relation generated by symplectomorphisms and Klein bottle Luttinger surgeries. We refer to the determination of these counts as the {\bf geography problem} for Klein bottle Luttinger surgeries. This terminology is motivated by the geography problem for symplectic four-manifolds, which asks which numerical topological invariants are realizable. For the Luttinger surgery question, we also incorporate the extra symplectic data $c_1\cdot[\omega],[\omega]^2$ and aim to count the equivalence classes realizing each fixed tuple.

The problem is already delicate enough in symplectic Kodaira dimension $-\infty$. The underlying manifolds are $X_n=\mathbb{CP}^2\#n\overline{\mathbb{CP}}^{2}$, $n\geq2$, and $X_{g,k}=(S^2\times\Sigma_g)\#k\overline{\mathbb{CP}}^{2}$, $g,k\geq1$. In Sections~\ref{subsec:rational-KB-graphs} and~\ref{sec:ruled-KB-graph}, we will introduce two surgery graphs $\mathbb G_{n;P,Q}^{\star}$ and $\mathbb G_{g,k;P,Q}^{\star}$ for rational and irrational ruled manifolds respectively, where $\star\in\{\mathrm{alg},\mathrm{geo}\}$.  For fixed $P\in\ZZ,Q\in\ZZ_{>0}$, where $P>0$ in the rational case, their vertices are classes $x\in H_2(X;\mathbb Z[1/2])$ in the $c_1$-symplectic cone satisfying $c_1\cdot x=P$ and $x^2=Q$. Both graphs contain Weyl edges corresponding to changes of marking induced by diffeomorphisms. Their additional Luttinger edges record algebraically permitted lattice reflections induced by virtual Luttinger surgeries and the subset admitting actual Lagrangian Klein bottle realizations, respectively. We write $\mathbb X(n,P,Q)$ and $\mathbb Y(n,P,Q)$, or $\mathbb X(g,k,P,Q)$ and $\mathbb Y(g,k,P,Q)$, for the unions of the connected components containing at least one integral vertex in the corresponding algebraic and geometric graphs, respectively. Our results of the KB Luttinger surgery geography problem in symplectic Kodaira dimension $-\infty$ can be summarized as the following.

\begin{theorem}
\label{thm:intro-KB-rational}
For each $n\geq 2$, $P,Q\in\mathbb Z_{>0}$ with realizable \(\widetilde{\mathcal A}=(9-n,\,3+n,\,P,\,Q,\,0),\)  
\[\#\pi_0\bigl(\mathbb X(n,P,Q)\bigr)\leq\mathcal N_{\widetilde{\mathcal A}}\leq\#\pi_0\bigl(\mathbb Y(n,P,Q)\bigr)<\infty.\]
Moreover,
\begin{itemize}[nosep]
\item $\mathcal N_{\widetilde{\mathcal A}}=1$ when $n=0$ or $1$.
\item For $2\leq n\leq10$, \(\mathcal N_{\widetilde{\mathcal A}}=\#\pi_0\bigl(\mathbb Y(n,P,Q)\bigr).\)
\item For $2\leq n\leq4$, \(\mathcal N_{\widetilde{\mathcal A}}=\#\pi_0\bigl(\mathbb X(n,P,Q)\bigr).\) In particular, when $n=2$, this count admits an explicit formula in terms of the numbers of integral solutions to Diophantine equations subject to specified inequalities and congruence conditions.
\end{itemize}
\end{theorem}

\begin{theorem}
\label{thm:intro-KB-ruled}
Let $g\geq1$, $k\geq 1$, $P\in\mathbb Z$, and $Q\in\mathbb Z_{>0}$, and suppose that
\[\widetilde{\mathcal A}=(8-8g-k,\,4-4g+k,\,P,\,Q,\,2g)\] 
is realizable. If $g=1$, further assume that $P>-k$. Then we have
\[\#\pi_0\bigl(\mathbb X(g,k,P,Q)\bigr)\leq\mathcal N_{\widetilde{\mathcal A}}=\#\pi_0\bigl(\mathbb Y(g,k,P,Q)\bigr)<\infty.\]
Moreover,
\begin{itemize}[nosep]
\item$\mathcal N_{\widetilde{\mathcal A}}=1$ when $k=0$.
\item For $1\leq k\leq3$, \(\mathcal N_{\widetilde{\mathcal A}}=\#\pi_0\bigl(\mathbb X(g,k,P,Q)\bigr).\)
In this range, the count admits an explicit formula in terms of the number of integral solutions to a quadratic Diophantine equation subject to specified inequalities.
\end{itemize}
\end{theorem}

In symplectic Kodaira dimension zero, we show fixed extended Auroux data can support infinitely many KB Luttinger equivalence classes.

\begin{theorem}
\label{thm:intro-KB-infty}
For every integer $q\geq1$ and $b\in\{2,3\}$, \(\mathcal N_{(0,\,0,\,0,\,2q,\,b)}=\infty.\)
\end{theorem}

The families in Theorem~\ref{thm:intro-KB-infty} have vanishing Euler characteristic. For positive Euler characteristic, we obtain a complementary result concerning the homology classes of embedded Lagrangian Klein bottles proved using Luttinger surgery.

\begin{theorem}[Homological non-triviality]\label{thm:KB-homological-essentiality}
Let $(X,\omega)$ be a closed symplectic four-manifold with $c_1\cdot[\omega]=0$ and positive Euler characteristic  $\chi(X)>0$. Then the $\ZZ_2$-homology class of any embedded Lagrangian Klein bottle $L\subset (X,\omega)$ is not the reduction of an integral torsion class. More precisely, \([L]_2\notin\rho\bigl(\operatorname{Tor}H_2(X;\ZZ)\bigr),\) where $\rho:H_2(X;\ZZ)\rightarrow H_2(X;\ZZ_2)$ is the coefficient reduction map.
\end{theorem}

The homological non-vanishing theorem of Shevchishin \cite{Shevchishin2009} and Nemirovski \cite[Theorem A]{Nemirovski2009} asserts that every embedded Lagrangian Klein bottle in a closed symplectic four-manifold with $c_1\cdot[\omega]>0$ represents a nonzero $\ZZ_2$-homology class. Theorem~\ref{thm:KB-homological-essentiality} extends this conclusion to the case $c_1\cdot[\omega]=0$ under the additional assumption $\chi(X)>0$, and strengthens homological non-vanishing by excluding the reductions of all integral torsion classes. It also extends Smirnov's result for symplectic homotopy K3 surfaces with vanishing first Chern class \cite[Theorem 5]{Smirnov}. In particular, for an Enriques surface it gives \([L]_2\notin\{0,\rho(\tau)\},\) where $\tau\in H_2(X;\ZZ)$ is the unique nonzero $\ZZ_2$-torsion class. 

\medskip
\noindent
{\bf Organization of the paper.}
Section \ref{section:torus-KB-comparison} compares the topological effects of torus and Klein bottle Luttinger surgery. Sections \ref{sec:rational surfaces} and \ref{sec:irrational-ruled} introduce the surgery graphs for the rational and irrational ruled cases, respectively, and establish the corresponding counting results in Theorems \ref{thm:intro-KB-rational} and \ref{thm:intro-KB-ruled}. Section \ref{sec:Kod-0} proves the homological non-triviality Theorem \ref{thm:KB-homological-essentiality} and constructs examples with infinite KB Luttinger counts of Theorem \ref{thm:intro-KB-infty} in symplectic Kodaira dimension zero. Section \ref{sec:KB-positive-Kodaira} aims to construct KB Luttinger surgeries changing symplectic Kodaira dimension from one to two. Finally, Appendix \ref{section:two notions} distinguishes sequential torus Luttinger surgery from surgery along a fixed disjoint collection of Lagrangian tori through the example of two symplectic forms on the same hyperelliptic surface.

\medskip
\noindent
{\bf Acknowledgements.}
We are grateful to Denis Auroux and Tian-Jun Li for their interest in this work, their encouragement, and many helpful discussions.

\medskip
\noindent
{\bf AI disclosure.} 
This project originated from the unexpected AI's discovery (independent with \cite{iida2026}) of various examples in Section \ref{sec:Auroux-question}, together with the clarification of the two notions discussed in Appendix~\ref{section:two notions}. After communicating these counterexamples to Denis Auroux in mid August 2026, we were encouraged to explore ways to salvage his question. The idea of investigating the Klein bottle version of Auroux's question was purely inspired by the illuminating paper \cite{AdaloglouEvans2026}. Besides cosmetic language and diagram editing, the other essential math inputs from LLM tools are: pointing out Remark \ref{rmk:Furuta-refinement}; bringing \cite{KirbyTaylor1990} and other relevant literature used in Section \ref{section:positive-Euler-Kod-0} to our attention; motivating the construction in Section \ref{sec:KB-positive-Kodaira} (see Remark \ref{rmk:cyclic-origination}). The authors are responsible for the mathematical content and the final text of the paper.

\section{A comparison between torus and Klein bottle Luttinger surgery}
\label{section:torus-KB-comparison}

In this section, we compare the topological and cohomological effects of torus and Klein bottle Luttinger surgery. Alternative proofs of some conclusions can also be found in \cite{HoLi2012,Smirnov,AdaloglouEvans2026}.  We first recall the definitions of two surgeries and fix notation that will be used throughout the paper. For detailed local constructions, we refer to \cite{Luttinger1995,EP94,AurouxDonaldsonKatzarkov2003} in the torus case and \cite{Nemirovski2009,Smirnov,AdaloglouEvans2026} in the Klein bottle case.

Let $L\subset(X,\omega)$ be an embedded Lagrangian torus or Klein bottle. Let $N$ be a closed Weinstein neighborhood of $L$, identified with a disk neighborhood of the zero section in $T^*L$. Set
\[C=X\setminus\operatorname{int}N,\qquad Y=\partial N.\]
We have a canonical identification $Y\cong \partial C$ using the original gluing map. Luttinger surgery removes $N$ and reglues it by a boundary twist $f:Y\to \partial C$, producing
\[X'=C\cup_f N.\]
Let $\mu\subset Y$ be an oriented meridian circle of $L$ and take $\mu'=f(\mu)\subset \partial C$ to be the new meridian circle inside $X'$. By the canonical identification between $\partial C$ and $Y$, we can view $[\mu']$ as a class in $H_1(Y;\ZZ)$. When $L$ is a torus, choose an oriented essential simple closed curve $\gamma\subset L$, and let $\gamma'\subset Y$ be its push-off in the Lagrangian framing. For an integer $k$, the gluing map $f$ of a $k$-fold Luttinger twist in the $\gamma$-direction sends
\[[\mu]\longmapsto[\mu']=[\mu]+k[\gamma']\in H_1(Y;\ZZ).\]
When $L$ is a Klein bottle, view it as a nontrivial $S^1$-bundle over $S^1$ and let $\gamma\subset L$ be a fiber circle. Denote by $\gamma'\subset Y$ its $\mathbb C$-normal push-off in the sense of \cite{Nemirovski2009}. Consider the orientation double cover $\pi:T^2\to L$ and let $\tilde{\gamma}$ be a component of $\pi^{-1}(\gamma)$. Then the $1$-fold torus Luttinger twist in the $\tilde{\gamma}$-direction can be chosen to be equivariant with respect to the cotangent $\ZZ_2$-action induced by the deck transformation of $\pi$. It descends to the Klein bottle Luttinger twist $f$ which sends
\[[\mu]\longmapsto [\mu']=[\mu]+[\gamma']\in H_1(Y;\ZZ).\]
In both cases, the local twist can be chosen to give a symplectic identification of boundary collars, so $X'$ carries a symplectic form $\omega'$ agreeing with $\omega$ on $C$. The zero section of the reglued Weinstein neighborhood gives a Lagrangian torus or Klein bottle $L'\subset(X',\omega')$. Although Klein bottle Luttinger surgery can also be discussed for any integer coefficient $k$ \cite[\S 2.1]{Nemirovski2009}, throughout this paper, we only refer to the $k=1$ surgery as in \cite{AdaloglouEvans2026}. 

For the manifolds $X$ and $X'$, we write
\begin{align*}
    b_j(X)=\dim_{\mathbb R}H_j(X;\mathbb R),\qquad
h_j(X)=\dim_{\mathbb Z_2}H_j(X;\mathbb Z_2),
\end{align*}
and $Q_X$ and $Q_{X'}$ for their intersection forms on the lattices 
\[\Lambda=H_2(X;\mathbb Z)/\operatorname{Tor},\qquad\Lambda'=H_2(X';\mathbb Z)/\operatorname{Tor}.\]
Denote by \(i:C\hookrightarrow X,\, i':C\hookrightarrow X'\) the inclusions of the common Weinstein neighborhood complement $C$ and $\Lambda_C\subset \Lambda,\,\Lambda'_C\subset \Lambda'$ the images of $H_2(C;\ZZ)$ under $i_*,\,i'_*$.

\begin{lemma}\label{lemma:surgery-characteristic-numbers}
Both torus and Klein bottle Luttinger surgery preserve Euler characteristic $\chi$ and signature $\sigma$. Consequently, they preserve $c_1^2$ and $c_2$, and satisfy
\[b_2(X')-b_2(X)=2\bigl(b_1(X')-b_1(X)\bigr),\qquad b_2^\pm(X')-b_2^\pm(X)=b_1(X')-b_1(X).\]
\end{lemma}
\begin{proof}
Euler characteristic is additive under gluing along the closed three-manifolds $Y$ and $\partial C$, and the gluing pieces $C$ and $N$ are the same. Hence $\chi(X')=\chi(X)$. In addition, Novikov additivity similarly gives $\sigma(X')=\sigma(X)$.

For any symplectic four-manifold $X$, we have \(c_2(X)=\chi(X)\) and \(c_1(X)^2=2\chi(X)+3\sigma(X).\) The remaining equations thus follow from \(\chi(X)=2-2b_1(X)+b_2(X)\) and \(b_2^\pm(X)=\frac{b_2(X)\pm\sigma(X)}{2}.\)
\end{proof}

For any coefficient ring $R$, we will frequently use the following isomorphisms combining excision with Poincar\'e--Lefschetz duality:
\begin{equation}\label{equ:Excision and Poincare dual}
\begin{aligned}
   H_j(X,C; R)\cong H_j(N,Y;R)\cong H^{4-j}(N;R)\cong H^{4-j}(L;R),\\
   H^j(X,C; R)\cong H^j(N,Y;R)\cong H_{4-j}(N;R)\cong H_{4-j}(L;R).
\end{aligned}
\end{equation}
Note that \eqref{equ:Excision and Poincare dual} still holds if we replace $X$ by $X'$.

The effect of Luttinger surgery on first homology can be read from the relative exact sequence for the pair $(X,C)$. By \eqref{equ:Excision and Poincare dual}, we have
\[
H_1(X,C;R)=0,
\qquad
H_2(X,C;R)\cong H^2(L;R)=
\begin{cases}
 R,&L\text{ is a torus},\\
 R/2R,&L\text{ is a Klein bottle}.
\end{cases}
\]
The relevant part of the long exact sequence for the pair $(X,C)$ is therefore
\begin{equation}\label{equ:LES}
    H_2(X;R)
\longrightarrow H^2(L;R)\cong
H_2(X,C;R)
\xrightarrow{\;\partial\;}
H_1(C;R)
\longrightarrow H_1(X;R)
\longrightarrow0.
\end{equation}
 Let $D\subset N$ be an oriented normal disk, with $\partial D=\mu\subset Y$. Its relative class generates $H_2(N,Y;R)$ as an $R$-module since the intersection number between $D$ and $L$ is $1\in R$. Under excision, the connecting homomorphism therefore sends this generator to $j_*[\mu]\in H_1(C;R)$. Thus
\(
\operatorname{im}\partial=\langle j_*[\mu]\rangle,
\)
where $j:Y\hookrightarrow C$ is the inclusion.
Applying the same argument after surgery, with $\mu'$ the
meridian of the new filling viewed in $Y$, gives
\[
H_1(X;R)
\cong
\frac{H_1(C;R)}
{\langle j_*[\mu]\rangle},
\qquad
H_1(X';R)
\cong
\frac{H_1(C;R)}
{\langle j_*[\mu']\rangle}.\]
Thus the change in first homology is determined by replacing the old meridian relation given by $j_*[\mu]$ with the new one given by $j_*[\mu']$. Taking $R=\RR$ yields the following lemma of the pattern for $b_1$ change. This recovers \cite[Proposition~4.4]{HoLi2012} in the torus Luttinger surgery case.
\begin{lemma}\label{lemma:surgery-betti-change}
For both torus and KB Luttinger surgeries, let \(m=j_*[\mu],\,a=j_*[\gamma']\in H_1(C;\mathbb R).\)  Then
\[b_1(X')-b_1(X)=
\begin{cases}
-1,&m=0\text{ and }m+ka\neq0,\\
 1,&m\neq0\text{ and }m+ka=0,\\
 0,&\text{otherwise}.
\end{cases}
\]
In particular, for Klein bottle surgery, we always have $b_1(X')=b_1(X)$ since $m=a=0$.
\end{lemma}

Another way to see the Betti number change is by checking whether the homology classes of the Lagrangians $L$ and $L'$ are trivial or not. Let us define (similarly for $L'\subset X'$)
\[\varepsilon_X^{\mathbb R}(L)=
\begin{cases}
1,&L\text{ is a torus and }
   [L]=0\in H_2(X;\mathbb R),\\
0,&\text{otherwise},
\end{cases}\quad
\varepsilon_X^{\mathbb Z_2}(L)=
\begin{cases}
1,&[L]_2=0\text{ in }H_2(X;\mathbb Z_2),\\
0,&[L]_2\neq0\text{ in }H_2(X;\mathbb Z_2).
\end{cases}
\]

\begin{lemma}\label{lem:betti-change-Lagrangian-class}
For both torus and KB Luttinger surgeries, 
\[b_1(X')-b_1(X)=\varepsilon_X^{\mathbb R}(L)-\varepsilon_{X'}^{\mathbb R}(L'),\qquad h_1(X')-h_1(X)=\varepsilon_X^{\mathbb Z_2}(L)-\varepsilon_{X'}^{\mathbb Z_2}(L').\]
\end{lemma}

\begin{proof}
First take $R=\mathbb R$ in \eqref{equ:LES}. If $L$ is a torus, the first map is then $x\mapsto x\cdot[L]$. By nondegeneracy of the intersection pairing on $H_2(X;\mathbb R)$, this map vanishes exactly when $[L]$ is torsion. Its cokernel, and hence $\operatorname{im}\partial$, therefore has rank $\varepsilon_X^{\mathbb R}(L)$. If $L$ is a Klein bottle, then $H^2(L;\mathbb R)=0$, so $\operatorname{im}\partial$ has rank zero. Thus, in either case, \(b_1(C)=b_1(X)+\varepsilon_X^{\mathbb R}(L).\)

Next take $R=\mathbb Z_2$. In either case, $H^2(L;\mathbb Z_2)\cong\mathbb Z_2$, and the first map is by taking the intersection number with $[L]_2$. By nondegeneracy of the mod-two intersection pairing, this map is surjective exactly when $[L]_2\neq0$. Consequently, \(h_1(C)=h_1(X)+\varepsilon_X^{\mathbb Z_2}(L).\) 

Applying the same arguments to $X'$ and subtracting gives the two equalities.
\end{proof}

\begin{remark}\label{remark:h1-increase}
   For the Klein bottle Luttinger surgery, Nemirovski \cite[Theorem~2.2]{Nemirovski2009} and Shevchishin \cite[\S 3]{Shevchishin2009} proved that \([L]_2=0 \Rightarrow[L']_2\neq0.\) Hence, in this case, one always has \(h_1(X')=h_1(X)+1\).
\end{remark}

\begin{example}\label{example:torus-betti-change}
Consider \(T^4=(\mathbb R/\mathbb Z)^4\) with \(\omega=dx_1\wedge dy_1+dx_2\wedge dy_2\) and the homologically essential Lagrangian torus \(L=\{y_1=y_2=0\}.\) Write \(T^4=S^1_{x_2}\times\bigl(T^2_{x_1,y_1}\times S^1_{y_2}\bigr)\) so that $L$ is viewed as the product of $S^1_{x_2}$ with the $x_1$-circle in the torus fiber over $y_2=0$. Luttinger surgery along $L$ in the $x_1$-direction is therefore the product with $S^1_{x_2}$ of Dehn surgery along this fiber circle. For $k\geq 1$, performing Luttinger surgery on $(L,x_1,k)$ will induce a change of the monodromy of the torus bundle over $S^1_{y_2}$ given by the $k$-th power of a Dehn twist $D_{x_1}$ along the $x_1$-circle. The resulting manifold is \(W_k\cong S^1_{x_2}\times M_{D^k_{x_1}},\) where $M_{D^k_{x_1}}$ is the mapping torus \(\bigl(T^2\times[0,1]\bigr)/\bigl((z,1)\sim(D_{x_1}^kz,0)\bigr);\) see the example at the end of \cite[Section~2.1]{AurouxDonaldsonKatzarkov2003}. We have \(H_1(W_k;\mathbb Z)\cong\mathbb Z^3\oplus\mathbb Z/k\mathbb Z.\) Thus $b_1$ decreases from $4$ to $3$, $b_2$ decreases from $6$ to $4$, and both $b_2^+$ and $b_2^-$ decrease from $3$ to $2$. These $W_k$'s are the underlying smooth manifolds of primary Kodaira surfaces; see \cite[Table~1(b)]{Geiges1992}. In particular, $W_1$ is diffeomorphic to the Kodaira--Thurston manifold.
\end{example}

The next lemma concerns the change of the intersection forms.

\begin{lemma}\label{lem:intersection-form}
Let $(X',\omega')$ be obtained from a closed symplectic four-manifold $(X,\omega)$ by Luttinger surgery along a Lagrangian torus or Klein bottle $L$. If $L$ is a torus, we further assume $\varepsilon_X^\RR(L)=\varepsilon_{X'}^\RR(L')=1$. Then there is an isometry \(\Phi:H_2(X;\mathbb R)\rightarrow H_2(X';\mathbb R)\) such that
\begin{itemize}[nosep]
    \item $\Phi\circ i_*=i'_*$,
    \item $\Phi(\operatorname{PD}([\omega]))=\operatorname{PD}([\omega'])$,
    \item $\Phi(\operatorname{PD}(c_1(X,\omega)_\RR))=\operatorname{PD}(c_1(X',\omega')_\RR)$,
    \item $\Phi(\Lambda_C)=\Lambda_C'$.
\end{itemize}
In the torus case, \(\Lambda=\Lambda_C,\Lambda'=\Lambda_C'\) and $\Phi$ preserves the integral intersection forms; in the Klein bottle case, \([\Lambda:\Lambda_C]=[\Lambda':\Lambda'_C]\in\{1,2\}\), and if $\varepsilon_X^{\ZZ_2}(L)=1$, then both indices are $1$ and $\Phi$ preserves the integral intersection forms.
\end{lemma}

\begin{proof}

Applying \eqref{equ:Excision and Poincare dual} with $j=2$ and $R=\mathbb R$. In the Klein bottle case, we obtain \(H_2(X,C;\mathbb R)\cong H^2(L;\mathbb R)=0\); in the torus case, our assumption $\varepsilon_X^\RR(L)=1$ implies that $H_2(X;\RR)\to H_2(X,C;\RR)$ is the zero map.
The relative exact sequence for the pair $(X,C)$ therefore shows that
\(i_*:H_2(C;\mathbb R)\rightarrow H_2(X;\mathbb R)\) is surjective. The same argument applies to $i'_*$. Note that \(i_*(a)\cdot i_*(b)=i'_*(a)\cdot i'_*(b)\) for all $a,b\in H_2(C;\mathbb R)$ since the intersection numbers of classes coming from $C$ can be computed in $\operatorname{int}C$.
If $i_*(a)=0$, then $i'_*(a)$ is orthogonal to the image of $i'_*$. Since $i'_*$ is surjective and the intersection pairing on $X'$ is nondegenerate, we obtain $i'_*(a)=0$. Interchanging $X$ and $X'$ gives
\(\ker i_*=\ker i'_*.\) Indeed, the Mayer--Vietoris sequence can be used to further show that both kernels are equal to $\operatorname{im}(H_2(Y;\mathbb R)\to H_2(C;\mathbb R))$.
Consequently, \(\Phi(i_*(a)):=i'_*(a)\) defines an isometry which maps $\Lambda_C$ to $\Lambda_C'$ by the definition. 
Let $(\alpha,\alpha')$ denote either  $([\omega],[\omega'])$ or $(c_1(X,\omega)_\RR,c_1(X',\omega')_\RR)$.
Since the symplectic forms and symplectic tangent bundles agree on $C$, we have \(i^*\alpha=i'^*\alpha'\).
Thus, for every $b\in H_2(C;\mathbb R)$, we have
\[
\Phi(\operatorname{PD}(\alpha))\cdot i'_*(b)=\Phi(\operatorname{PD}(\alpha))\cdot \Phi(i_*(b))=
\operatorname{PD}(\alpha)\cdot i_*(b)=
\langle i^*\alpha,b\rangle=
\langle i'^*\alpha',b\rangle=
\operatorname{PD}(\alpha')\cdot i'_*(b).
\]
Surjectivity of $i'_*$ and nondegeneracy of the
intersection pairing give
\(\Phi(\operatorname{PD}(\alpha))=\operatorname{PD}(\alpha').\)
This proves the preservations of the Poincar\'e dual of the symplectic classes and the first Chern classes.

Now, consider the integral relative exact sequence. In the torus case, the relative homomorphism
\(
H_2(X;\mathbb Z)
\rightarrow H_2(X,C;\mathbb Z)\cong\mathbb Z
\) sends $A$ to $A\cdot[L]$.
Since $[L]=0\in H_2(X;\mathbb R)$, $A\cdot[L]=0$ for every $A\in H_2(X;\mathbb Z)$.
Exactness therefore implies that $H_2(C;\mathbb Z)\to H_2(X;\mathbb Z)$ is surjective.
Passing to the torsion-free quotient gives $\Lambda_C=\Lambda$.
Applying the same argument to $X'$ and $L'$ gives $\Lambda_C'=\Lambda'$.
In the Klein bottle case, since $H_2(X,C;\mathbb Z)\cong \ZZ_2$, we have $[\Lambda:\Lambda_C]\in\{1,2\}$.
The same argument gives
$[\Lambda':\Lambda_C']\in\{1,2\}$.
Since $\Lambda$ and $\Lambda'$ are unimodular and $\Lambda_C$ is isomorphic to $\Lambda_C'$ via $\Phi$,
\[[\Lambda:\Lambda_C]^2=|\det Q_{\Lambda_C}|=|\det Q_{\Lambda_C'}|=[\Lambda':\Lambda_C']^2.
\]
Thus the two indices are equal.
Finally, note that the map
\(H_2(X;\mathbb Z)\rightarrow H_2(X,C;\mathbb Z)\cong\mathbb Z_2\)
is given by mod-two intersection with $[L]_2$.
If $[L]_2=0$, this map vanishes, so
$H_2(C;\mathbb Z)\to H_2(X;\mathbb Z)$ is surjective.
It follows that $\Lambda_C=\Lambda$ and $\Lambda_C'=\Lambda'$.
\end{proof}

\begin{example}\label{example:KB-parity-change}
KB Luttinger surgery may replace one lattice by another containing the isomorphic index-two sublattice. Indeed, Adaloglou--Evans show that Luttinger surgery along any Lagrangian Klein bottle $L$ in $S^2\times S^2$ must produce $\CP^2\#\overline{\CP}^2$; see \cite[Proposition~3.1]{AdaloglouEvans2026}.
Their integral intersection forms are even and odd respectively. For $S^2\times S^2$, let $A,B$ be the two $S^2$-factor classes with \( [L]_2= B\bmod2.\)
The identification $\Phi$ in Lemma \ref{lem:intersection-form} satisfies
\[
\Phi(A)=\frac{H+E}{2},\qquad \Phi(B)=H-E,
\]
where $H,E$ are the line and exceptional classes after surgery.
The two integral lattices
\[
\Lambda=\mathbb ZA\oplus\mathbb ZB,\qquad
\Lambda'=\mathbb ZH\oplus\mathbb ZE
\]
contain the index-two sublattices
\[\Lambda_C=\mathbb Z(2A)\oplus\mathbb ZB,\qquad \Lambda_C'=\mathbb Z(H+E)\oplus\mathbb Z(H-E)\] such that $\Phi(\Lambda_C)=\Lambda_C'$.
\end{example}

\begin{lemma}\label{lem:pullback-intersection-lattices}
Under the hypotheses of Lemma \ref{lem:intersection-form}, assume that $H_1(X;\ZZ)$ is torsion-free.
Put \(\Gamma:=\Phi^{-1}(\Lambda')\subset H_2(X;\RR).\)
Then
\(
\Lambda\cap\Gamma=\Lambda_C.
\)
\end{lemma}

\begin{proof}
We only need to consider the case when $L$ is a Klein bottle and \([\Lambda:\Lambda_C]=[\Lambda':\Lambda_C']=2.\) Taking $R=\ZZ$ in \eqref{equ:LES}, we have
\[
H_2(X;\ZZ)\longrightarrow\ZZ_2
\longrightarrow H_1(C;\ZZ)
\longrightarrow H_1(X;\ZZ)\longrightarrow0.
\]
The first map is surjective by $[\Lambda:\Lambda_C]=2$.
Hence \(H_1(C;\ZZ)\cong H_1(X;\ZZ)\) is torsion-free.

The inclusion $\Lambda_C\subset\Lambda\cap\Gamma$ is obvious. Now let $x\in\Lambda\cap\Gamma$. Choose closed oriented embedded surfaces $S\subset X$ and $S'\subset X'$ representing $x$ and $\Phi(x)$ respectively. We may assume they are both transverse to $Y$ and cut them along $Y$.
Their portions $S_C$ and $S'_C$ in $C$ represent relative two-cycles in $(C,Y)$. For every $y\in H_2(C;\ZZ)$, we have
\[
[S_C]\cdot y
=x\cdot i_*y
=\Phi(x)\cdot i'_*y
=[S'_C]\cdot y,
\]
because $\Phi$ is an isometry and $\Phi\circ i_*=i'_*$.
Since $H_1(C;\ZZ)$ is torsion-free, Poincar\'e--Lefschetz duality and the universal coefficient
theorem identify
\[
H_2(C,Y;\ZZ)
\cong H^2(C;\ZZ)\cong
\operatorname{Hom}\bigl(H_2(C;\ZZ),\ZZ\bigr)
\]
by the intersection pairing. Thus $S_C$ and $S'_C$ represent the same integral relative class so that their boundaries are homologous in $Y$. Set
\(z:=[\partial S_C]=[\partial S_C']\in H_1(Y;\ZZ).\)
Write 
\[H_1(Y;\ZZ)\cong H_1(L;\ZZ)\oplus \ZZ_2\langle[\mu]\rangle=\ZZ\langle[\ell]\rangle\oplus\ZZ_2\langle [\gamma'] \rangle\oplus\ZZ_2\langle[\mu]\rangle.\]
The kernel of the inclusion from $Y$ into the original filling $N\subset X$ is $\langle[\mu]\rangle$.
Viewed from the same boundary $Y\subset C$, the kernel for the inclusion from $Y$ into the new filling $N\subset X'$ is \(f_*\langle[\mu]\rangle=\langle[\mu]+b\rangle.\)
Since the two boundary cycles $\partial S_C$ and $\partial S_C'$ bound the portions $S_N$ and $S_N'$ of $S$ and $S'$ in their respective fillings, we obtain \(z\in\langle[\mu]\rangle\cap\langle[\mu]+[\gamma']\rangle=\{0\}.\)
Hence, we can choose an integral two-chain $D$ in $Y$ with $\partial D=\partial S_C$. Then $S_C-D$ is a two-cycle in $C$. Since $S=S_C+S_N$, its difference from $S$ is \(S-(S_C-D)=S_N+D,\) which is a two-cycle in $N$. Since $H_2(N;\ZZ)=0$, the cycle $S_C-D$ represents $x$ in $X$. Therefore $x\in\Lambda_C$, proving the reverse inclusion.
\end{proof}

\begin{corollary}\label{cor:surgery-symplectic-pairings}
Let $(X',\omega')$ be obtained from a closed symplectic
four-manifold $(X,\omega)$ by Luttinger surgery along a
Lagrangian torus or Klein bottle. Then
\[c_1(X',\omega')\cdot[\omega']=c_1(X,\omega)\cdot[\omega],\qquad[\omega']^2=[\omega]^2.\]
\end{corollary}

\begin{proof}
For the Klein bottle case, both identities follow immediately from Lemma~\ref{lem:intersection-form}. Suppose that $L$ is a torus. The relative exact sequence gives $\operatorname{im}(H_2(C;\mathbb R)\xrightarrow{i_*}H_2(X;\mathbb R))=[L]^\perp$, and similarly $\operatorname{im}i'_*=[L']^\perp$. Since $L$ is Lagrangian,
\(
\langle c_1(X,\omega),[L]\rangle=\langle[\omega],[L]\rangle=0.
\)
We may therefore choose $a,b\in H_2(C;\mathbb R)$ such that
\(
i_*(a)=\operatorname{PD}(c_1(X,\omega)),
\,
i_*(b)=\operatorname{PD}([\omega]).
\)
Since the symplectic forms and symplectic tangent bundles agree
on $C$, for any $i_*'(c)\in [L']^\perp$, we have 
\[\begin{aligned}
    i'_*(a)\cdot i_*'(c)&=i_*(a)\cdot i_*(c)=\operatorname{PD}(c_1(X,\omega))\cdot i_*(c)=\operatorname{PD}(c_1(X',\omega'))\cdot i'_*(c),\\
     i'_*(b)\cdot i_*'(c)&=i_*(b)\cdot i_*(c)=\operatorname{PD}([\omega])\cdot i_*(c)=\operatorname{PD}([\omega'])\cdot i'_*(c).
\end{aligned}
\]
Consequently, both \( i'_*(a)-\operatorname{PD}(c_1(X',\omega'))\) and \(i'_*(b)-\operatorname{PD}([\omega'])\) belong to \(([L']^\perp)^\perp=\mathbb R[L'].\) Since \( [L']^2=0,\langle c_1(X',\omega'),[L']\rangle=0,\langle[\omega'],[L']\rangle=0,\) these differences do not contribute to the pairings. It follows that
\[c_1(X',\omega')\cdot[\omega']=i'_*(a)\cdot i'_*(b)=i_*(a)\cdot i_*(b)=c_1(X,\omega)\cdot[\omega],\quad[\omega']^2=i'_*(b)^2=i_*(b)^2=[\omega]^2.\]
\end{proof}

\begin{corollary}\label{cor:torus-period-lattice}
Let $(X',\omega')$ be obtained from a closed integral symplectic four-manifold $(X,\omega)$ by Luttinger surgery along a Lagrangian torus $L$, and suppose that \(b_2^+(X)=b_2^+(X')=1.\) Then there is a lattice isometry between $\Lambda$ and $\Lambda'$ which maps $\operatorname{PD}([\omega])$ to $\operatorname{PD}([\omega'])$ and $\operatorname{PD}(c_1(X,\omega))$ to $\operatorname{PD}(c_1(X',\omega'))$.
\end{corollary}

\begin{proof}
The Lagrangian torus $L$ satisfies \([L]^2=0,\, \langle[\omega],[L]\rangle=0.\)
Since $b_2^+(X)=1$, the orthogonal complement of $\operatorname{PD}([\omega])$ is negative definite. Thus $[L]=0$ in real homology. Similarly, $[L']=0$. The conclusion now follows from Lemma~\ref{lem:intersection-form}.
\end{proof}

\section{KB Luttinger surgery on rational manifolds}
\label{sec:rational surfaces}
In this section, we extend the discussion in Example \ref{example:KB-parity-change} to rational manifolds. Throughout this section, we write
\[
X_n=\CC\PP^2\#n\overline{\CC\PP}^{2},
\qquad
(h;e_1,\ldots,e_n)=hH-\sum_{k=1}^n e_kE_k,
\]
where $H,E_1,\ldots,E_n$ form the standard basis of $\Lambda=H_2(X_n;\ZZ)$, with
\[
H^2=1,\qquad E_k^2=-1,\qquad
H\cdot E_k=E_k\cdot E_\ell=0\quad(k\ne\ell).
\]
Let $c_1=3H-E_1-\cdots-E_n$ be a distinguished element in the lattice. We use the same notation for homology classes and their Poincar\'e duals.
For the symplectic form $\omega_{h;e_1,\cdots,c_n}$ on $X_n$, the parameters $h,e_1,\ldots,e_n$ are the symplectic areas of these classes.
Write
\[
P=3h-\sum_{k=1}^n e_k,
\qquad
Q=h^2-\sum_{k=1}^n e_k^2.
\]
Our goal is to understand the number $\mathcal{N}_{(9-n,\,3+n,\,P,\,Q,\,0)}$.

\subsection{The Adaloglou--Evans twist}\label{subsec:AE-twist}
For $S^2\times S^2$ with $A=[S^2\times\{p\}],B=[\{p\}\times S^2]$, let $\omega_{a,b}$ be the symplectic form which have symplectic areas $a$ and $b$ on $A$ and $B$ respectively.
 When $0<b<2a$, it contains a visible Lagrangian Klein bottle $L$ in the class $B\bmod 2$ using the standard toric picture. By \cite[\S3.1, 3.2]{AdaloglouEvans2026},
Klein bottle Luttinger surgery along $L$ turns
\((S^2\times S^2,\omega_{a,b})\) into \(\left(X_1,\omega_{a+b/2;a-b/2}\right),\) which follows from the homology correspondence
\begin{equation}\label{eq:AE-basic-homology}
    2A\longleftrightarrow H+E,
\qquad
B\longleftrightarrow H-E.
\end{equation}

Now one can perform symplectic blowups on both $S^2\times S^2$ and $X_1$ supported away from the surgery region by taking the same symplectic balls in their common complement $C$. After at least one such blowup, the manifolds on both sides have the same diffeomorphism type, since
\[
(S^2\times S^2)\#(n-1)\overline{\CC\PP}^{2}
\cong X_n,
\qquad n\geq2.
\]

\begin{proposition}
\label{prop:AE-reflection}
Suppose $(X_n,\omega)$ is obtained by a blowup of $(S^2\times S^2,\omega_{a,b})$, where $0<b<2a$, along $n-1$ pairwise disjoint symplectic balls avoiding a Lagrangian Klein bottle $L$ in the class $B\bmod2$.
Let $C$ be the exceptional class of one distinguished blowup, and let $C_3,\ldots,C_n$ be the remaining exceptional classes. Identify the homology by
\begin{equation}\label{eq:AE-source-marking}
H=A+B-C,
\qquad E_1=B-C,
\qquad E_2=A-C,
\qquad E_k=C_k\quad(k\geq3).
\end{equation}
Then Klein bottle Luttinger surgery along $L$ produces a symplectic manifold $(X_n,\omega')$ satisfying
\begin{equation}\label{eq:AE-reflection}
[\omega']=R_{v_{12}}([\omega]),
\qquad
R_{v_{12}}(x)=x+\frac{x\cdot v_{12}}{2}v_{12},
\qquad
v_{12}=H-2E_1-E_2.
\end{equation}
Equivalently, if $[\omega]=(h;e_1,e_2,\ldots,e_n)$, then the symplectic period change is given by
\begin{equation}\label{eq:AE-coordinate-formula}
(h;e_1,e_2,e_3,\ldots,e_n)
\longmapsto
\left(
\frac{3h-2e_1-e_2}{2};
 h-e_1-e_2,
\frac{h-2e_1+e_2}{2},
 e_3,\ldots,e_n
\right).
\end{equation}
\end{proposition}

\begin{proof}
Since the blowup balls are disjoint from $L$, Klein bottle Luttinger surgery commutes with the blowups. Thus $(X_n,\omega')$ can be obtained by first performing the Luttinger surgery on $(S^2\times S^2,\omega_{a,b})$, which produces \(\left(X_1,\omega_{a+b/2;a-b/2}\right),\) and then performing the same $n-1$ blowups.
Hence, we can choose a standard basis $H',E_1',\ldots,E_n'$ of $H_2(X_n;\ZZ)$ such that $H'$ and $E_2'$ are the total transforms of the line and exceptional classes in $X_1$, $E_1'$ is the distinguished blowup class, and $E_k'$ corresponds to $C_k$ for $k\geq3$. By \eqref{eq:AE-basic-homology}, the real homology isometry $\Phi$ in Lemma \ref{lem:intersection-form} satisfies
\[
\Phi(A)=\frac{H'+E_2'}2,\qquad
\Phi(B)=H'-E_2',\qquad
\Phi(C)=E_1',\qquad
\Phi(C_k)=E_k'\quad(k\geq3).
\]
Write $[\omega']=(h';e_1',\ldots,e_n')$, and let $c$ be the
symplectic area of $C$. The identification \eqref{eq:AE-source-marking} gives
\[a=h-e_1,\qquad b=h-e_2,\qquad c=h-e_1-e_2.\]
Since $\Phi$ preserves symplectic periods, we have
\[
h'=a+\frac b2,\qquad
e_1'=c,\qquad
e_2'=a-\frac b2,\qquad
e_k'=e_k\quad(k\geq3).
\]
Substituting $a,b,c$ by $h,e_1,e_2$ proves \eqref{eq:AE-coordinate-formula}. To see \eqref{eq:AE-reflection}, apply $\Phi$ to \eqref{eq:AE-source-marking}:
\begin{align*}
\Phi(H)&=\frac32H'-E_1'-\frac12E_2',\\
\Phi(E_1)&=H'-E_1'-E_2',\\
\Phi(E_2)&=\frac12H'-E_1'+\frac12E_2',\\
\Phi(E_k)&=E_k'\qquad(k\geq3).
\end{align*}
After identifying the markings $\{H,E_1,\ldots,E_n\}$ with $\{H',E_1',\ldots,E_n'\}$, these agree with the formula of the reflection along the $(-4)$-class $v_{12}$.
\end{proof}

\begin{definition}
We call the KB Luttinger surgery on the rational manifold $X_n$ which turns $\omega$ to $\omega'$ in Proposition~\ref{prop:AE-reflection} a {\bf geometric Adaloglou--Evans (AE) twist} and the real isometry $R_{v_{12}}:\Lambda\otimes\RR\to \Lambda\otimes \RR$ an {\bf algebraic AE twist}.
\end{definition}

\begin{figure}[htbp]
\centering
\begin{tikzpicture}[x=0.4cm,y=0.4cm,font=\footnotesize]
    \fill[gray!20] (2,0) -- (6,0) -- (6,4) -- cycle; \fill[gray!20] (0,9) -- (0,10) -- (1,10) -- cycle;

    \draw[dashed,gray] (2,0) -- (6,0) -- (6,4);
    \draw[dashed,gray] (0,9) -- (0,10) -- (1,10);
    \draw[thick] (0,0) -- (2,0) -- (6,4) -- (6,10) -- (1,10) -- (0,9) -- cycle;

    \draw[thick] (0.5,0) -- (5.5,10) node[midway,right] {$\ell$};
    \fill (0.5,0) circle (1.2pt);
    \fill (5.5,10) circle (1.2pt);

    \node at (4.7,1.3) {$4$};
    \node[font=\scriptsize,inner sep=0pt] at (0.33,9.67) {$1$};
    \node[below left]  at (0,0)  {$(0,0)$};
    \node[below right] at (6,0)  {$(6,0)$};
    \node[above left]  at (0,10) {$(0,10)$};
    \node[above right] at (6,10) {$(6,10)$};
\end{tikzpicture}
\caption{The moment rectangle with corner triangles of capacities $4$ and $1$ removed. The segment $\ell$, supporting the visible Lagrangian Klein bottle, avoids both triangles.}
\label{fig:AE-Karshon-Kessler}
\end{figure}
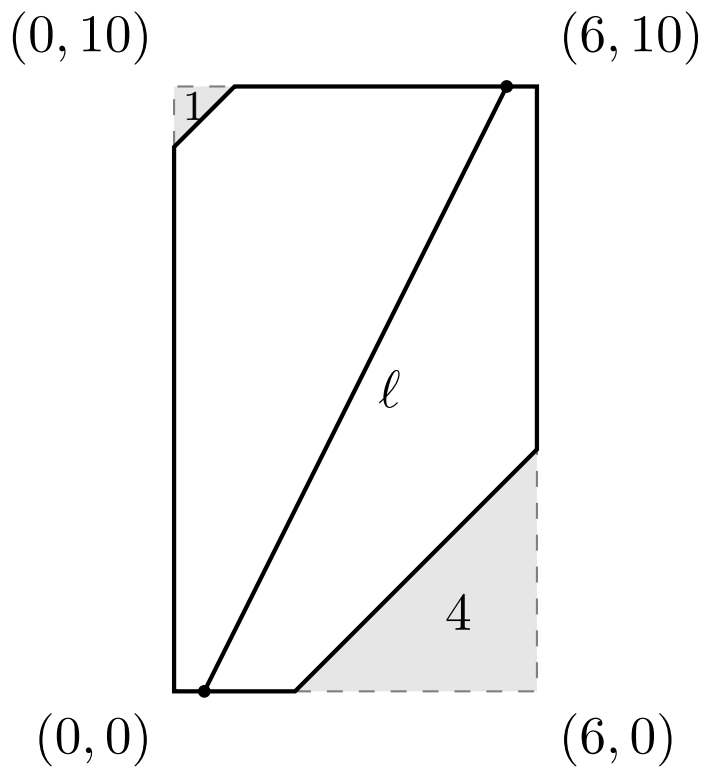

\begin{example}
\label{example:AE-Karshon-Kessler}
Consider two symplectic forms on $X_3$ with classes
\([\omega_0]=(12;6,2,1),\,[\omega_1]=(11;4,1,1).\) Their Auroux data agree:
\[\mathcal{A}(X_3,\omega_0)=\mathcal{A}(X_3,\omega_1)=(6,6,27,103).
\]
They occur in \cite[\S~1.1]{KK17} as an example of non-symplectomorphic forms on $X_3$. By \cite[Theorem~1.2]{HoLi2012}, no sequence of torus Luttinger surgeries connects these two forms, whereas one Klein bottle Luttinger surgery does.
First, they are related by a single algebraic AE twist:
\[
R_{v_{12}}([\omega_0])
=[\omega_0]-(H-2E_1-E_2)
=11H-4E_1-E_2-E_3
=[\omega_1].
\]
To see its geometric realizability, start with
$(S^2\times S^2,\omega_{6,10})$. Let $[0,6]\times[0,10]$ be its moment rectangle and take the visible Lagrangian Klein bottle over
\[
\ell=\left\{\left(\frac12+\frac t2,t\right):0\leq t\leq10\right\}
\subset[0,6]\times[0,10].
\]
Perform toric blowups of capacity $4$ at the lower-right corner and capacity $1$ at the upper-left corner. Both are disjoint from $\ell$; see Figure \ref{fig:AE-Karshon-Kessler}. Let $C$ be the exceptional class of capacity $4$ and $C_3$ the class of capacity $1$. By \eqref{eq:AE-source-marking}, its symplectic period is therefore $[\omega_0]=(12;6,2,1)$.
\end{example}
The $\ZZ_2$-homology class $[L]_2$ of the Lagrangian Klein bottle in Proposition \ref{prop:AE-reflection} is $H+E_2$, which is exactly the mod $2$ reduction of $v_{12}$. One can choose a different marking on $H_2(X_n;\ZZ)$ to change the class of $[L]_2$ so that the surgery effect is described by the reflection along a different class. More precisely, let $W_n\subset O(\Lambda)$ denote the {\bf Weyl group} generated by the reflections \(R_{\alpha_0},\ldots,R_{\alpha_{n-1}}\) (if $n=2$, only $R_{\alpha_1}$) where 
\[\alpha_0=H-E_1-E_2-E_3,\quad\alpha_i=E_i-E_{i+1},\,1\leq i\leq n-1,\qquad
R_{\alpha_i}(x)=x+(x\cdot\alpha_i)\alpha_i.
\]
 Changing the marking by $w\in W_n$ sends $[L]_2$ to $w([L]_2)$ and conjugates the reflection by \(wR_{12}w^{-1}=R_{w(v_{12})},\) where $w(v_{12})$ is another $(-4)$-class which pairs $0$ with $c_1$. For example, $w(v_{12})=v_{ij}=H-2E_i-E_j$ if $w$ swaps $E_1,E_2$ with $E_i,E_j$; when $n\geq4$, $w(v_{12})=E_3+E_4-E_1-E_2$ if $w$ denotes the reflection along $H-E_1-E_3-E_4.$

However, not every Lagrangian Klein bottle arises as the hypothesis of Proposition \ref{prop:AE-reflection}, even after a change of marking. For example, two Lagrangian blowups of a Lagrangian sphere in the class $H-E_1-E_2-E_3$ produce a characteristic Lagrangian Klein bottle $L\subset X_5$, with \([L]_2=H+E_1+\cdots+E_5 \bmod2\).
Such a Klein bottle has odd mod-two intersection with every exceptional class, whereas the construction in Proposition~\ref{prop:AE-reflection} leaves exceptional spheres disjoint from $L$.
To extend the reflection description beyond this situation, we establish a general integral $(-4)$-lifting lemma.

\begin{lemma}
\label{lem:KB-minus-four-lift}
For $n\ge2$, let $L\subset(X_n,\omega)$ be a Lagrangian Klein bottle in a symplectic rational four-manifold satisfying $c_1(X_n,\omega)\cdot[\omega]>0$.
Then there exists $v\in\Lambda$ satisfying
\[
v\bmod2=[L]_2,
\qquad
v^2=-4,
\qquad
c_1\cdot v=0.
\]
\end{lemma}

\begin{proof}
By Shevchishin--Nemirovski's theorem \cite{Shevchishin2009,Nemirovski2009}, $[L]_2\ne0$. Suppose first that $[L]_2\ne c_1\bmod2$, and suppose
\[
[L]_2=\epsilon H+\sum_{i\in I}E_i\bmod2,
\qquad
\epsilon\in\{0,1\},
\qquad
I\subset\{1,2,\ldots,n\}.
\]
We say such a $\ZZ_2$-class is of type $(\epsilon,|I|)$. Audin's congruence \cite{Audin} says $|I|\equiv\epsilon\pmod4$. If $\epsilon=0$, then $|I|\ge4$. Choosing distinct $i,j,k\in I$, the Cremona transformation $R_{H-E_i-E_j-E_k}$ turns $[L]_2$ into a class of type $(1,|I|-3)$. If $\epsilon=1$, then $n>|I|\ge5$. Choosing distinct $i,j\in I$ and $k\in \{1,2,\ldots,n\}\setminus I$, $R_{H-E_i-E_j-E_k}$ turns $[L]_2$ into a class of type $(0,|I|-1)$.
We may use these reflections to decrease $|I|$ until reaching a class of type $(1,1)$. Hence, after further performing a permutation, we may assume there is a $w\in W_n$ such that \([L]_2=w(H+E_2)\bmod2.\)
The integral $(-4)$-class \(v=w(H-2E_1-E_2)\) therefore satisfies our requirement.

Suppose now that $[L]_2=c_1\bmod2$. Audin's congruence gives $n\equiv1\pmod4$. We first show that $n\not\equiv1\pmod8$. Let $(X',\omega')$ be obtained by KB Luttinger surgery along $L$ and put $\Gamma=\Phi^{-1}(\Lambda')$, where $\Phi$ is the real isometry of Lemma~\ref{lem:intersection-form}. The sublattice induced by the surgery complement is then given by \(
\Lambda_C=\{x\in\Lambda:x\cdot c_1\equiv0\pmod2\}.\)
Since $c_1$ is characteristic, this is an even sublattice of index two in $\Lambda$. Since $\Lambda'\cong\Lambda$ is odd, we can choose $z\in\Gamma$ with odd square. By Lemma \ref{lem:pullback-intersection-lattices}, we have $\Gamma\cap \Lambda=\Lambda_C$ which implies that $z\notin\Lambda$. Since $[\Gamma:\Lambda_C]=[\Lambda:\Lambda_C]=2,$  $2z\in \Lambda_C\subset\Lambda.$ It follows that $2z\in\Lambda\setminus2\Lambda$.
For every $x\in\Lambda_C$, integrality of $\Gamma$ gives \((2z)\cdot x=2(z\cdot x)\equiv0\pmod{2}.\)
Since $\Lambda$ is unimodular and $2z\notin2\Lambda$, there exists $y\in\Lambda$ with $(2z)\cdot y\equiv1\pmod{2}$.
Necessarily $y\notin\Lambda_C$, so $c_1\cdot y\equiv1\pmod{2}$ by the definition of $\Lambda_C$.
As $[\Lambda:\Lambda_C]=2$, we have $\Lambda=\Lambda_C+\mathbb Z y$, and therefore \((2z-c_1)\cdot x\equiv0\pmod{2}\) for every $x\in\Lambda$. Unimodularity of $\Lambda$ now implies $2z-c_1\in2\Lambda$, so we can write $z=c_1/2+a$ for some $a\in\Lambda$. Since $a^2+a\cdot c_1$ is even,
\[
\frac{9-n}{4}
=\frac{c_1^2}{4}
=z^2-a^2-a\cdot c_1
\in2\ZZ+1,
\]
and hence $n\equiv5\pmod8$. Consequently, we may write $n=8q+5$ and set
\[
v=(2q+1)H-(2q-1)E_1
-\sum_{i=2}^{6q+5}E_i
+\sum_{i=6q+6}^{8q+5}E_i.
\]
One can easily verify $v\bmod2=c_1\bmod2=[L]_2$, $v^2=-4$ and $v\cdot c_1=0$. This completes the proof.
\end{proof}

\begin{remark}\label{rmk:Furuta-refinement}
The above proof excludes characteristic Lagrangian Klein bottles in a symplectic rational manifold $X_n$ with $n\equiv 1\pmod8$ and $c_1\cdot[\omega]>0.$ One can even strengthen the argument by using Furuta's $\frac{10}{8}$-theorem \cite{Furuta2001} to remove the assumption $c_1\cdot[\omega]>0$.
Indeed, if $n=8q+1$ with $q>0$, the preceding lattice argument shows that the surgery manifold $X'$ has even intersection form. Moreover, $H_1(C;\ZZ)\cong H_1(X;\ZZ)=0$ by the relative exact sequence. The surjectivity of $H_1(Y;\ZZ)\to H_1(U';\ZZ)$ then gives $H_1(X';\ZZ)=0$ by Mayer--Vietoris. Thus $X'$ is spin, and \cite{Furuta2001} implies
\[
n+1=b_2(X')
\ge \frac54|\sigma(X')|+2
=\frac54(n-1)+2,
\]
forcing $n\le1$, a contradiction. This also implies that a non-characteristic hypothesis is necessary for \cite[Proposition 3.1]{DaiHoLi2019}. 
\end{remark}
Now we define
\[
\Delta_{c_1}^{-4}
=
\{v\in\Lambda:v^2=-4,\ c_1\cdot v=0\},
\]
and for any $\ZZ_2$-class $\beta$, write
\[
\Delta_{c_1}^{-4}(\beta)
=
\{v\in\Delta_{c_1}^{-4}:v\bmod2=\beta\}.
\]
Lemma~\ref{lem:KB-minus-four-lift} therefore shows that \(\Delta_{c_1}^{-4}([L]_2)\ne\varnothing\) for every Lagrangian Klein bottle $L$ in a $c_1$-positive symplectic rational manifold. For $v\in\Delta_{c_1}^{-4}$, define
\[
R_v:\Lambda\otimes\RR\longrightarrow\Lambda\otimes\RR,
\qquad
R_v(x)=x+\frac{x\cdot v}{2}v.
\]
This is an involutive real isometry fixing $c_1$.

\begin{lemma}
\label{lem:KB-reflection-marking}
With the same hypotheses and notations of
Lemma~\ref{lem:KB-minus-four-lift}, for every $v\in\Delta_{c_1}^{-4}([L]_2)$, we have
\(
R_v(\Lambda_C)=\Lambda_C,
\,
R_v(\Gamma)=\Lambda.
\)
Consequently,
\(
\mathfrak m_v
:=R_v\circ\Phi^{-1}:\Lambda'\to\Lambda
\)
is an integral marking satisfying
\(
\mathfrak m_v\!\left(
\operatorname{PD}c_1(X',\omega')
\right)=c_1
\)
and
\(
(\mathfrak m_v\otimes\RR)
\bigl(\operatorname{PD}[\omega']\bigr)
=
R_v\bigl(\operatorname{PD}[\omega]\bigr).
\)
\end{lemma}

\begin{proof}
By Lemmas \ref{lem:intersection-form} and \ref{lem:pullback-intersection-lattices}, we have
\[
\Lambda_C=\{x\in\Lambda:x\cdot v\equiv0\pmod2\},
\qquad
[\Lambda:\Lambda_C]=[\Gamma:\Lambda_C]=2,
\qquad
\Lambda\cap\Gamma=\Lambda_C.
\]
For every $x\in\Lambda_C$, we have $x\cdot v\in2\ZZ$. Since $v\in\Lambda_C$ and $R_v^2=\mathrm{id}$, this gives \(R_v(\Lambda_C)=\Lambda_C.\)
Now choose $z\in\Gamma\setminus\Lambda_C$. Then $2z\in\Lambda_C$ and $z\notin\Lambda$, so $2z\in\Lambda\setminus2\Lambda$. We have $2z\cdot x=2(z\cdot x)\equiv 0\pmod2$ for any $x\in\Lambda_C$, but $2z\cdot y\equiv 1\pmod 2$ for some $y\in\Lambda\setminus\Lambda_C$ since $\Lambda$ is unimodular. It follows that $(2z-v)\cdot x\equiv 0\pmod 2$ for any $x\in\Lambda$ so that $2z-v\in2\Lambda$. Hence,
\(
\Gamma=\Lambda_C+\ZZ z,\) where
\(
z=a+v/2
\) for some $a\in\Lambda$. Using $v^2=-4$ and $c_1\cdot v=0$, we obtain
\[
0\equiv z^2-c_1\cdot z
=(a^2-c_1\cdot a)+a\cdot v-1
\equiv a\cdot v-1\pmod2.
\]
Thus $a\cdot v$ is odd so that \(\Lambda=\Lambda_C+\ZZ a.\)
We compute
\[
R_v(z)
=a+\frac{a\cdot v-1}{2}v
\equiv a\pmod{\Lambda_C}.
\]
Therefore,
\[
R_v(\Gamma)
=\Lambda_C+\ZZ R_v(z)
=\Lambda_C+\ZZ a
=\Lambda.
\]
It follows that $\mathfrak m_v=R_v\circ\Phi^{-1}$ is an
isometry from $\Lambda'$ to $\Lambda$. The identities for the first Chern class and the symplectic class follow from Lemma \ref{lem:intersection-form} and $R_v(c_1)=c_1$.
\end{proof}

\begin{remark}\label{rmk:geometric-marking}

If $v$ and $\widetilde v$ are two classes in $\Delta^{(-4)}_{c_1}([L]_2)$, then \(R_v(\Lambda)=R_{\widetilde v}(\Lambda)=\Gamma.\)
Thus $R_{\widetilde v}R_v$ is an automorphism of $\Lambda$ fixing $c_1$, and the two periods $R_v\bigl(\operatorname{PD}[\omega]\bigr)$ and $R_{\widetilde{v}}\bigl(\operatorname{PD}[\omega]\bigr)$ differ by this change of marking. However, Lemma \ref{lem:KB-reflection-marking} does not assert that a particular marking $\mathfrak m_v:\Lambda'\to\Lambda$ is induced by a diffeomorphism from $X'_n$ to $X_n$, when $n>9$. Indeed, by \cite[Theorem 2.8]{LiLi2002}, this is equivalent to requiring the period $R_v\bigl(\operatorname{PD}[\omega]\bigr)$ in the $c_1$-symplectic cone.
\end{remark}

\subsection{The surgery graphs}
\label{subsec:rational-KB-graphs}

We now associate two graphs to the Klein bottle Luttinger geography problem on rational manifolds. The algebraic graph records the lattice transformations permitted by Lemma \ref{lem:KB-reflection-marking}, while the geometric graph records reflection moves admitting Lagrangian Klein bottle representatives. Their connected components give lower and upper bounds, respectively, for the number of Klein bottle Luttinger equivalence classes.

Consider the positive forward cone
\[
\mathscr P_n^+
=
\{x\in\Lambda\otimes{\RR}:x^2>0,\ x\cdot H>0\},
\]
and denote by $\mathcal C_{c_1}$ the $c_1$-symplectic cone in
$\mathscr P_n^+$:
\[
\mathcal C_{c_1}
=
\left\{
\operatorname{PD}[\omega]\in\mathscr P_n^+:
\omega\text{ is symplectic with first Chern class }c_1
\right\}.
\]
We distinguish the $c_1$-preserving lattice automorphism group
\[
O^+(\Lambda,c_1)
=
\left\{
A\in O(\Lambda):
A(c_1)=c_1,\ A(\mathscr P_n^+)=\mathscr P_n^+
\right\}
\]
from its subgroup induced by diffeomorphisms:
\[
D_{c_1}
=
\left\{
A\in O^+(\Lambda,c_1):
A=f_*\text{ for some }
f\in\operatorname{Diff}^+(X_n)
\right\}.
\]
Every element in $D_{c_1}$ preserves $\mathcal C_{c_1}$ by definition. For every $n$, $D_{c_1}$ is equal to the Weyl group $W_n$; see \cite[Propositions~4.7 and~4.10]{LiWu2012}. For $n\le9$, it follows from Wall's theorem \cite{Wall1964} that $O^+(\Lambda,c_1)=D_{c_1}=W_n$.
For $n=10$, we can write \(\Lambda=\mathbb Zc_1\oplus c_1^\perp,\) where $c_1^\perp$ is isomorphic to the Enriques lattice, whose isometry group is known to be generated by Weyl reflections and $\{\pm\operatorname{Id}\}$; see \cite[\S 3]{Dolgachev2016}. Hence we still have $O^+(\Lambda,c_1)=D_{c_1}=W_{10}$. When $n>10$, $O^+(\Lambda,c_1)$ is strictly larger than $D_{c_1}$ and there exist elements in $O^+(\Lambda,c_1)\setminus D_{c_1}$ which do not preserve $\mathcal{C}_{c_1}$. 

Let $\ZZ[1/2]$ be the ring of dyadic rational numbers and put \(\Lambda[1/2]=\Lambda\otimes_{\ZZ}\ZZ[1/2].\) For any pair of positive integers $(P,Q)$ such that the extended Auroux data $(9-n,3+n,P,Q,0)$ is realizable, define
\[
\mathbb{V}_{n;P,Q}
=
\left\{
x\in\mathcal C_{c_1}\cap\Lambda[1/2]:
c_1\cdot x=P,\ x^2=Q
\right\},
\qquad
\mathbb{I}_{n;P,Q}=\mathbb{V}_{n;P,Q}\cap\Lambda.
\]

\begin{definition}[KB Luttinger surgery graphs for rational manifolds]
\label{def:rational-KB-graphs}
The graphs $\mathbb G_{n;P,Q}^{\mathrm{alg}}$ and $\mathbb G_{n;P,Q}^{\mathrm{geo}}$ are simple graphs with common vertex set $\mathbb V_{n;P,Q}$. Both contain the {\bf Weyl edges} \(z\leftrightarrow w(z)\) for any vertex $z$ and $w\in W_n$ such that $z\neq w(z)$, and the {\bf Luttinger edges} between distinct vertices $x$ and $y$ which are defined as follows.
\begin{enumerate}[nosep]
\item
The {\bf algebraic graph} $\mathbb G_{n;P,Q}^{\mathrm{alg}}$ contains an edge  whenever \(y=A R_v(x)\) for some $A\in O^+(\Lambda,c_1)$ and $v\in\Delta_{c_1}^{(-4)}$.

\item
The {\bf geometric graph} $\mathbb G_{n;P,Q}^{\mathrm{geo}}$ contains an edge whenever \(y= R_v(x)\) for some $v\in\Delta_{c_1}^{(-4)}$ and there exists a Lagrangian Klein bottle $L\subset(X_n,\omega)$ satisfying $\operatorname{PD}[\omega]=x$ and \( [L]_2=v\bmod2.\)
\end{enumerate}
 Each unordered pair $\{x,y\}$ admits at most one edge, even if it satisfies both edge conditions or admits several presentations.
\end{definition}

It is straightforward to see that the geometric graph is a subgraph of the algebraic graph. For an algebraic Luttinger edge, only $x$ and $A R_v(x)$ must belong to $\mathbb V_{n;P,Q}$; the intermediate class $R_v(x)$ need not belong to the $c_1$-symplectic cone $\mathcal C_{c_1}$. In contrast, a geometric Luttinger edge requires $R_v(x)\in\mathcal C_{c_1}$ as well as a Lagrangian Klein bottle with the prescribed mod-$2$ class. Notice that any symplectic form representing a class in $\mathcal C_{c_1}$ has first Chern class $c_1$, since symplectic cones associated to distinct canonical classes are disjoint; see \cite[\S 2]{LiLi2002}. In particular, the symplectic form in Definition \ref{def:rational-KB-graphs} has the distinguished element $c_1\in\Lambda$ as its first Chern class. 

For $\star\in\{\mathrm{alg},\mathrm{geo}\}$, define
\[N_n^\star(P,Q)=\#\bigl\{T\in\pi_0(\mathbb G_{n;P,Q}^\star):T\cap\mathbb I_{n;P,Q}\ne\varnothing\bigr\}.\]
Here $T$ denotes a connected component of the full graph. These are exactly the counts of all connected components in the graphs $\mathbb X(n,P,Q)$ and $\mathbb Y(n,P,Q)$ introduced in Section \ref{sec:intro-KB-Luttinger}.

\begin{proposition}
\label{prop:rational-KB-graph-bounds}
For every $n\ge2$ and realizable $P,Q\in\mathbb Z_{>0}$,
\[N_n^{\mathrm{alg}}(P,Q)\le \mathcal N_{(9-n,\,3+n,\,P,\,Q,\,0)}\le N_n^{\mathrm{geo}}(P,Q)<\infty.\]
\end{proposition}

\begin{proof}
Since $P>0$, every symplectic four-manifold with the prescribed extended Auroux data must be diffeomorphic to $X_n$ by Liu--Ohta--Ono's theorem \cite{Liu1996,OhtaOno1996}. By Lemma \ref{lemma:surgery-betti-change} and Corollary \ref{cor:surgery-symplectic-pairings},  every intermediate manifold in a surgery sequence is also a symplectic rational manifold.
Every such a manifold admits a marking induced by a diffeomorphism to $X_n$ that identifies its first Chern class with the distinguished element $c_1\in \Lambda$ and its symplectic period with a class in $\mathcal C_{c_1}$; see \cite[\S 4.1.1]{LiWu2012}. Two such markings differ by an element of $D_{c_1}=W_n$. Moreover, symplectic forms on $X_n$ with the same period class are symplectomorphic through a diffeomorphism acting trivially on homology; see \cite{KK17}.
Thus a vertex determines a symplectic manifold up to symplectomorphism, every integral object being counted is represented by a vertex of $\mathbb I_{n;P,Q}$, and changing a marking does not change its component in either graph because of the Weyl edges.

For the lower bound, consider an actual KB Luttinger surgery that turns $(X_n,\omega)$ to $(X',\omega')$, where $\operatorname{PD}[\omega]=x\in\mathcal{C}_{c_1}.$ By Lemmas \ref{lem:KB-minus-four-lift} and \ref{lem:KB-reflection-marking}, there is a class $v\in\Delta_{c_1}^{(-4)}([L]_2)$ and a marking \(\mathfrak m_v:H_2(X';\mathbb Z)\rightarrow\Lambda\) such that \(\mathfrak m_v\bigl(\operatorname{PD}c_1(X',\omega')\bigr)=c_1,\,\mathfrak m_v\bigl(\operatorname{PD}[\omega']\bigr)=R_v(x).\) Although $\mathfrak m_v$ might not be induced by a diffeomorphism, we can choose another diffeomorphism-induced marking $\mathfrak n:H_2(X';\mathbb Z)\rightarrow\Lambda$ taking the first Chern class to $c_1$ and the period into $\mathcal C_{c_1}$, and put \(y=\mathfrak n\bigl(\operatorname{PD}[\omega']\bigr),\, A=\mathfrak n\circ \mathfrak m_v^{-1}.\) Then $A\in O^+(\Lambda,c_1)$ and \(y=A R_v(x).\) Furthermore, since $R_v$ and $A$ preserve $\Lambda[1/2]$, in a sequence of surgeries with endpoints equipped with integral symplectic forms, the intermediate manifolds all have dyadic periods. Thus all manifolds in an actual KB Luttinger equivalence class are contained in a single component of $\mathbb G^{\operatorname{alg}}_{n;P,Q}$. On the other hand, every component meeting $\mathbb I_{n;P,Q}$ must contain at least one KB Luttinger equivalence class. This proves the lower bound.

For the upper bound, it suffices to show that every geometric edge gives an actual KB Luttinger equivalence. This holds for Weyl edges because $W_n=D_{c_1}$. For a geometric Luttinger edge $y=R_v(x)$, suppose $L\subset(X_n,\omega)$ satisfies $[L]_2=v\bmod 2$. Perform surgery along $L$ to obtain $(X',\omega')$ and set $y=\mathfrak m_v\bigl(\operatorname{PD}[\omega']\bigr)$, where $\mathfrak m_v$ is the marking given by Lemma \ref{lem:KB-reflection-marking}. By Definition \ref{def:rational-KB-graphs}, $y\in\mathcal C_{c_1}$ and $c_1\cdot y=P>0$. By \cite[Theorem~2.8]{LiLi2002}, the marking $\mathfrak m_v$ is induced by a diffeomorphism between $X'$ and $X_n$. $(X',\omega')$ consequently represents $y$ up to symplectomorphism. The components of $\mathbb G_{n;P,Q}^{\mathrm{geo}}$ meeting $\mathbb I_{n;P,Q}$ thus map onto the actual equivalence classes counted by $\mathcal N_{(9-n,3+n,P,Q,0)}$, proving the upper bound.

Finally, let $T$ be a component counted by $N_n^{\mathrm{geo}}(P,Q)$. The integral vertex of $T$ can be connected to a reduced integral vertex $x=hH-\sum_{i=1}^n e_iE_i$ in the sense of \cite{KK17} via Weyl edges. This means  $e_1\ge\cdots\ge e_n\ge1$, $h\ge e_1+e_2+e_3$ for $n\geq 3$ and $h>e_1+e_2$ for $n=2$. Set $m=e_3$ for $n\geq 3$ and $m=0$ for $n=2$. Since $e_i\le m$ for $i\ge3$, the reduced condition gives
\[\sum_{i=1}^n e_i(e_i-m)\le e_1(e_1+e_2-2m)\le e_1(h-3m).\]
Using $P=3h-\sum_i e_i$ and $Q=h^2-\sum_i e_i^2$, we obtain
\[Q-mP=h(h-3m)-\sum_i e_i(e_i-m)\ge (h-e_1)(h-3m).\]
Since $P\ge1$, $h-e_1\ge1$, and $h\ge3m$, it follows that \[Q\ge m+(h-3m)=h-2m\ge h/3.\] Thus there are only finitely many possibilities for the integer $h$ when $Q$ is fixed. It follows that the choices of $e_1,\ldots,e_n$ are also finitely many. Consequently $N_n^{\mathrm{geo}}(P,Q)<\infty$.
\end{proof}

The bound $N_n^\star(P,Q)$ can be obtained by some alternative versions of the graphs. For instance, one may consider the quotient graph by using the vertex set $\mathbb V_{n;P,Q}/W_n$, with edges induced from the original graph. The quotient induces a bijection on connected components since all vertices in the orbit $W_n\cdot v$ are already connected by Weyl edges in the original graph. The same counts can therefore be obtained by counting components meeting the image of $\mathbb I_{n;P,Q}$.

When $n\not\equiv5\pmod8$, the algebraic graph has the same connected components as its subgraph consisting of Weyl edges and only Luttinger edges $x\leftrightarrow AR_{v_{12}}(x)$, with $A\in O^+(\Lambda,c_1)$ and the algebraic AE twist $R_{v_{12}}$. Indeed, for a Luttinger edge $y=AR_v(x)$, the proof of Lemma \ref{lem:KB-minus-four-lift} gives $w\in W_n$ such that $u=w(v_{12})\equiv v\pmod2$. The  isometry $B=R_vR_u$ belongs to $O^+(\Lambda,c_1)$ so that \(y=(ABw)R_{12}(w^{-1}x).\) Thus the edge can be replaced by the path $x\leftrightarrow w^{-1}x\leftrightarrow y$ in the subgraph.

When $2\leq n\leq 10$, the algebraic graph can be also replaced, without changing its connected components, by the subgraph consisting of the Weyl edges and only  edges $s\leftrightarrow R_v(x)$. In this range, \(  O^+(\Lambda,c_1)=D_{c_1}=W_n.\) If $y=A R_v(x)$ is an arbitrary algebraic Luttinger edge, then $A\in W_n$ and \(R_v(x)=A^{-1}y\in\mathcal C_{c_1}.\) The edge can therefore be replaced by the path \(x\leftrightarrow R_v(x)\leftrightarrow y\) in the subgraph. For the geometric graph, if $x\in\mathbb V_{n;P,Q}$ and $v\in\Delta_{c_1}^{(-4)}$ lifts the $\ZZ_2$-class of an actual Lagrangian Klein bottle for a symplectic form representing $x$, then $R_v(x)\in\mathbb V_{n;P,Q}$ automatically. Thus, we have the following fact.
\begin{corollary}
\label{cor:rational-KB-at-most-ten}
When $2\le n\le10$, \(\mathcal N_{(9-n,\,3+n,\,P,\,Q,\,0)}=N_n^{\mathrm{geo}}(P,Q).\)
\end{corollary}

\subsection{The count in small Euler numbers}
\label{subsec:two-three-four-point-blowups}
We retain the notations from Sections \ref{subsec:AE-twist}, \ref{subsec:rational-KB-graphs} and focus on the case when $2\leq n\leq4$. In this range, every algebraic Luttinger edge between vertices in the $c_1$-symplectic cone is geometrically realizable. It follows that the algebraic graph  $\mathbb G_{n;P,Q}^{\mathrm{alg}}$ also determines the actual count of the KB Luttinger equivalence classes by Corollary \ref{cor:rational-KB-at-most-ten}.

\begin{proposition}
\label{prop:small-blowup-realization}
Let $2\leq n\leq4$ and $v\in\Delta_{c_1}^{(-4)}$. Suppose that \(x,R_vx\in\mathcal C_{c_1}\) and \(R_vx\neq x.\) Then there is a symplectic form $\omega$ with $\operatorname{PD}[\omega]=x$ and a Lagrangian Klein bottle $L\subset(X_n,\omega)$ satisfying \([L]_2=v\bmod2.\)
\end{proposition}

\begin{proof}
Writing $v=dH-\sum_{i=1}^n m_iE_i$, the identities $\sum_i m_i=3d$ and $\sum_i m_i^2=d^2+4$ imply $(9-n)d^2\le4n$ by Cauchy--Schwarz, hence $|d|\le1$ for $2\le n\le4$, so, up to sign and permutation, $v$ equals $v_{12}=H-2E_1-E_2$ or, when $n=4$, $E_3+E_4-E_1-E_2=R_{H-E_1-E_3-E_4}(v_{12})$, and therefore $v=\pm w(v_{12})$ for some $w\in W_n$. Since the Weyl group $W_n$ preserves the $c_1$-symplectic cone $\mathcal{C}_{c_1}$, we may reduce to the case when $v=v_{12}$. Let us write $R$ for $R_{v_{12}}$ below. Put $r=n-1$ and use the identification \[
A=H-E_1,\qquad B=H-E_2,\qquad
C_1=H-E_1-E_2,\qquad
C_j=E_{j+1}\quad(2\leq j\leq r).
\]
Set \(a=x\cdot A,\,b=x\cdot B,\, p_i=x\cdot C_i.\)
Since $C_i$, $A-C_i$ and $B-C_i$ are exceptional classes, 
\begin{equation}
\label{eq:small-blowup-source}
0<p_i<\min\{a,b\}.
\end{equation}
Since $v_{12}=2C_1-B$, the reflection $R$ acts by
\[
R(B)=B,\qquad
R(C_1)=B-C_1,\qquad
R(C_j)=C_j\quad(j>1),\qquad
R(A)=A-C_1+\frac12B.
\]
Since $Rx\in\mathcal{C}_{c_1}$, evaluating $Rx$ on the exceptional class $E_2=A-C_1$ gives
\begin{equation}
\label{eq:small-blowup-target}
a-\frac b2>0.
\end{equation}
Let $G\subset W_n$ be the subgroup generated by reflections along $(-2)$-classes $C_i-C_j$ and $B-C_i-C_j$ for all distinct $i,j$. Note that elements in $G$ fix $[L]_2=B\bmod 2.$ We can use $R_{C_i-C_j}$ to permute the $p_i$'s and $R_{B-C_i-C_j}$ to change the periods by \((a,p_i,p_j)\mapsto(a+b-p_i-p_j, b-p_j,b-p_i),\) leaving $b$ and the remaining $p_k$ unchanged. We can eventually find some $g\in G$ such that the periods of $gx$ satisfy
\begin{equation}
\label{eq:small-blowup-normal-form}
p_1\geq p_2\geq\cdots\geq p_r,\qquad p_1+p_2\leq b.
\end{equation}
This is because $W_n$ is finite, we may minimize $a$ over the $G$-orbit and then permute the $C_i$. Moreover, $R$ exchanges the $(-2)$-classes $C_1-C_j$ and $B-C_1-C_j$ and fixes the other generating $(-2)$-classes of $G$. Hence, we have $RgR\in G$. This implies that \(Rgx=(RgR)Rx\in\mathcal{C}_{c_1}\). We will treat the class $gx$ satisfying \eqref{eq:small-blowup-normal-form} in the following discussion.

Consider the moment rectangle $[0,a]\times[0,b]$. By \eqref{eq:small-blowup-source} and \eqref{eq:small-blowup-target}, we can choose
\[0<s<\min\left\{a-\frac b2,\ a-p_1\right\}.\]
The segment
\[\ell=\{(s+\frac{t}{2},t):0\leq t\leq b\}\]
supports the visible Lagrangian Klein bottle in the class $B\bmod2$ as shown in Example \ref{example:AE-Karshon-Kessler}. Place a ball of capacity $p_1$ at the lower-right corner. Its moment triangle is disjoint from $\ell$, because on $\ell$
\[a-(s+\frac{t}{2})+t=a-s+\frac t2>p_1.\]
For $n=2$, this is the only required ball. If $n\geq3$, place a second ball of capacity $p_2$ at the upper-left corner. The two corner triangles are obviously disjoint. The second triangle also avoids $\ell$, since
\[(s+\frac{t}{2})+b-t=s+b-\frac t2\geq s+\frac b2>p_2.\]
Here $p_2\leq b/2$ follows from \eqref{eq:small-blowup-normal-form}. For $n=3$, these are all the required two balls. It remains to place the third ball when $n=4$. If the class $gx$ has $p_1=p_2=p_3=b/2$, it is fixed by both $G$ and $R$ so that $gRx=(gRg)gx=gx$, contrary to $Rx\neq x$. Therefore \(p_3<b/2,\, p_2+p_3<b,\) and we can choose
\begin{equation}
\label{eq:small-blowup-third-ball}
p_3<q<\min\{a,\frac{b}{2},b-p_2\}.
\end{equation}
Consider the triangle
\[T_q=\operatorname{conv}\{(0,0),(0,q),(q,2q)\}.\]
Its primitive edge directions are $(1,2)$ and $(0,1)$ so that $T_q$ is integral affine equivalent to the standard moment triangle of capacity $q$. Also, it is straightforward to verify that $T_q$ avoids both corner triangles and lies strictly to the left of $\ell$ by using \eqref{eq:small-blowup-third-ball}. One can then use either Traynor's trick \cite{Traynor} to pack a ball of capacity $q$ into the preimage of a small translation of $T_q$; or view $T_q$ as a Symington triangle to perform an almost toric blowup \cite{Symington}.

Therefore, for the symplectic class $gx$, we can always find a Lagrangian Klein bottle within the class $B\bmod2$. Since $g(B)=B$, such a Lagrangian Klein bottle also exists for the symplectic class $x.$

\end{proof}

\begin{remark}
\label{rmk:Rv-fixed-points}
For $n=2,3$, the preceding proof also works when $R_vx=x$. For $n=4$, the condition $R_vx\neq x$ does not affect the counting using algebraic graph, since an algebraic Luttinger edge
$x\leftrightarrow AR_vx$ is already a Weyl edge when $R_vx=x$.
\end{remark}

We write $R$ for the algebraic AE twist $R_{v_{12}}$ and consider a simplified graph for counting.

\begin{definition}
\label{def:small-blowup-AE-graph}
For $2\leq n\leq4$, let $\mathbb G'_{n;P,Q}$ be the simple graph on $\mathbb V_{n;P,Q}$ consisting of all the Weyl edges $x\leftrightarrow w(x)$, $w\in W_n$, and the edges $x\leftrightarrow Rx$, whenever the endpoints are distinct and belong to $\mathbb V_{n;P,Q}$.
\end{definition}

\begin{proposition}
\label{prop:small-blowup-count}
For $2\leq n\leq4$ and $P,Q\in\mathbb Z_{>0}$,
\[
N_n^{\mathrm{alg}}(P,Q)
=\mathcal N_{(9-n,\,3+n,\,P,\,Q,\,0)}
=N_n^{\mathrm{geo}}(P,Q)
=\#\bigl\{T\in\pi_0(\mathbb G'_{n;P,Q}):
T\cap\mathbb I_{n;P,Q}\neq\varnothing\bigr\}.
\]
\end{proposition}

\begin{proof}
The reductions at the end of Section \ref{subsec:rational-KB-graphs} show that
$\mathbb G'_{n;P,Q}$ and $\mathbb G^{\mathrm{alg}}_{n;P,Q}$ have the same connected components since $O^+(\Lambda,c_1)=W_n$, and every $(-4)$-reflection is Weyl conjugate to $R$. By Proposition~\ref{prop:small-blowup-realization} and Definition \ref{def:rational-KB-graphs}, we have
\[\mathbb G'_{n;P,Q}\subset\mathbb G^{\mathrm{geo}}_{n;P,Q}\subset\mathbb G^{\mathrm{alg}}_{n;P,Q}.\]
All three graphs therefore have the same connected components and the claim follows from Corollary \ref{cor:rational-KB-at-most-ten}.
\end{proof}

We write $W_2=\{\mathrm{id},S\}$ where $S$ is the permutation on the exceptional areas $e_1$ and $e_2$.

\begin{proposition}
\label{prop:two-point-paths}
The graph $\mathbb G'_{2;P,Q}$ has no cycles and every vertex has degree at most two. Consequently, each connected component is a (possibly infinite) chain.
\end{proposition}

\begin{proof}
Each vertex has at most one $S$-edge and one $R$-edge since $S^2=R^2=\operatorname{id}$. A cycle would therefore alternate between the two types and give \((RS)^m x=x\) for some $m\geq1$. Put $T=RS$ and write $T^jx=(h_j;e_{1,j},e_{2,j})$. The reflection formula \eqref{eq:AE-coordinate-formula} gives
\[2h_{j+1}=-3h_j+2P+e_{1,j},\qquad e_{1,j+1}=P-2h_j.\]
Eliminating $e_{1,j+1}$ yields
\[2h_{j+2}+3h_{j+1}+2h_j=3P.\]
Subtracting consecutive identities gives
\[2(h_{j+3}-h_{j+2})+3(h_{j+2}-h_{j+1})+2(h_{j+1}-h_j)=0.\]
If $T^mx=x$, the sequence $h_j$ is periodic. If it were nonconstant, we could rescale all $h_{\star+1}-h_\star$ to integers by multiplying a common factor, such that at least one, say $h_{j+2}-h_{j+1}$, becomes an odd integer. This contradicts the above identity modulo two. Thus $h_j$ is constant so that $7h_j=3P$ and $e_{1,j}=e_{2,j}=P/7$. However, this implies that $x=(P/7)c_1$, which is fixed by both $S$ and $R$ and cannot lie on a cycle.
\end{proof}

Although the vertex set $\mathbb V_{n;P,Q}$ allows non-integral points, when $n=2$, the components meeting $\mathbb I_{2;P,Q}$ can be computed entirely within the integral subgraph. We now express the $n=2$ KB Luttinger count as a difference of two finite counts of solutions of Diophantine equations subject to certain restrictions. Define
\[
\mathcal T^{\mathrm W}_{2;P,Q}
=
\left\{
(h;e_1,e_2)\in\mathbb Z_{>0}^3
\;\middle|\;
\begin{aligned}
&3h-e_1-e_2=P,\quad h^2-e_1^2-e_2^2=Q,\\
&h>e_1+e_2,\quad e_1\geq e_2
\end{aligned}
\right\},
\]
and
\[
\mathcal T^{R}_{2;P,Q}
=
\left\{
(h;e_1,e_2)\in\mathbb Z_{>0}^3
\;\middle|\;
\begin{aligned}
&3h-e_1-e_2=P,\quad h^2-e_1^2-e_2^2=Q,\\
&h>2e_1+e_2,\quad h\equiv e_2\pmod2
\end{aligned}
\right\}.
\]

\begin{theorem}
\label{thm:two-point-arithmetic-count}
For $P,Q\in\mathbb Z_{>0}$,
\[\mathcal N_{(7,5,P,Q,0)}=\#\mathcal T^{\mathrm W}_{2;P,Q}-\#\mathcal T^{R}_{2;P,Q}.\]
\end{theorem}

\begin{proof}
We first show that integral vertices in the same component of $\mathbb G'_{2;P,Q}$ are connected by a path of integral vertices. Take a shortest path, whose edges alternate between $S$ and $R$. If an $R$-edge first leaves the integral vertices, then \eqref{eq:AE-coordinate-formula} gives
\[h,e_2\in\mathbb Z+\frac12,\qquad e_1\in\mathbb Z.\]
More generally, suppose the current vertex has the form
\[(h;e_1,e_2)=\left(\frac{a}{2^m};\frac{b}{2^{m-1}},\frac{c}{2^m}\right),\qquad a,c\text{ odd},\quad b\in\mathbb Z.\]
Applying $S$ followed by $R$ gives
\[RS(h;e_1,e_2)=\left(\frac{3a-2b-2c}{2^{m+1}};\frac{a-2b-c}{2^m},\frac{a+2b-2c}{2^{m+1}}\right).\]
The numerators of the first and third coordinates are again odd, so the same pattern holds with $m$ replaced by $m+1$. So, every subsequent vertex remains non-integral, contradicting the integrality of the endpoint.

Consequently, the components of the subgraph containing only integral vertices in $\mathbb I_{2;P,Q}$ one-to-one correspond to the components counted in Proposition~\ref{prop:small-blowup-count}. This subgraph is acyclic by Proposition~\ref{prop:two-point-paths}, so its number of components is its number of vertices minus its number of edges. Subtracting the number of $S$-edges from the number of vertices gives the number of Weyl orbits. Each orbit has a unique representative with $e_1\geq e_2$, so this contribution is $\#\mathcal T^{\mathrm W}_{2;P,Q}$. It remains to count the $R$-edges. For an integral vertex, $Rx$ is integral precisely when $h\equiv e_2\pmod2$. The reflection formula \eqref{eq:AE-coordinate-formula} gives
\[h(Rx)-h(x)=\frac{h-2e_1-e_2}{2}.\]
Hence, exactly one of the endpoints of an $R$-edge with smaller $h$ satisfies $h>2e_1+e_2$. Conversely, the inequality $h>2e_1+e_2$ implies $Rx\in\mathcal C_{c_1}$ and $Rx\neq x$. This shows that $\mathcal T^{R}_{2;P,Q}$ is one-to-one correspondence to all $R$-edges with integral endpoints. Subtracting the count of these $R$-edges gives the counting formula.
\end{proof}

Using Theorem \ref{thm:two-point-arithmetic-count}, we can construct both families of extended Auroux data whose KB Luttinger count is one, or greater than one. For $m\geq 2$, we have
\[\mathcal N_{(7,\,5,\,7m,\,7m^2-2,\,0)}=1;\]
and for $m\geq 5$, we have
\[\mathcal N_{(7,\,5,\,7m,\,7m^2-14,\,0)}=2.\]
Indeed, the defining equations of $P,Q$ give $e_1+e_2=3h-P$ and $e_1^2+e_2^2=h^2-Q$. When $P=7m$,
\[(e_1-e_2)^2=2(h^2-Q)-(3h-7m)^2=2(7m^2-Q)-7(h-3m)^2,\]
which rearranges to
\[7(h-3m)^2+(e_1-e_2)^2=2(7m^2-Q).\]
For $Q=7m^2-2$, the right hand side is $4$, so the unique solution with $e_1\geq e_2$ is
\[(h;e_1,e_2)=(3m;m+1,m-1).\]
It lies in $\mathcal C_{c_1}$ for $m\geq2$. The exceptional areas are in the opposite parity with $h$. Hence $\mathcal T^R_{2;P,Q}$ is empty so that the count is one. For $Q=7m^2-14$, the right-hand side is $28$, so the only solutions are
\[(3m-2;m-3,m-3),\qquad (3m+2;m+3,m+3).\]
Both lie in $\mathcal C_{c_1}$ for $m\geq5$. Again, $h-e_2$ is odd so that $\mathcal T^R_{2;P,Q}$ is empty. The counting formula therefore gives two.

\section{KB Luttinger surgery on irrational ruled manifolds}
\label{sec:irrational-ruled}

We now fix $g\geq1$ and consider irrational ruled manifolds of genus $g$. Their minimal models have exactly two diffeomorphism types. In this section, we write $X_0=S^2\times\Sigma_g$ for the trivial $S^2$-bundle and $X_1=S^2\tilde{\times}\Sigma_g$ for the twisted $S^2$-bundle. Choose a fiber class $F$ and a section class $S_\varepsilon$ with
\[F^2=0,\qquad S_\varepsilon\cdot F=1,\qquad S_\varepsilon^2=-\varepsilon,\qquad \varepsilon\in\{0,1\}.\]
On the $k$-fold blowup, write $E_1,\ldots,E_k$ for exceptional classes. We have
\[\operatorname{PD}(c_1)=2S_\varepsilon+(2-2g+\varepsilon)F-\sum_{i=1}^k E_i,\qquad \operatorname{PD}[\omega]=aS_\varepsilon+\left(h+\frac{\varepsilon a}{2}\right)F-\sum_{i=1}^k\delta_iE_i.\]
Thus $a=[\omega]\cdot F$, $\delta_i=[\omega]\cdot E_i$, and
\[P=c_1\cdot [\omega]=2h-2(g-1)a-\sum_{i=1}^k\delta_i,\qquad Q=[\omega]^2=2ah-\sum_{i=1}^k\delta_i^2.\]
Our goal is to investigate the count $\mathcal{N}_{\widetilde{\mathcal A}}$ for the extended Auroux data 
\begin{equation}
\label{eq:ruled-data}
\widetilde{\mathcal{A}}
=(8-8g-k,\ 4-4g+k,\ P,\ Q,\ 2g).
\end{equation}
\subsection{Existence of Lagrangian KB in minimal models}

In $S^2\times S^2$, the existence of a Lagrangian Klein bottle in a prescribed $\ZZ_2$-class imposes a cohomological restriction on the symplectic form \cite{AdaloglouEvans2026}. However, such a  restriction disappears when the base genus becomes positive.

\begin{lemma}
\label{lem:ruled-construction}
Every minimal irrational ruled symplectic four-manifold contains a Lagrangian Klein bottle in the mod-two fiber class.
\end{lemma}

\begin{proof}
 Let $g\geq1$ and $\varepsilon\in\{0,1\}$. Every symplectic class $aS_\varepsilon+\left(h+\frac{\varepsilon a}{2}\right)F$ on $X_\varepsilon$ has $a,h>0$, and the classification of Lalonde--McDuff~\cite{LalondeMcDuff1996Ruled} identifies any symplectic form in the fixed cohomology class up to symplectomorphism. Therefore, it suffices to construct a symplectic form $\omega$ admitting a Lagrangian KB in class $F\bmod 2$ for each positive $a,h$.
 
We first construct a model for the universal cover of $X_\varepsilon$. Choose two elements in $SO(3)$
\[\rho(x,y,z)=(x,-y,-z),\qquad \tau(x,y,z)=(-x,-y,z).\]
Note that $\rho$ and $\tau$ commute. With the standard presentation
\[\pi_1(\Sigma_g)=\Big\langle a_1,b_1,\ldots,a_g,b_g\ \Big|\ \prod_{i=1}^g[a_i,b_i]=1\Big\rangle,\]
we can define a homomorphism $r_\varepsilon:\pi_1(\Sigma_g)\to SO(3)$ by assigning
\[r_\varepsilon(a_1)=\rho,\qquad r_\varepsilon(b_1)=\tau^\varepsilon,\qquad r_\varepsilon(a_i)=r_\varepsilon(b_i)=1\quad(i\geq2).\]
Let $\widetilde{\Sigma}_g$ be the universal cover of $\Sigma_g$. The quotient of the action \(\gamma\cdot(u,v)=(\gamma u,r_\varepsilon(\gamma)v)\)
on $\widetilde\Sigma_g\times S^2$ gives a sphere bundle
\[E_\varepsilon=(\widetilde\Sigma_g\times S^2)/\pi_1(\Sigma_g)\longrightarrow\Sigma_g.\]
We can identify the bundle type $E_\varepsilon\cong X_\varepsilon$ using Stiefel--Whitney classes. Let $\alpha,\beta\in H^1(\Sigma_g;\mathbb Z_2)$ be the mod-two duals of $a_1,b_1$. Since $\rho$ and $\tau$ are diagonal, the associated rank-three real vector bundle of $r_\varepsilon$ splits as \(\ell_{\varepsilon\beta}\oplus\ell_{\alpha+\varepsilon\beta}\oplus\ell_\alpha\) where $\ell_\delta$ denotes the real line bundle with $w_1(\ell_\delta)=\delta$. The Whitney sum formula thus gives
\[w_2(\ell_{\varepsilon\beta}\oplus\ell_{\alpha+\varepsilon\beta}\oplus\ell_\alpha)=(\varepsilon\beta)\cup(\alpha+\varepsilon\beta)+(\varepsilon\beta)\cup\alpha+(\alpha+\varepsilon\beta)\cup\alpha=\varepsilon\beta\cup\alpha.\]

Now choose an $SO(3)$-invariant area form $\sigma$ on $S^2$ and an area form $\eta$ on $\Sigma_g$, normalized by
\(\int_{S^2}\sigma=\int_{\Sigma_g}\eta=1.\)
Let $\widetilde\eta$ be the pullback of $\eta$ to the $\widetilde{\Sigma}_g$. For $a,h>0$, the product form
\[\widetilde\omega=a\,\operatorname{pr}_{S^2}^*\sigma+h\,\operatorname{pr}_{\widetilde\Sigma_g}^*\widetilde\eta\]
on $\widetilde{\Sigma}_g\times S^2$ is invariant under the action by $\pi_1(\Sigma_g)$ and therefore descends to a symplectic form $\omega$ on $X_\varepsilon$ such that $[\omega]\cdot F=a,\,[\omega]^2=2ah$. Hence \(\operatorname{PD}[\omega]=aS_\varepsilon+\left(h+\frac{\varepsilon a}{2}\right)F.\)

Finally, we construct the Lagrangian KB by choosing an embedded loop in $\Sigma_g$ representing $a_1$. The restriction of the $S^2$-bundle to this loop is the mapping torus of $\rho$. The equator $\{z=0\}\subset S^2$ is preserved by $\rho$, and in its angular coordinate $\rho|_C$ acts by $\theta\mapsto-\theta$. Therefore,
\[L=\bigl([0,1]\times C\bigr)/\bigl((1,v)\sim(0,\rho(v))\bigr)\subset X_\varepsilon\]
is an embedded Klein bottle. It is straightforward to see $\omega|_L=0$ in a local trivialization so that $L$ is a Lagrangian. Moreover, let $D_+=\{z\geq0\}\subset S^2$ be the upper hemisphere. Since $\rho(D_+)$ is the lower hemisphere, the image of $[0,1]\times D_+$ in the mapping torus defines a mod-two three-chain whose boundary is $L$ together with one fiber sphere. Hence $[L]_2=F\bmod2$.
\end{proof}

\subsection{Irrational ruled AE twists}
We now generalize the phenomenon of Adaloglou--Evans twist \cite{AdaloglouEvans2026} from the rational setting to the irrational ruled setting.
\begin{proposition}
\label{prop:ruled-preservation}
Let $(X',\omega')$ be obtained by KB Luttinger surgery along $L$ in a symplectic irrational ruled manifold $(X,\omega)$. Then $X'$ is also irrational ruled. Moreover, the isometry $\Phi$ of Lemma~\ref{lem:intersection-form} satisfies \(\Phi(F)=F',\) where $F'$ denotes the fiber class of $X'$.
\end{proposition}

\begin{proof}
By \eqref{equ:Excision and Poincare dual}, \(H^1(X,C;\mathbb Z)= H^2(X,C;\mathbb Z)=H^1(X',C;\mathbb Z)= H^2(X',C;\mathbb Z)=0.\) By the relative exact sequence, the restriction maps $i^*$ and $i'^*$ are therefore both isomorphisms in degree one. Choose $\alpha,\beta\in H^1(X;\mathbb Z)$ such that $\alpha\cup\beta=\operatorname{PD}(F)$. Let $\alpha',\beta'\in H^1(X';\mathbb Z)$ be the unique classes with $i'^*\alpha'=i^*\alpha$ and $i'^*\beta'=i^*\beta$, and put $f'=\alpha'\cup\beta'$. For every $z\in H_2(C;\mathbb R)$, naturality of the cup product and Lemma \ref{lem:intersection-form} give
\[\langle f',i'_*z\rangle=\langle\alpha\cup\beta,i_*z\rangle=F\cdot i_*z=\Phi(F)\cdot i'_*z.\]
Since $i'_*:H_2(C;\mathbb R)\to H_2(X';\mathbb R)$ is surjective, $\operatorname{PD}(f')=\Phi(F)$ in real homology. $\Phi$ in Lemma \ref{lem:intersection-form} also identifies the symplectic and first Chern classes, so
\[(f')^2=F^2=0,\qquad f'\cdot[\omega']=\langle[\omega],F\rangle>0,\qquad c_1(X',\omega')\cdot f'=\langle c_1(X,\omega),F\rangle=2.\]
Let $\pi:X'\to M$ be the smooth blowdown map to a symplectic minimal model $(M,\omega_M)$. Since $\pi^*$ is an isomorphism in degree one, both $\alpha'$ and $\beta'$ are pullbacks. Hence $f'=\pi^*f_M$ for their corresponding cup product $f_M\in H^2(M;\mathbb Z)$. The blowup formulas give
\[f_M^2=0,\qquad f_M\cdot[\omega_M]>0,\qquad c_1(M,\omega_M)\cdot f_M=2.\]
If $M$ were neither rational nor ruled, then, since $b_2^+(M)=1$, \cite[Theorems~A and~B]{Liu1996} would give \(c_1(M,\omega_M)^2\geq0,\,c_1(M,\omega_M)\cdot[\omega_M]\leq0.\) The light cone lemma would then imply $c_1(M,\omega_M)\cdot f_M\leq0$, a contradiction. It follows that $M$ is irrational ruled since $b_1(M)=b_1(X')=2g>0$. So, $X'$ is irrational ruled as well. 

Finally, all degree one cohomology classes on $X'$ come from the pullbacks from the base $\Sigma_g$, so their cup product $f'$ equals $m\operatorname{PD}(F')$ for some $m\in\mathbb Z$. Pairing with $c_1(X',\omega')$ gives $2=2m$. Thus $m=1$ and $\Phi(F)=F'$.

\end{proof}

\begin{corollary}\label{cor:ruled-nonvanishing}
    Let $L$ be a Lagrangian Klein bottle in a symplectic irrational ruled manifold $(X,\omega)$. Then \([L]_2\neq0,\, F\cdot[L]_2=0\). 
\end{corollary}
\begin{proof}
    If $[L]_2=0$, Remark~\ref{remark:h1-increase} would give $h_1(X')=h_1(X)+1$, contradicting that $X'$ is also irrational ruled.
To see that $F\cdot[L]_2=0$, write $\operatorname{PD}(F)=\alpha\cup\beta$ as in Proposition \ref{prop:ruled-preservation}. Recall that cup product on $H^1(L;\mathbb Z)$ vanishes. Hence $(\alpha\cup\beta)|_L=0$. Reducing modulo two and evaluating on $[L]_2$ gives $F\cdot[L]_2=0$.
\end{proof}

\begin{remark}
    No sign assumption on $c_1\cdot[\omega]$ is required in Corollary \ref{cor:ruled-nonvanishing}.
\end{remark}

\begin{corollary}
\label{cor:ruled-parity}
Let $(X_\varepsilon,\omega)$ be a minimal symplectic irrational ruled manifold. Then, KB Luttinger surgery along any Lagrangian Klein bottle produces $(X_{1-\varepsilon},\omega')$. 
Moreover, under the identification $\Phi$ in Lemma \ref{lem:intersection-form},
\[\Phi^{-1}H_2(X_{1-\varepsilon};\mathbb Z)=\mathbb ZF\oplus\mathbb Z\left(S_\varepsilon+\frac12F\right),\]
and the change of the symplectic classes is given by 
\[[\omega]=aS_\varepsilon+\left(h+\frac{\varepsilon a}{2}\right)F\mapsto[\omega']=aS_{1-\varepsilon}+\left(h+\frac{(1-\varepsilon) a}{2}\right)F'.\]
\end{corollary}

\begin{proof}
Proposition~\ref{prop:ruled-preservation} gives $[L]_2=F\bmod2$. In the notation of Lemmas \ref{lem:intersection-form} and \ref{lem:pullback-intersection-lattices},
\[\Lambda_C=\mathbb ZF\oplus\mathbb Z(2S_\varepsilon),\qquad[\Gamma:\Lambda_C]=2,\qquad\Gamma\cap\Lambda=\Lambda_C.\]
Choose some $z\in \Gamma\setminus \Lambda_C$. Since $2z\in\Lambda_C$, we may write $2z=\alpha F+2\beta S_\varepsilon$ for $\alpha,\beta\in\ZZ$. Since $2z\notin 2\Lambda$, $\alpha$ is odd. Notice that $F/2\notin\Gamma$. Otherwise, there will be the contradiction
\[1=\Phi^{-1}\big(c_1(X',\omega')\big)\cdot(F/2)\equiv(F/2)^2=0\pmod 2.\] It follows that $\beta$ must also be odd. Hence, we have \[\Gamma=\Lambda_C\oplus\ZZ(z)=\mathbb ZF\oplus\mathbb Z(S_\varepsilon+\frac12F),\qquad z^2=\alpha\beta-\varepsilon\beta^2\equiv 1-\varepsilon\pmod 2.\]
Thus the parity of the intersection form must be changed after the KB Luttinger surgery, which turns $X_\varepsilon$ to $X_{1-\varepsilon}$. Finally, by computing the self-intersection number, we have
\[S_{1-\varepsilon}=\Phi\big(S_\varepsilon+(\varepsilon-\frac12)F\big).\]
Together with $\Phi(F)=F'$, we can then conclude the symplectic period change formula.
\end{proof}

\begin{remark}
The above argument actually applies for $g=0$ as well since we may assume $[L]_2=F\bmod2$ after choosing the ruling. For $S^2\times S^2$, this recovers \cite[Proposition~3.1]{AdaloglouEvans2026} without a Maslov index computation.
\end{remark}

\begin{corollary}
\label{cor:ruled-minimal-count}
Let $g\geq1$ and suppose that the extended Auroux data
\[\widetilde{\mathcal A}=(8-8g,\,4-4g,\,P,\,Q,\,2g)\]
is realizable. If $g=1$, further assume that $P>0$. Then $\mathcal N_{\widetilde{\mathcal A}}=1$.
\end{corollary}

\begin{proof}
We first show that every realization of $\widetilde{\mathcal{A}}$ is a minimal irrational ruled surface of genus $g$. The data $\widetilde{\mathcal A}$ give $b_2^+=b_2^-=1$. If $g\geq2$, we can obtain its minimal model by at most one blowdown. So its first Chern class of the minimal model has square at most $8-8g+1<0$. It is therefore ruled by \cite[Theorem~A]{Liu1996}. If $g=1$, the assumption $P=c_1\cdot[\omega]>0$ implies that the realization is rational or ruled by \cite[Theorem~B]{Liu1996}. In either case, we see that the underlying manifold is a ruled surface. Also note that $b_1=2g$ implies that the base genus is $g$, and $b_2=2$ implies that the ruled surface is minimal.

Now, one can compute the fiber area $a$ satisfies
\[P=\frac{Q}{a}-2(g-1)a.\]
The right hand side is strictly decreasing for $a>0$, so $P,Q$ determine $a$ and $h=Q/(2a)$. These parameters determine the symplectomorphism class within each bundle type by \cite{LalondeMcDuff1996Ruled}. Lemma~\ref{lem:ruled-construction} and Corollary \ref{cor:ruled-parity} connect the two bundle types with these same parameters by one KB Luttinger surgery. Thus all integral realizations belong to a single KB Luttinger equivalence class.
\end{proof}

\begin{remark}
When $g=1$ and $P=0$, $\widetilde{\mathcal{A}}$ can not be realized by a ruled surface. However, it can be realized as certain symplectic $T^2$-bundles over $T^2$; see Proposition \ref{prop:SCY-infinite-KB-count}.
\end{remark}

\subsection{The reflection formula on blowups}

Assume now $k\geq1$. The two bundle types become diffeomorphic after one blowup. So we use the fixed model
\[X_{g,k}=(S^2\times\Sigma_g)\#k\overline{\mathbb{CP}}^{2},\qquad\Lambda=\mathbb ZS\oplus\mathbb ZF\oplus\bigoplus_{i=1}^k\mathbb ZE_i,\qquad c_1=2S+(2-2g)F-\sum_{i=1}^k E_i\]
with $S^2=F^2=0$ and $S\cdot F=1$. Analogous to the rational case, we define
\[\Delta_{c_1,F}^{(-4)}=\{v\in\Lambda:v^2=-4,\ c_1\cdot v=F\cdot v=0\},\]
and its subset for any $\ZZ_2$-class $\beta$
\[\Delta_{c_1,F}^{(-4)}(\beta)=\{v\in\Delta_{c_1,F}^{(-4)}:v\bmod2=\beta\}.\]

\begin{lemma}
\label{lem:ruled-KB-lift}
For every Lagrangian Klein bottle $L$ in $X_{g,k}$,
\[[L]_2=F\bmod 2\quad\text{or}\quad[L]_2=\delta F+\sum_{i\in I}E_i,\qquad\delta\in\{0,1\},\quad |I|=4.\]
In particular, the set \(\Delta_{c_1,F}^{(-4)}([L]_2)\) is nonempty.
\end{lemma}

\begin{proof}
We consider the Luttinger surgery along $L$ and retain the notations in previous sections. By Corollary \ref{cor:ruled-nonvanishing}, we can write $[L]_2=\delta F+\sum_{i\in I}E_i\neq0$. Audin's congruence \cite{Audin} gives $|I|\equiv0\pmod4$. For the integral lift $b=\delta F+\sum_{i\in I}E_i$, Lemmas \ref{lem:intersection-form} and \ref{lem:pullback-intersection-lattices} give
\[\Lambda_C=\{x\in\Lambda:x\cdot b\equiv0\pmod2\},\qquad[\Gamma:\Lambda_C]=2,\qquad\Gamma\cap\Lambda=\Lambda_C.\]
Since $\Phi(F)=F'$ by Proposition \ref{prop:ruled-preservation}, $\Gamma\cap F^\perp\cong F'^\perp\subset \Lambda'$. Hence, $\Gamma\cap F^\perp$ contains $k$ pairwise orthogonal $(-1)$-classes. For any $z\in (\Gamma\cap F^\perp)\setminus \Lambda_C$, by the same argument in Lemmas \ref{lem:KB-minus-four-lift} and \ref{lem:KB-reflection-marking}, we have that $2z-b\in2\Lambda$. It follows that $z$ has half-integral $E_i$-coefficients for every $i\in I$ and zero coefficient for $S$. $z^2$ is therefore at most $-|I|/4$. Consequently, if $|I|\geq 8$, the $(-1)$-classes in $\Gamma\cap F^\perp$ can only be the classes in $\Lambda_C\cap F^\perp$, which are of the form $m F\pm E_j$, for some $m\in \ZZ$ and $j\notin I$. There are at most $k-|I|$ of these which are pairwise orthogonal, which is a contradiction. Thus $|I|=0$ or $4$, and the nonvanishing of $[L]_2$ gives $[L]_2=F\bmod2$ in the first case.

To see $\Delta_{c_1,F}^{(-4)}([L]_2)\neq\varnothing$, notice that the class $F\bmod2$ lifts to $F-2E_1$; if $|I|=4$, choose signs $m_i\in\{\pm1\}$ so that $d=\frac12\sum_{i\in I}m_i\equiv\delta\pmod2$. Then \(v=dF-\sum_{i\in I}m_iE_i\) has $v^2=-4$, $c_1\cdot v=F\cdot v=0$, and the required reduction.
\end{proof}

Now we can consider the reflection map $R_v:x\mapsto x+\frac{(x\cdot v)}{2}v$ similar to the rational case.

\begin{lemma}
\label{lem:ruled-reflection}
For the KB Luttinger surgery on $L\subset(X,\omega)$ and any $v\in\Delta_{c_1,F}^{(-4)}([L]_2)$, the reflection $R_v$ satisfies \(R_v(\Lambda_C)=\Lambda_C,\, R_v(\Gamma)=\Lambda.\)
Consequently $\mathfrak m_v=R_v\Phi^{-1}:\Lambda'\to\Lambda$ is an integral isometry preserving $c_1$ and $F$, and the new period is $R_v(\operatorname{PD}[\omega])$.
\end{lemma}

\begin{proof}
\(R_v(\Lambda_C)=\Lambda_C\) and \(R_v(\Gamma)=\Lambda\) follow from the same lattice argument of Lemma~\ref{lem:KB-reflection-marking}. The remaining claims follow from Lemma \ref{lem:intersection-form} and $R_v(c_1)=c_1$, $R_v(F)=F$.
\end{proof}

We provide the explicit reflection formulas. Write \(v=dF-\sum_i m_iE_i\in\Delta^{(-4)}_{c_1,F}.\) By Lemma \ref{lem:ruled-KB-lift}, $v$ is either $\pm(F-2E_i)$ or a class with four coefficients $m_i=\pm1$. For $\operatorname{PD}[\omega]=aS+hF-\sum_i\delta_iE_i$, put
\[s=\operatorname{PD}[\omega]\cdot v=da-\sum_i m_i\delta_i.\] Then
\begin{equation}
\label{eq:ruled-reflection-coordinates}
a'=a,\qquad h'=h+\frac{ds}{2},\qquad\delta_i'=\delta_i+\frac{m_i s}{2}.
\end{equation}
In particular, $R_{F-2E_i}$ sends
\[\delta_i\longmapsto a-\delta_i,\qquad h\longmapsto h+\frac a2-\delta_i\]
and fixes the other coefficients.

\subsection{The surgery graphs and the actual count}\label{sec:ruled-KB-graph}
We now define the surgery graphs using the dyadic lattice similar to Section \ref{subsec:rational-KB-graphs}. For $X_{g,k}$, its $c_1$-symplectic cone in the component $x\cdot F>0$ is
\begin{equation}
\label{eq:ruled-cone}
\mathcal C_{g,k}
=\bigg\{x=aS+hF-\sum_i\delta_iE_i:
a>0,\ x^2>0,\ 0<\delta_i<a\ \text{for every }i\bigg\},
\end{equation}
since all exceptional classes which have intersection pairing one with $c_1$ are either $E_i$ or $F-E_i$. Two symplectic forms with the same cohomology class are symplectomorphic through a diffeomorphism acting trivially on $H_2$; see \cite[Lemma~2.11]{HK19}. Consider the {\bf Weyl group} \[W_{g,k}=\bigl\langle R_{E_i-E_j},\,R_{F-E_i-E_j}\mid 1\le i<j\le k\bigr\rangle\] generated by the reflections along the $(-2)$-classes $E_i-E_j$ and $F-E_i-E_j$. It is equal to the group $O(\Lambda,c_1,F)$ of lattice automorphisms preserving both $c_1$ and $F$, and can be induced by orientation-preserving diffeomorphisms; see \cite[\S 4.2]{LiWu2012}. In particular, every such lattice automorphism preserves $\mathcal{C}_{g,k}$. Fix $P\in\mathbb Z$ and $Q\in\mathbb Z_{>0}$. We let
\[
\begin{split}
\mathbb V_{g,k;P,Q}
&=\{x\in\mathcal C_{g,k}\cap\Lambda[1/2]:
c_1\cdot x=P,\ x^2=Q\},\\
\mathbb I_{g,k;P,Q}
&=\mathbb V_{g,k;P,Q}\cap\Lambda.
\end{split}
\]

\begin{definition}[KB Luttinger surgery graphs for irrational ruled manifolds]
\label{def:ruled-graphs}
Define two simple graphs $\mathbb G^{\mathrm{alg}}_{g,k;P,Q}$ and $\mathbb G^{\mathrm{geo}}_{g,k;P,Q}$ with the same vertex set $\mathbb V_{g,k;P,Q}$. Both contain the {\bf Weyl edges} $x\leftrightarrow w(x)$ for $w\in W_{g,k}$, with distinct endpoints. Their additional {\bf Luttinger edges} are as follows:
\begin{enumerate}[nosep]
\item The {\bf algebraic graph} $\mathbb G^{\mathrm{alg}}_{g,k;P,Q}$ contains $x\leftrightarrow R_v(x)$ whenever $v\in\Delta_{c_1,F}^{(-4)}$ and both endpoints are distinct vertices.
\item The {\bf geometric graph} $\mathbb G^{\mathrm{geo}}_{g,k;P,Q}$ contains such an edge if, in addition, a symplectic form in class $x$ admits a Lagrangian Klein bottle $L$ with $[L]_2=v\bmod2$.
\end{enumerate}
\end{definition}
As before, for $\star\in\{\mathrm{alg},\mathrm{geo}\}$, define
\[N^\star_{g,k}(P,Q)=\#\{T\in\pi_0(\mathbb G^\star_{g,k;P,Q}):T\cap\mathbb I_{g,k;P,Q}\neq\varnothing\},\] where the path $T$ may go through non-integral dyadic vertices. These are exactly the counts of all connected components in the graphs $\mathbb X(g,k,P,Q)$ and $\mathbb Y(g,k,P,Q)$ introduced in Section \ref{sec:intro-KB-Luttinger}. One could instead include algebraic Luttinger edges $x\leftrightarrow AR_v(x)$ with $A\in O(\Lambda,c_1,F)$, exactly as in Definition \ref{def:rational-KB-graphs}. However, this is unnecessary since we know $O(\Lambda,c_1,F)= W_{g,k}$ preserves the $c_1$-symplectic cone, which gives the same count of connected components.

\begin{proposition}
\label{prop:ruled-count}
Let $g\geq1$, $k\geq1$, $P\in\mathbb Z$ and $Q\in\mathbb Z_{>0}$. If $g=1$, further assume $P>-k$. For realizable $\widetilde{\mathcal{A}}$ given by \eqref{eq:ruled-data},
\[N^{\mathrm{alg}}_{g,k}(P,Q)\leq \mathcal N_{\widetilde{A}}=N^{\mathrm{geo}}_{g,k}(P,Q)<\infty.\] 
\end{proposition}

\begin{proof}
The data \eqref{eq:ruled-data} force $b_2^+=1$ and $b_2^-=k+1$. If an integral realization has a non-ruled minimal model $(X_{\operatorname{min}},\omega_{\operatorname{min}})$ after $t$ symplectic blowdowns, then $t\leq k+1$ and \cite[Theorem A]{Liu1996} gives
\[0\leq c_{1,\min}^2=8-8g-k+t\leq9-8g.\]
This is impossible for $g\geq2$. If $g=1$, it implies that $t\geq k$. The minimal model satisfies $c_{1,\min}\cdot[\omega_{\min}]\leq0$ by \cite[Theorem B]{Liu1996}, and each exceptional area is a positive integer. Therefore
\[P=c_{1,\min}\cdot[\omega_{\min}]-\sum_{j=1}^t\omega(E_j)\leq-t\leq-k.\]
This excludes non-ruled realizations when $P>-k$. By Proposition~\ref{prop:ruled-preservation}, the Luttinger equivalence sequence stays inside the genus $g$ irrational ruled category.

The classification result of symplectic structures recalled above identifies each vertex with a symplectomorphism class. Also, different markings differ by $W_{g,k}$ and hence lie in the same graph component. Consider an actual Luttinger surgery on a symplectic irrational ruled manifold $(X,\omega)$ with $[\omega]=x$. By Proposition~\ref{prop:ruled-preservation} and Lemma \ref{lem:ruled-reflection}, the surgered manifold $(X',\omega')$ admits the marking $\mathfrak m_v$ with period $R_v(x)$, preserving $c_1,F$. Choose also a diffeomorphism $X'\cong X$ inducing a marking $\mathfrak n:\Lambda'\to\Lambda$. Then \(A=\mathfrak n\circ\mathfrak m_v^{-1}\in O(\Lambda,c_1,F)=W_{g,k}.\)
Consequently $\mathfrak m_v$ is itself induced by a diffeomorphism from $X'$ to $X$ and $R_v(x)=A^{-1}\mathfrak n(\operatorname{PD}[\omega'])$ lies in $\mathcal C_{g,k}$. The actual surgery therefore gives a geometric Luttinger edge or leaves the vertex unchanged followed by a Weyl edge. Moreover, formula \eqref{eq:ruled-reflection-coordinates} shows that a surgery sequence starting at an integral period has dyadic periods throughout.

Conversely, similar to the proof of Proposition \ref{prop:rational-KB-graph-bounds}, we know a Weyl edge can be realized by a symplectomorphism; a geometric Luttinger edge can be realized by an actual KB surgery since elements in $O(\Lambda,c_1,F)$ are all induced by diffeomorphisms. This proves $\mathcal N_{\widetilde{\mathcal{A}}}=N^{\mathrm{geo}}_{g,k}$. Since the geometric graph is a subgraph of the algebraic graph, the inequality $N^{\mathrm{alg}}_{g,k}(P,Q) \leq \mathcal N_{\widetilde{A}}$ follows.

It remains to show that the integral vertex set is finite.
Eliminating $h$ from the defining equations for $P,Q$ gives
\begin{equation}
\label{eq:ruled-finiteness}
Q=aP+2(g-1)a^2+\sum_i\delta_i(a-\delta_i).
\end{equation}
At an integral vertex, $a,h,\delta_i$ are integers and $1\leq\delta_i\leq a-1$. If $g\geq2$, \eqref{eq:ruled-finiteness} bounds $a$, and hence all $\delta_i$ and $h$. For $g=1$, set $t_i=\min\{\delta_i,a-\delta_i\}$, $T=\sum_i t_i$, and $U=\sum_i t_i^2$. Then
\[Q=a(P+T)-U,\qquad U\leq\frac{aT}{2}.\]
If $T+2P\geq1$, then $Q\geq a/2$. Otherwise $T\leq-2P$, while $Q>0$ forces the integer $P+T\geq1$. Hence
\[a=\frac{Q+U}{P+T}\leq Q+T^2\leq Q+4P^2.\]
Thus $a$ is still bounded, proving finiteness.
\end{proof}

The following geometric realization should be compared with Proposition \ref{prop:small-blowup-realization}.

\begin{proposition}
\label{prop:ruled-small-realization}
Let $g\geq1$ and $k\geq1$. Every $x\in\mathcal C_{g,k}$ has a symplectic representative containing a Lagrangian Klein bottle $L$ with $[L]_2=F\bmod 2$. The symplectic manifold after the KB Luttinger surgery along $L$ can be represented by $R_{F-2E_i}(x)$ for any $1\leq i\leq k$. In particular, if $k\leq3$, every algebraic Luttinger edge is geometric.
\end{proposition}

\begin{proof}
For each $i$, choose $D_i\in\{E_i,F-E_i\}$ so that
\[t_i=x\cdot D_i=\min\{\delta_i,a-\delta_i\}.\]
Blowing down disjoint representatives of the exceptional classes $D_i$ gives a minimal model $X_\varepsilon$ with symplectic class 
\[\operatorname{PD}[\omega_{\operatorname{min}}]=aS+\big(H+\frac{\varepsilon a}{2}\big)F,\qquad H=\frac{x^2+\sum_i t_i^2}{2a}.\]
It suffices to realize the Lagrangian KB in a form with fiber area $a$, exceptional areas $t_i$, and square $x^2$. To achieve this, we use the symplectic fiber sum construction similar to \cite[Proposition 3.1]{DaiHoLi2019}.

For the first summand, start with the model of $X_\varepsilon$ from Lemma~\ref{lem:ruled-construction} with parameters $(h,a)=(H,2r)$, where $r>0$ is a small number chosen later. Let $L$ be the equatorial Klein bottle over $\gamma$. The north and south poles trace out a horizontal symplectic surface $B$ disjoint from $L$, which is connected and projects to $\Sigma_g$ with degree two. Thus $B$ has genus $2g-1$ and area $2H$. In fact, $[B]=[S_\varepsilon]+[S_{-\varepsilon}]$, so $[B]^2=0$. Put
\[J=\{i:t_i=a/2\},\qquad m=|J|.\]
For each $i\in J$, choose a ball with capacity $r$ centered at a point of $B$ on the fibers over points outside $\gamma$. For small enough $r$, we can choose $m$ disjoint such balls since the base area is fixed and the fiber area is $2r$. The proper transform $B'$ of $B$ thus has square $-m$ and area $u=2H-mr$.

For the second summand, take another ruled surface $Y$ of genus $2g-1$ with a section $Z$ of square $m$. Set $b=a/2-r$. For sufficiently small $r$, we have $u>0$, $0<t_i<b$ for every $i\notin J$, and
\[2bu-mb^2-\sum_{i\notin J}t_i^2=x^2-(4rH-mr^2)>0.\]
By the description of the symplectic cone in \eqref{eq:ruled-cone}, we can choose a symplectic form on the $(k-m)$-fold blowup with area $u$ on $Z$, area $b$ on the fiber, and exceptional area $t_i$ for each $i\notin J$. By \cite[Theorem~2.13]{DorfmeisterLi2010}, we may choose the symplectic form on the blowup such that the section surface $Z$ is a symplectic submanifold.

Now take the symplectic fiber sum 
\[\big(X_\varepsilon\# m\overline{\mathbb{CP}}^{2}\big)\#_{B'=Z}\big(Y\# (k-m)\overline{\mathbb{CP}}^{2}\big),\] which is diffeomorphic to $X_{g,k}$ by \cite[Remarks~1.3 and~1.5]{Usher2009}. The fiber sum region is disjoint from $L$, so $L$ remains Lagrangian. The fiber sphere in the summed manifold is obtained by gluing a fiber sphere of the first summand with two fiber spheres of the second summand. Also, each exceptional sphere in the first summand indexed by $J$ is glued with one fiber sphere in the second summand to form an exceptional sphere after summing. In the summed manifold, the symplectic areas of the fiber class and exceptional classes indexed by $J$ are therefore \(2r+2b=a\) and \( r+b=a/2.\) The remaining exceptional areas are $t_i$, and the squares of the two summands add to $x^2$. Hence the resulting symplectic period is $x$.

Since all exceptional spheres are away from $L$, we see that $[L]_2=F\bmod 2$ by Lemma \ref{lem:ruled-construction}. Each $F-2E_i$ is therefore a $(-4)$-lift of $[L]_2$. By Lemma~\ref{lem:ruled-reflection}, we can choose a marking for every $i$ such that $R_{F-2E_i}(x)$ denotes the symplectic class after the surgery.

Finally, when $k\leq3$, the proof of Lemma \ref{lem:ruled-KB-lift} gives
\[\Delta_{c_1,F}^{(-4)}=\{\pm(F-2E_i):1\leq i\leq k\}.\]
In particular, all these classes reduce to $F\bmod 2$. Thus the reflections realized above exhaust the algebraic Luttinger edges in this range.
\end{proof}
For $1\leq k\leq 3$, let us now introduce the set
\[\mathcal T_{g,k}(P,Q)=\left\{(a,t_1,\ldots,t_k)\in\mathbb Z^{k+1}\;\middle|\;
\begin{aligned}
&1\leq t_1\leq\cdots\leq t_k\leq a/2,\\
&Q=aP+2(g-1)a^2+\sum_{i=1}^k t_i(a-t_i)
\end{aligned}
\right\}.\]

\begin{theorem}\label{thm:ruled-small-exact-count}
Let $g\geq1$, $1\leq k\leq3$, $P\in\ZZ$, $Q\in\ZZ_{>0}$. If $g=1$, further assume $P>-k$. Then
\[\mathcal{N}_{\widetilde{\mathcal A}}=N^{\mathrm{alg}}_{g,k}(P,Q)=N^{\mathrm{geo}}_{g,k}(P,Q)=\#\mathcal T_{g,k}(P,Q).\]
\end{theorem}

\begin{proof}
By Propositions \ref{prop:ruled-count} and \ref{prop:ruled-small-realization}, we already have $\mathcal{N}_{\widetilde{\mathcal A}}=N^{\mathrm{alg}}_{g,k}(P,Q)=N^{\mathrm{geo}}_{g,k}(P,Q)$.

To identify them with $\#\mathcal T_{g,k}(P,Q)$, we explain how a component in the graph corresponds to an element in $\mathcal T_{g,k}(P,Q)$. First, by eliminating $h$ from the defining equations of $P,Q$, we have the displayed equation in the definition of $\mathcal{T}_{g,k}(P,Q)$. By \eqref{eq:ruled-reflection-coordinates}, a Weyl wedge permutes the exceptional areas and may replace some $(\delta_i,\delta_j)$ by $(a-\delta_j,a-\delta_i)$, and a Luttinger edge turns exactly one $\delta_i$ into $a-\delta_i$. Hence, every edge preserves the fibre area $a$ and the unordered collection
\[\bigl\{t_i=\min\{\delta_i,a-\delta_i\}:1\leq i\leq k\bigr\}.\] A component of the graph is therefore determined by $a$ and this unordered collection of $t_i$'s and distinct tuples in $\mathcal T_{g,k}(P,Q)$ belong to distinct components. 

Conversely, for each element in $\mathcal{T}_{g,k}(P,Q)$, choose $\delta_i\in\{t_i,a-t_i\}$ and set
\[h=\frac{P+2(g-1)a+\sum_i\delta_i}{2}.\]
These periods have square $Q>0$ and satisfy $0<\delta_i<a$, so they belong to $\mathcal C_{g,k}$. If $a$ is odd, replacing one $\delta_i$ by $a-\delta_i$ switches the parity of the numerator, so one choice gives integral $h$. If $a$ is even, the equation gives $\sum_i t_i\equiv Q\equiv P\pmod2$, and $h$ is integral with $\delta_i=t_i$. Thus every listed tuple gives exactly one component meeting $\mathbb I_{g,k;P,Q}$. We therefore have $\#\mathcal T_{g,k}(P,Q)=\mathcal{N}_{\widetilde{\mathcal A}}.$
\end{proof}

As an application of Theorem \ref{thm:ruled-small-exact-count}, we can obtain both families of extended Auroux data with $\mathcal{N}_{\widetilde{\mathcal A}}=1$ and $\mathcal{N}_{\widetilde{\mathcal A}}>1$. Let $g\geq1$ and $1\leq k\leq3$, then we have
\[\mathcal N_{(8-8g-k,\,4-4g+k,\,0,\,18(g-1)+2k,\,2g)}=1.\]
Indeed, the counting equation gives
\[18(g-1)+2k=2(g-1)a^2+\sum_{i=1}^k t_i(a-t_i)\geq 2(g-1)a^2+k(a-1),\] which forces $a\leq3$. The value $a=2$ gives $8(g-1)+k<18(g-1)+2k$, so the unique solution is
\((a,t_1,\ldots,t_k)=(3,1,\ldots,1).\)
For $k=2$, we also have the similar example
\[\mathcal N_{(6-8g,\,6-4g,\,0,\,8g-6,\,2g)}=1.\]
Here the same inequality gives
\[8g-6\geq 2(g-1)a^2+2(a-1)\geq8g-6,\] forcing the unique solution $(a,t_1,t_2)=(2,1,1)$.

To obtain an example with the count greater than one, choose an integer $r\geq2$ such that \[P_r=kr(r-2)-2(g-1)(4r+1)\geq0,\qquad Q_r=2rP+\bigl(8(g-1)+k\bigr)r^2.\] We claim that
\[\mathcal N_{(8-8g-k,\,4-4g+k,\,P_r,\,Q_r,\,2g)}=2.\]
To see this, first note that every such tuple satisfies
\[aP+2(g-1)a^2+k(a-1)\leq Q\leq aP+\bigr(2(g-1)+\frac{k}{4}\bigr)a^2.\]
Both bounds are strictly increasing in $a\geq2$. The upper bound equals $Q_r$ at $a=2r$, and the lower bound equals $Q_r$ at $a=2r+1$. Thus $a\in\{2r,2r+1\}$, and equality case gives precisely the two solutions \((a,t_1,\ldots,t_k)=(2r,r,\ldots,r),\,(2r+1,1,\ldots,1).\)

\section{KB Luttinger surgery in Kodaira dimension zero}\label{sec:Kod-0}

Recall that a symplectic Calabi--Yau surface is a closed symplectic four-manifold with torsion first Chern class, or equivalently a minimal symplectic four-manifold of Kodaira dimension zero \cite{Li2006Quaternionic}. The following proposition shows that Klein bottle Luttinger surgery is a symplectic Calabi--Yau surgery, giving the Klein bottle analogue of \cite[\S~4.3]{HoLi2012}.

\begin{proposition}
\label{prop:klein-bottle-canonical-order}
Let $(X',\omega')$ be obtained from a closed symplectic four-manifold $(X,\omega)$ by Luttinger surgery along a Lagrangian Klein bottle $L$. For every positive integer $m$,
\[m\,c_1(X,\omega)=0\quad\Longleftrightarrow\quad m\,c_1(X',\omega')=0.\]
\end{proposition}

\begin{proof}
We use the notation from Section \ref{section:torus-KB-comparison}. \eqref{equ:Excision and Poincare dual} gives
\[H^2(X,C;\mathbb Z)\cong H_2(L;\mathbb Z)=0,\qquad H^2(X',C;\mathbb Z)\cong H_2(L';\mathbb Z)=0.\]
Thus the restriction maps
\[i^*:H^2(X;\mathbb Z)\longrightarrow H^2(C;\mathbb Z),\qquad(i')^*:H^2(X';\mathbb Z)\longrightarrow H^2(C;\mathbb Z)\]
are both injective. Since the symplectic forms agree on the common complement $C$, $i^*TX$ and $(i')^*TX'$ are isomorphic as symplectic vector bundles so that \[c_1(i^*TX)=c_1((i')^*TX')\in H^2(C;\ZZ)\] By the naturality of the first Chern class,
\[i^*c_1(X,\omega)=(i')^*c_1(X',\omega')\in H^2(C;\mathbb Z).\]
Consequently, for every positive integer $m$,
\[m\,c_1(X,\omega)=0\Longleftrightarrow m\,i^*c_1(X,\omega)=0\Longleftrightarrow m\,(i')^*c_1(X',\omega')=0\Longleftrightarrow m\,c_1(X',\omega')=0.\]
\end{proof}

We next establish two results according to the Euler characteristic: homological non-triviality of Lagrangian Klein bottles when $\chi(X)>0$, and the existence of fixed extended Auroux data admitting infinitely many KB Luttinger equivalence classes when $\chi(X)=0$.

\subsection{Positive Euler characteristic: homological non-triviality}\label{section:positive-Euler-Kod-0}

Smirnov \cite[\S6]{Smirnov} proved that every embedded Lagrangian Klein bottle in a symplectic homotopy K3 surface with $c_1=0$ has nonzero mod-two homology class. His proof combines the homological effect of KB Luttinger surgery obtained by Nemirovski \cite{Nemirovski2009} with the gauge theoretic results of Taubes \cite{TaubesConstraints} and Morgan--Szab\'o \cite{MorganSzabo}. With the Enriques surface case in mind, we extend Smirnov's result by combining the preceding proposition with a covering argument. Let us first recall the result of Li.

\begin{theorem}[\cite{Li2006Quaternionic}]\label{thm:Li-Kod-0}
    Let $(X,\omega)$ be a symplectic Calabi--Yau surface. Then its Euler charateristic number $\chi(X)\in\{0,12,24\}$.
\end{theorem}
Applying this theorem to finite covers gives the integral homology information needed below.

\begin{corollary}
\label{cor:SCY-positive-euler-homology}
Let $(X,\omega)$ be a symplectic Calabi--Yau surface with
$\chi(X)>0$. Then
\[H_1(X;\ZZ)\cong\operatorname{Tor}H_2(X;\ZZ)\cong
\begin{cases}
0, & \chi(X)=24,\\
\ZZ_2, & \chi(X)=12.
\end{cases}\]
If $\chi(X)=12$, then $\sigma(X)=-8$ and the unique nonzero torsion class in $H_2(X;\ZZ)$ is $\operatorname{PD}(c_1(X,\omega))$.
\end{corollary}

\begin{proof}
Every finite covering space $p:\widetilde X\to X$ equipped with the pullback symplectic form $p^*\omega$ has torsion first Chern class. Theorem~\ref{thm:Li-Kod-0} therefore gives \(\deg(p)\,\chi(X)=\chi(\widetilde X)\leq24.\) If $b_1(X)>0$, a surjection $\pi_1(X)\to\mathbb Z$ yields connected cyclic covers of every finite degree $d$, contradicting $d\chi(X)\le24$ for sufficiently large $d$. Thus $b_1(X)=0$, and the finitely generated abelian group $H_1(X;\mathbb Z)$ is finite. Moreover, applying this inequality to the covering space corresponding to the kernel of the abelianization homomorphism $\pi_1(X)\to H_1(X;\ZZ)$ yields \(|H_1(X;\ZZ)|\,\chi(X)\leq24.\) Thus $H_1(X;\ZZ)=0$ if $\chi(X)=24$, while $H_1(X;\ZZ)$ is either zero or $\ZZ_2$ if $\chi(X)=12$. In the latter case, \(0=c_1(X,\omega)^2=2\chi(X)+3\sigma(X)\) gives $\sigma(X)=-8$. Rokhlin’s congruence then implies that $X$ is not spin, so $c_1(X,\omega)\bmod2=w_2(X)\neq0$. Since $c_1(X,\omega)$ is torsion, the universal coefficient theorem rules out $H_1(X;\ZZ)=0$. Hence $H_1(X;\ZZ)\cong\ZZ_2$. the Finally, the assertions about $\operatorname{Tor}H_2(X;\ZZ)$ follow from Poincar\'e duality and the universal coefficient theorem.
\end{proof}

Next, we need Guillou--Marin's refinement of Rokhlin's congruence which has been mentioned in \cite{Smirnov}. Let us closely follow the exposition of Kirby--Taylor \cite{KirbyTaylor1990}. Write $V=H_1(F;\ZZ_2)$ and let $\bullet:V\times V\to\ZZ_2$ be the mod-two intersection pairing. A {\bf quadratic enhancement} of $\bullet$ is a function $q:V\to\ZZ_4$ satisfying
\begin{equation}\label{equ:quadratic-enhancement}
    q(x+y)\equiv q(x)+q(y)+2\cdot x\bullet y\pmod4.
\end{equation}
Here $2\cdot:\ZZ_2\to\ZZ_4$ denotes the homomorphism sending $1\in \ZZ_2$ to $2\in\ZZ_4$.
A $\operatorname{Pin}^{-}$ structure on $F$ is equivalently a spin structure on $TF\oplus\det(TF)$. It determines such an enhancement by the $\text{mod } 4$ count of the right half twists of the normal line $\nu_{K\subset F}$ along any embedded circle $K\subset F$, measured in an odd framing of $\nu_{K\subset F}\oplus\det(TF)|_K$ relative to the induced spin structure. This gives a bijection between $\operatorname{Pin}^{-}$ structures on $F$ and quadratic enhancements of its intersection pairing
\cite[Theorem~3.2]{KirbyTaylor1990}.
The {\bf Brown--Arf invariant} $\beta(q)\in\ZZ_8$ is the algebraic invariant associated to $q$ defined by the Gauss sum
\begin{equation}\label{equ:Brown}
    \sum_{x\in V}\exp\!\left(\frac{2\pi iq(x)}{4}\right)=\sqrt{|V|}\exp\!\left(\frac{2\pi i\beta(q)}{8}\right).
\end{equation}

Now let $X$ be a closed oriented smooth four-manifold and $F\subset X$ an embedded characteristic surface, namely $[F]_2=\operatorname{PD}(w_2(X))$.
For such a surface, one can choose a $\operatorname{Pin}^-$ structure on $X\setminus F$ which does not extend across any component of $F$.
This choice, called a characterization of $(X,F)$ in \cite{KirbyTaylor1990}, induces a $\operatorname{Pin}^{-}$ structure on $F$ and hence a quadratic enhancement $q$. Denote by $F\bullet F$ the normal Euler number of $F\subset X$. When $F$ is non-orientable, this is $\frac{1}{2}e(\varphi^*\nu_{F\subset X})$ for the orientable double cover $\varphi:\widetilde{X}\to X$.
The Guillou--Marin congruence formulated in \cite[Theorem 6.3]{KirbyTaylor1990}, with a reverse of sign on $\beta(q)$ (see the footnote after
\cite[Theorem 6.18]{ConwayOrsonPowell2025}), is the following.
\begin{theorem}[\cite{KirbyTaylor1990}]\label{thm:Guillou-Marin}
\(\sigma(X)\equiv F\bullet F+2\beta(q)\pmod{16}.\)
\end{theorem}

This formulation requires no assumption on $H_1(X;\ZZ)$ (compare \cite[\S5]{Matsumoto}). As a consequence, we obtain the following obstruction which will be used to exclude the nonzero torsion reduction in the Enriques case.

\begin{corollary}
\label{cor:characteristic-klein-bottle}
Let $X$ be a closed connected oriented smooth four-manifold with $\sigma(X)\equiv8\pmod{16}$. Then there is no characteristic Klein bottle in $X$ with normal Euler number zero.
\end{corollary}

\begin{proof}
Suppose there exists such a characteristic Klein bottle. Theorem \ref{thm:Guillou-Marin} gives $\beta(q)=4\in\ZZ_8$. The $\ZZ_2$-intersection form of a Klein bottle has a basis $a,b$ with $a^2=b^2=1,\,a\cdot b=0$. By \eqref{equ:quadratic-enhancement}, we have $q(a+b)=q(a)+q(b)\in\ZZ_4$ and $q(a),q(b)\in\{1,-1\}\subset\ZZ_4$. Thus \eqref{equ:Brown} gives 
\[\exp\left(\frac{\pi i\beta(q)}{4}\right)=\frac{1}{2}\left(1+i^{q(a)}+i^{q(b)}+i^{q(a)+q(b)}\right).\]
 Consequently,
\[\beta(q)=
\begin{cases}
2,  & q(a)=q(b)=1,\\
0,  & q(a)=-q(b),\\
-2, & q(a)=q(b)=-1,
\end{cases}
\qquad\text{in }\ZZ_8.\]
In particular, $\beta(q)\neq4 \in \ZZ_8$, a contradiction.
\end{proof}

Now we prove the homological non-triviality result of Lagrangian Klein bottle in symplectic four-manifolds with $c_1\cdot [\omega]=0$.
\begin{proof}[Proof of Theorem~\ref{thm:KB-homological-essentiality}]
Suppose that \([L]_2=\rho(\tau)\) for some \(\tau\in\operatorname{Tor}H_2(X;\ZZ).\)
We first show that $c_1(X,\omega)$ is torsion.
Otherwise, choose a closed real two-form $\eta$ such that $c_1(X,\omega)\cdot[\eta]>0$.
Since a Weinstein neighborhood $N$ of $L$ satisfies
$H^2(N;\mathbb R)=0$, a cutoff modification by an exact form allows us to arrange that $\eta$ vanishes near $L$.
For sufficiently small $t>0$, the form
$\omega_t=\omega+t\eta$ is symplectic, $L$ remains Lagrangian,
and $c_1(X,\omega_t)=c_1(X,\omega)$. Hence
\[
c_1(X,\omega_t)\cdot[\omega_t]
=t\,c_1(X,\omega)\cdot[\eta]>0.
\]
By Liu--Ohta--Ono's theorem \cite{Liu1996,OhtaOno1996}, $X$ is rational or ruled. In particular, $H_2(X;\ZZ)$ is torsion-free so that $\tau=0$ and $[L]_2=0$.
This contradicts Nemirovski--Shevchishin's theorem \cite{Nemirovski2009,Shevchishin2009}. Therefore $(X,\omega)$ is a symplectic Calabi--Yau surface.

Now, suppose first that $[L]_2=0$. Perform the Klein bottle Luttinger surgery along $L$ to obtain $(X',\omega')$. Proposition~\ref{prop:klein-bottle-canonical-order} shows that $(X',\omega')$ is also a symplectic Calabi--Yau surface. By Lemma \ref{lemma:surgery-characteristic-numbers} and Remark \ref{remark:h1-increase}, we have
\[\chi(X')=\chi(X)>0,\qquad h_1(X')=h_1(X)+1.\]
On the other hand, Corollary~\ref{cor:SCY-positive-euler-homology} determines $h_1$ from in the case of positive Euler characteristic, so $h_1(X')=h_1(X)$, a contradiction.

It remains to consider $[L]_2\neq0$. In this case, $\tau\neq0$ and Corollary~\ref{cor:SCY-positive-euler-homology} gives
\[\chi(X)=12,\qquad\sigma(X)=-8,\qquad\tau=\operatorname{PD}(c_1(X,\omega)).\]
Consequently,
\([L]_2=\rho(\tau)=\operatorname{PD}(w_2(X)).\)
So $L$ is a characteristic surface. Since $L$ is Lagrangian, $\nu_{L\subset X}\cong TL$. It follows that $L\bullet L=0$ since the Klein bottle admits a nonwhere vanishing vector field. This contradicts Corollary~\ref{cor:characteristic-klein-bottle}.
\end{proof}

\begin{remark}
    In symplectic Kodaira dimension zero, torus bundles over torus are examples with Euler characteristic zero, which can admit nullhomologous Lagrangian Klein bottle as indicated by \cite[\S 1.6]{Shevchishin2009}.
\end{remark}

\subsection{Euler characteristic zero: infinite KB Luttinger counts}
We now show that a fixed extended Auroux data can support infinitely many KB Luttinger equivalence classes by using symplectic Calabi--Yau surfaces with Euler characteristic zero.

\begin{proposition}
\label{prop:SCY-infinite-KB-count}
For every integer $q\geq1$ and $b\in\{2,3\}$,
\[
\mathcal N_{(0,0,0,2q,b)}=\infty.
\]
Moreover, every integral symplectic manifold realizing these data is symplectic Calabi--Yau.
\end{proposition}

\begin{proof}
Let $(X,\omega)$ realize $(0,0,0,2q,b)$. If $c_1$ were non-torsion, a sufficiently small perturbation $\omega_\varepsilon$ of $\omega$ would satisfy $c_1\cdot[\omega_\varepsilon]>0$ without changing $c_1$. The Liu--Ohta--Ono theorem \cite{Liu1996,OhtaOno1996} would then imply that $X$ is rational or ruled, and hence that $b_2^+(X)=1$. However, it then follows from the light cone lemma that $c_1$ is proportional to $[\omega]$, which contradicts $[\omega]^2>0$. Thus $(X,\omega)$ must be symplectic Calabi--Yau.

Now fix $q\geq1$. For each odd prime $p$ and $b\in\{2,3\}$, let $M_{p,b}$ be the mapping torus of the linear self-diffeomorphism of $T^2$ given by the matrix $A_{p,b}$, where
\[A_{p,2}=
\begin{pmatrix}
p+2&-1\\
1&0
\end{pmatrix},
\qquad
A_{p,3}=
\begin{pmatrix}
1&p\\
0&1
\end{pmatrix}.\]
Set \(X_{p,b}=S^1\times M_{p,b}.\) In the notation of \cite{Geiges1992}, the natural $T^2$-bundle $X_{p,b}\to T^2$ has the bundle data \(\bigr\{A_{p,b},I,(0,0)\bigr\},\) where the Euler class is $(0,0)$ since the origin in each torus fiber defines a section. For $b=2$, the monodromy has $\operatorname{tr}(A_{p,2})>2$ so that $X_{p,2}$ belongs to \cite[Table~1(g)]{Geiges1992}. For $b=3$, the manifold $X_{p,3}$ is the primary Kodaira surface $W_p$ in Example \ref{example:torus-betti-change}. The Wang sequence and K\"unneth formula give
\[H_1(X_{p,b};\ZZ)\cong\ZZ^2\oplus\operatorname{coker}(A_{p,b}-I)\cong\ZZ^b\oplus\ZZ_p,\]
and the intersection form
\[H^2(X_{p,b};\ZZ)/\operatorname{Tor}\cong U^{\oplus (b-1)}\] where $U$ is the hyperbolic lattice. Choose integral classes $u,v$ in one hyperbolic summand so that \(u^2=v^2=0,\, u\cdot v=1.\) Then $\alpha_q=u+qv$ is integral and satisfies $\alpha_q^2=2q$. By Geiges's description of the symplectic cones of torus bundles \cite[Theorem 2]{Geiges1992}, $\alpha_q$ is represented by an integral symplectic form $\omega_{p,b,q}$. Symplectic Kodaira dimension only depends on the underlying smooth four-manifold \cite{HoLi2012}, so $c_1(\omega_{p,b,q})$ must be torsion. Together with $\chi(X_{p,b})=0$, this gives
\[\widetilde{\mathcal A}(X_{p,b},\omega_{p,b,q})=(0,0,0,2q,b).\]

Finally, we claim that distinct primes give distinct KB Luttinger equivalence classes. To see this, let $X'$ be obtained from $X_{p,b}$ by surgery along a Lagrangian Klein bottle $L$. Maintain the notions in Section \ref{section:torus-KB-comparison} and take the dyadic ring $R=\ZZ[1/2]$ in \eqref{equ:Excision and Poincare dual}. We have
\[H_i(X_{p,b},C;R)\cong H^{4-i}(L;R)=0,\qquad i=1,2,\]
because $H^3(L;R)=0$ and $H^2(L;R)=R/2R=0$. The same argument applies to $X'$, so the long exact sequences \eqref{equ:LES} of the pairs yield
\[H_1(X_{p,b};R)\cong H_1(C;R)\cong H_1(X';R).\]
In particular, the torsion subgroup with odd order of $H_1(-;\ZZ)$ is unchanged under KB Luttinger surgery. This subgroup is $\ZZ_p$ for $X_{p,b}$, thus proving the theorem.
\end{proof}

\section{KB Luttinger surgery in positive Kodaira dimensions}
\label{sec:KB-positive-Kodaira}

In this section, we give examples of Klein bottle Luttinger surgeries which change symplectic Kodaira dimension from one to two.

\begin{theorem}
\label{thm:KB-change-Kodaira}
For every pair of integers $m>n>0$ and $r\in\{2,3\}$, there exist closed integral symplectic four-manifolds $(X_1,\omega_1)$ and $(X_2,\omega_2)$, related by a single Klein bottle Luttinger surgery, such that $(X_1,\omega_1)$ is minimal with Kodaira dimension one, $(X_2,\omega_2)$ is non-minimal with Kodaira dimension two, and have extended Auroux data
\[\widetilde{\mathcal A}(X_1,\omega_1)=\widetilde{\mathcal A}(X_2,\omega_2)=(0,48,-2rm,4rmn,0).\]
\end{theorem}

\subsection{Symplectic cyclic branched covers}
\label{subsec:KB-branched-cover-formulas}

Let $(X,\omega)$ be a closed symplectic four-manifold and let $B\subset X$ be an embedded symplectic surface with $\operatorname{PD}[B]=rA$, where $r\geq2$ and $A\in H^2(X;\mathbb Z)$. Choose a complex line bundle $\mathcal L\to X$ with $c_1(\mathcal L)=A$ and a transverse section $s$ of $\mathcal L^{\otimes r}$ with zero set $B$. Then,
\begin{equation}\label{eq:cyclic-construction}
    \widehat X=\bigl\{(x,w):w\in\mathcal L_x,\;w^{\otimes r}=s(x)\bigr\}\subset\mathcal L
\end{equation} 
defines a cyclic branched cover $\pi: \widehat X\to X$ with ramification of order $r$ along $B$ (\cite[\S 7.1]{GSbook}). One can either use symplectic cuts (\cite{Gompf1998}) or a local perturbation (\cite[\S 6]{Auroux2000}, \cite[\S 2.5]{Pasquotto2004}) to construct a symplectic form $\widehat\omega$ on $\widehat X$ which agrees with $\pi^*\omega$ away from an arbitrarily small neighborhood of $\pi^{-1}(B)$ and such that
\[[\widehat\omega]=\pi^*[\omega].\]
In particular, a Lagrangian submanifold $L\subset (X,\omega)$ disjoint from a neighborhood of $B$ has Lagrangian inverse image $\widehat L=f^{-1}(L)\subset (\widehat X,\widehat\omega)$.

The Auroux data of $(\widehat X,\widehat\omega)$ can be computed from some well-known formulas about the symplectic branched covers, which can be found in \cite{GSbook,Gompf1998,Hamilton2012}:
\begin{equation}
\label{eq:branch-cover-change}
\begin{aligned}
c_1(\widehat X)&=\pi^*\bigl(c_1(X)+(1-r)A\bigr),\\
c_1(\widehat X)^2&=r\bigl(c_1(X)+(1-r)A\bigr)^2,\\
c_2(\widehat X)&=rc_2(X)-(r-1)\chi(B),\\
c_1(\widehat X)\cdot[\widehat\omega]&=r\bigl(c_1(X)+(1-r)A\bigr)\cdot[\omega],\\
[\widehat\omega]^2&=r[\omega]^2.\\
\end{aligned}
\end{equation}

\subsection{A criterion for changing Kodaira dimension}
\label{subsec:KB-general-lifting}
\begin{lemma}
\label{lem:KB-zero-canonical-square}
Let $(X,\omega)$ be a closed symplectic four-manifold with \(c_1(X)^2=0,\,c_1(X)\cdot[\omega]<0.\) Then the symplectic Kodaira dimension $\kappa(X,\omega)=1$ if $(X,\omega)$ is minimal; $\kappa(X,\omega)=2$ if it is non-minimal.
\end{lemma}

\begin{proof}
The minimal case follows from the definition. In the non-minimal case, the minimal model $(X_{\min},\omega_{\min})$ has $c_1(X_{\min})^2>0$. Hence, it must have Kodaira dimension either $-\infty$ or $2$. If its Kodaira dimension were $-\infty$, it would be rational or ruled. $c_1(X_{\min})^2>0$ would force it to be rational, so $b_2^+(X)=1$ and \(c_1(X_{\min})\cdot[\omega_{\min}]>0.\)
Thus $c_1(X_{\min})$ and $[\omega_{\min}]$ lie in the same component of the positive cone. By the blowup formula, so do $c_1(X)$ and $[\omega]$. Since $c_1(X)^2=0$, we obtain $c_1(X)\cdot[\omega]\geq 0$ by the light cone lemma, a contradiction. Hence the Kodaira dimension must be $2$.
\end{proof}
Our setup to apply Lemma \ref{lem:KB-zero-canonical-square} is as follows. Let $L_Y\subset(Y,\omega_Y)$ be a Lagrangian Klein bottle with a Weinstein neighborhood $U$. Suppose that \(\pi_0:(X_0,\omega_0)\rightarrow(Y,\omega_Y)\) is a symplectic cyclic branched cover of degree $r\geq2$, whose branch locus is disjoint from $U$, and that $\pi^{-1}_0(U)\to U$ is the trivial $r$-fold cover. Write \(\pi_0^{-1}(L_Y)=L_1\sqcup\cdots\sqcup L_r.\)  Perform Luttinger surgeries on $L_1,\ldots,L_{r-1}$ to obtain $(X_2,\omega_2)$, and then the surgery on $L_r$ to obtain $(X_1,\omega_1)$.

\begin{proposition}
\label{prop:KB-lifting-criterion}
Assume that
\begin{enumerate}[nosep]
\item $c_1(X_1)=kx$ for some $x\in H^2(X_1;\mathbb Z)$ and $k\geq2$, with \(x^2=0,\, x\cdot[\omega_1]<0;\)
\item $X_0$ contains a symplectic exceptional sphere disjoint from $L_1,\ldots,L_{r-1}$.
\end{enumerate}
Then $(X_1,\omega_1)$ is minimal with Kodaira dimension one, and $(X_2,\omega_2)$ is non-minimal with Kodaira dimension two.
\end{proposition}

\begin{proof}
An embedded symplectic exceptional sphere has first Chern pairing equal to $1$ by adjunction formula. The divisibility in the first assumption therefore implies that $(X_1,\omega_1)$ is minimal. It also gives the assumptions in Lemma~\ref{lem:KB-zero-canonical-square}, which implies that
$\kappa(X_1,\omega_1)=1$.

Shrink the surgery neighborhoods so that the exceptional sphere in the second assumption is disjoint from those of $L_1,\ldots,L_{r-1}$. It survives as a symplectic exceptional sphere in $X_2$, so $(X_2,\omega_2)$ is non-minimal. Corollary~\ref{cor:surgery-symplectic-pairings} gives \(c_1(X_2)^2=c_1(X_1)^2=0,\,c_1(X_2)\cdot[\omega_2]=c_1(X_1)\cdot[\omega_1]<0.\) The same lemma now gives $\kappa(X_2,\omega_2)=2$.
\end{proof}

The first assumption in Proposition \ref{prop:KB-lifting-criterion} can be checked downstairs. Let $(Z,\omega_Z)$ be the result of the chosen surgery on $L_Y$. Performing all $r$ lifting surgeries produces a cyclic branched cover \(\pi_1:(X_1,\omega_1)\rightarrow(Z,\omega_Z).\) These operations fit into the diagram
\begin{equation*}
\label{eq:KB-lifted-surgery-diagram}
\begin{tikzcd}[column sep=large, row sep=large]
(X_0,\omega_0)
  \arrow[d, "\pi_0"']
  \arrow[r, dashed, "{L_1,\ldots,L_{r-1}}"]
&
(X_2,\omega_2)
  \arrow[r, dashed, "L_r"]
&
(X_1,\omega_1)
  \arrow[d, "\pi_1"]
\\
(Y,\omega_Y)
  \arrow[rr, dashed, "L_Y"']
&&
(Z,\omega_Z).
\end{tikzcd}
\end{equation*}
Write $B_Z\subset(Z,\omega_Z)$ for the branch surface corresponding to $B\subset (Y,\omega_Y)$ under the identification of the common complement of the Luttinger surgery. Then \(\operatorname{PD}[B_Z]=rA_Z\) for some $A_Z\in H^2(Z;\ZZ)$ and \eqref{eq:branch-cover-change} gives
\[c_1(X_1)=\pi_1^*\bigl(c_1(Z)+(1-r)A_Z\bigr),\qquad[\omega_1]=\pi_1^*[\omega_Z].\]
Consequently, it suffices to find $z\in H^2(Z;\mathbb Z)$ and $k\geq2$ such that
\begin{equation}\label{eq:check-downstairs}
    c_1(Z)+(1-r)A_Z=kz,\qquad z^2=0,\qquad z\cdot[\omega_Z]<0.
\end{equation}

\subsection{The constructions using octic and sextic curves}
\label{subsec:KB-cyclic-plane-construction}

\begin{proof}[Proof of Theorem~\ref{thm:KB-change-Kodaira}]
Let \(Y:=\operatorname{Bl}_{p}\mathbb{CP}^2\cong\CP^2\#\overline{\CP}^2\) be the complex blowup at the point $p=[1:0:0].$ Denote the line class by $H$, and use $E$ for the exceptional curve or its class. We can choose a K\"ahler form $\omega_Y$ with \([\omega_Y]=(m+n)H-(m-n)E\) for which complex conjugation is anti-symplectic by averaging a K\"ahler representative with the negative of its conjugate pullback. The real locus \(L_Y\) is the real blowup of $\RR\PP^2\subset \CP^2$, which is a Lagrangian Klein bottle. By Section \ref{sec:rational surfaces}, the Luttinger surgery on $L_Y$ gives $(Z,\omega_Z)$ where $Z\cong S^2\times S^2$ with \([\omega_Z]=2n\alpha+m\beta,\, c_1(Z)=2\alpha+2\beta,\) for the two $S^2$-factor classes $\alpha,\beta\in H^2(Z;\ZZ)$.

For $d=r\ell$ with $r\geq 2$ and $\ell>0$, take the proper transform of the degree $d$ curve $F_d:=\{z_0^d+z_1^d+z_2^d=0\}\subset \CP^2$ as the branch locus $B\subset(Y,\omega_Y)$. Note that when $d$ is even, $F_d$ avoids $\mathbb{RP}^2\subset \CP^2$ so that $B$ avoids the Klein bottle $L_Y$. We have the $r$-fold cyclic branched cover $\pi_0:(X_0,\omega_0)\to (Y,\omega_Y)$ defined by the line bundle $\mathcal L$ over $Y$ with $c_1(\mathcal{L})=\ell H$ using the construction \eqref{eq:cyclic-construction}. 

We claim that when $\ell$ is even, $\pi_0^{-1}(L_Y)$ is a disjoint union of $r$ Lagrangian Klein bottles. Let $b:Y\to\mathbb{CP}^2$ denote the complex blowdown map, so that $\mathcal L=b^*\mathcal O(\ell)$ and $b(L_Y)=\mathbb{RP}^2$. Write $z_0^\ell+z_1^\ell+z_2^\ell$ for a section of $\mathcal{O}(\ell)$. The pullback section \(\tau=b^*(z_0^\ell+z_1^\ell+z_2^\ell)\) is nowhere vanishing along $L_Y$ and hence trivializes $\mathcal L|_{L_Y}$. In the induced trivialization of $\mathcal L^{\otimes r}|_{L_Y}$, the defining section of the cyclic cover is $f\tau^{\otimes r}$, where
\[f(y)=\frac{z_0^d+z_1^d+z_2^d}{(z_0^\ell+z_1^\ell+z_2^\ell)^r}>0,\qquad \text{for }b(y)=[z_0:z_1:z_2]\in\mathbb{RP}^2.\]
It follows that $f$ has globally defined $r$th roots on $L_Y$ so that the covering equation in \eqref{eq:cyclic-construction} has exactly $r$ solutions $\zeta f^{1/r}\tau$, indexed by $\zeta^r=1$. Their graphs give \(\pi_0^{-1}(L_Y)=L_1\sqcup\cdots\sqcup L_r,\) with each restriction $\pi_0|_{L_i}:L_i\to L_Y$ a diffeomorphism. Moreover, each $L_i$ is Lagrangian because $\omega_0=\pi_0^*\omega_Y$ near $\pi_0^{-1}(L_Y)$.

Next, in order to meet the requirements in \eqref{eq:check-downstairs}, we need $z^2=0,z\cdot[\omega_Z]<0$ and $k\geq 2$ such that
\[c_1(Z)+(1-r)A_Z=2\alpha+2\beta+(1-r)\frac{\ell}{2}(2\alpha+\beta)=kz.\]
It follows that we must have 
\[z=-\alpha,\qquad-k=(1-r)\ell+2,\qquad (1-r)\frac{\ell}{2}+2=0\]
So, \((r-1)\ell=4\) gives that the possibilities are
\begin{equation}
\label{eq:KB-two-cyclic-choices}
(r,\ell,d)=(2,4,8)\quad\text{or}\quad(3,2,6).
\end{equation}
Hence, for either choice in \eqref{eq:KB-two-cyclic-choices}, the first assumption in Proposition \ref{prop:KB-lifting-criterion} is satisfied. To check the second, note that \(\pi_0^{-1}(E)=E_1\sqcup\cdots\sqcup E_r\) is a union of symplectic exceptional spheres since $p\notin F_d$ and $\mathcal{L}|_E$ is trivial. Label the components so that $E_j$ meets $L_j$ along the lifted real exceptional circle. Then $E_r$ is disjoint from $L_1,\ldots,L_{r-1}$. Therefore, Proposition~\ref{prop:KB-lifting-criterion} proves the assertions about minimality and Kodaira dimension.

Finally, let us compute the extended Auroux data for the two choices in \eqref{eq:KB-two-cyclic-choices}. Notice that
\(\chi(B)=\chi(F_d)=3d-d^2.\)
Since $c_2(Y)=4$, \eqref{eq:branch-cover-change} gives
\[c_2(X_0)=4r-(r-1)(3d-d^2)=48.\]
Also,
\[\begin{aligned}
c_1(X_0)^2&=r\bigr(3H-E+(1-r)\ell H\bigr)^2=0,\\
c_1(X_0)\cdot[\omega_0]&=r\bigr(3H-E+(1-r)\ell H\bigr)\cdot\bigr((m+n)H-(m-n)E\bigr)=-2rm,\\
[\omega_0]^2&=r\bigr((m+n)H-(m-n)E\bigr)^2=r\bigl((m+n)^2-(m-n)^2\bigr)=4rmn.
\end{aligned}\]
Moreover, $b_1(X_0)=0$ follows from \cite[\S 3]{Rohlin71}. Corollary~\ref{cor:surgery-symplectic-pairings} and Lemma~\ref{lemma:surgery-betti-change} now give $\widetilde{\mathcal{A}}$ for both $(X_1,\omega_1)$ and $(X_2,\omega_2)$.
\end{proof}

\begin{remark}\label{rmk:cyclic-origination}
The manifolds in Theorem \ref{thm:KB-change-Kodaira} also provide counterexamples with $b_2^+>1$ to Auroux's original question concerning torus Luttinger surgery, since they have the same Auroux data but distinct symplectic Kodaira dimensions so that \cite{HoLi2012} works. We opt not to include it in Section \ref{sec:Auroux-question} since we have adapted the construction here to show that KB Luttinger surgery can change symplectic Kodaira dimension.
\end{remark}

\begin{appendices}
\section{Sequential and disjoint Luttinger surgery}\label{section:two notions}

In this appendix, we distinguish the two surgery relations appearing in Auroux's original uniqueness question \ref{question:Auroux} and in the formulation of McDuff--Salamon \cite[Problem 11]{McDuffSalamon2017}. The distinction concerns whether the surgery tori may be chosen successively in the intermediate manifolds or must form a fixed pairwise disjoint collection in the initial manifold.

\begin{definition}
\label{def:two-Luttinger-equivalences}
Two symplectic four-manifolds are {\bf sequentially Luttinger related} if there is a finite chain in which each term is obtained from the preceding one by a torus Luttinger surgery. They are {\bf disjointly Luttinger related} if there is a fixed pairwise disjoint collection of Lagrangian tori in the initial manifold such that surgeries on their disjoint neighborhoods produce the second manifold.
\end{definition}

Disjointly Luttinger related symplectic manifolds must also be sequentially Luttinger related. The following Proposition \ref{prop:hyperelliptic-sequential-not-disjoint} shows that the converse fails, even assuming the underlying smooth four-manifolds are the same. 

Recall that a hyperelliptic surface, also called a bielliptic surface, admits a presentation \((E\times F)/G,\) where $E,F$ are elliptic curves and $G$ acts by translations on $E$ and faithfully by automorphisms on $F$, with $F/G\cong\mathbb{CP}^1$. The classification of them gives seven types, corresponding to
\[G\cong\mathbb Z_2,\quad(\mathbb Z_2)^2,\quad\mathbb Z_4,\quad\mathbb Z_4\times\mathbb Z_2,\quad\mathbb Z_3,\quad(\mathbb Z_3)^2,\quad\mathbb Z_6;\]
see \cite[V.5]{BHPV}. For the $\mathbb Z_6$-type, one may take $F=\mathbb C/(\mathbb Z+\zeta\mathbb Z)$, where $\zeta=e^{\pi i/3}$, with a generator of $G$ acting by \((x,y)\mapsto(x+a,\zeta y),\) where $a\in E$ has exact order six. Denote by $X$ the underlying smooth four-manifold. Its intersection form $\Lambda$ is the rank $2$ hyperbolic lattice with a basis $s_1,s_2$ such that $s_1^2=s_2^2=0,\,s_1\cdot s_2=1$. Consider two integral symplectic forms $\omega_0,\,\omega_1$ such that \([\omega_0]=6s_1+s_2,\,[\omega_1]=2s_1+3s_2.\) Similar to the Enriques surface Example \ref{example:Auroux-Enriques}, the existence is again ensured by the ampleness criterion since it is a standard fact that generic complex structures on $X$ do not admit negative curves.

We now explain why they are not disjointly Luttinger related by modifying the proof of Corollary \ref{cor:torus-period-lattice}. Suppose that surgeries on a disjoint collection of Lagrangian tori in $(X,\omega_0)$ produced a symplectic manifold symplectomorphic to $(X,\omega_1)$. Since $b_2^+(X)=1$ and $H_2(X;\ZZ)$ has no torsion (\cite{Suwa}), all Lagrangian tori are nullhomologous. The relative homology sequence gives a surjection \(H_2(C;\mathbb Z)\rightarrow H_2(X;\mathbb Z)\). Choose lifts of an integral basis from $H_2(X;\ZZ)$ to $H_2(C;\ZZ)$. Their images in the surgered manifold have the same intersection numbers and symplectic periods. Hence, we have an integral lattice isomorphism which preserves the symplectic periods. This is the argument underlying Corollary~\ref{cor:torus-period-lattice}, but applied to multiple disjoint Lagrangian tori. Consequently, $[\omega_0]$ and $[\omega_1]$ would lie in the same $O(\Lambda)$-orbit. However, as in Example~\ref{example:Auroux-Enriques}, our choice the classes $6s_1+s_2,2s_1+3s_2$ implies that they must lie in distinct $O(\Lambda)$-orbits. This leads to a contradiction.

Next, we give a precise construction of $(X_0,\omega_0)$ and $(X_1,\omega_1)$ from Luttinger surgeries on two distinct Lagrangian tori $L_0,L_1$ respectively in a Kodaira--Thurston manifold $(W,\Omega_W).$ As a consequence, $(X_0,\omega_0)$ and $(X_1,\omega_1)$ are related by a two-step sequence of Luttinger surgeries. Notice that the intermediate manifold $W$ has $b_2^+(W)=2$, so the lattice obstruction above does not apply at every stage of this sequence.

\begin{proposition}
\label{prop:hyperelliptic-sequential-not-disjoint}
Let $X$ be the underlying smooth four-manifold of a hyperelliptic surface of $\mathbb Z_6$-type. There exist integral symplectic forms $\omega_0,\omega_1$ on $X$ with
\[\mathcal{A}(X,\omega_0)=\mathcal{A}(X,\omega_1)=(0,0,0,12)\]
such that $(X,\omega_0)$ and $(X,\omega_1)$ are sequentially Luttinger related but are not disjointly Luttinger related.
\end{proposition}

\begin{proof} 
\noindent{\bf Step 1. The construction.}
Give \(T^4=(\mathbb R/\mathbb Z)^4\) coordinates $(x_1,y_1,x_2,y_2)$ and the symplectic form \[\Omega=dx_1\wedge dy_1+6\,dx_2\wedge dy_2.\]
Consider the following three embedded Lagrangian tori
\[L=\{y_1=y_2=0\},\,L_0=\{x_1=0,\ y_2=1/3\},\,L_1=\{x_1-2x_2=0,\ y_1+3y_2=1/3\}.\]
So, $L$ is disjoint from $L_0\cup L_1$, while $L_0$ and $L_1$ have two intersection points; see Figure \ref{fig:hyperelliptic-three-tori}.
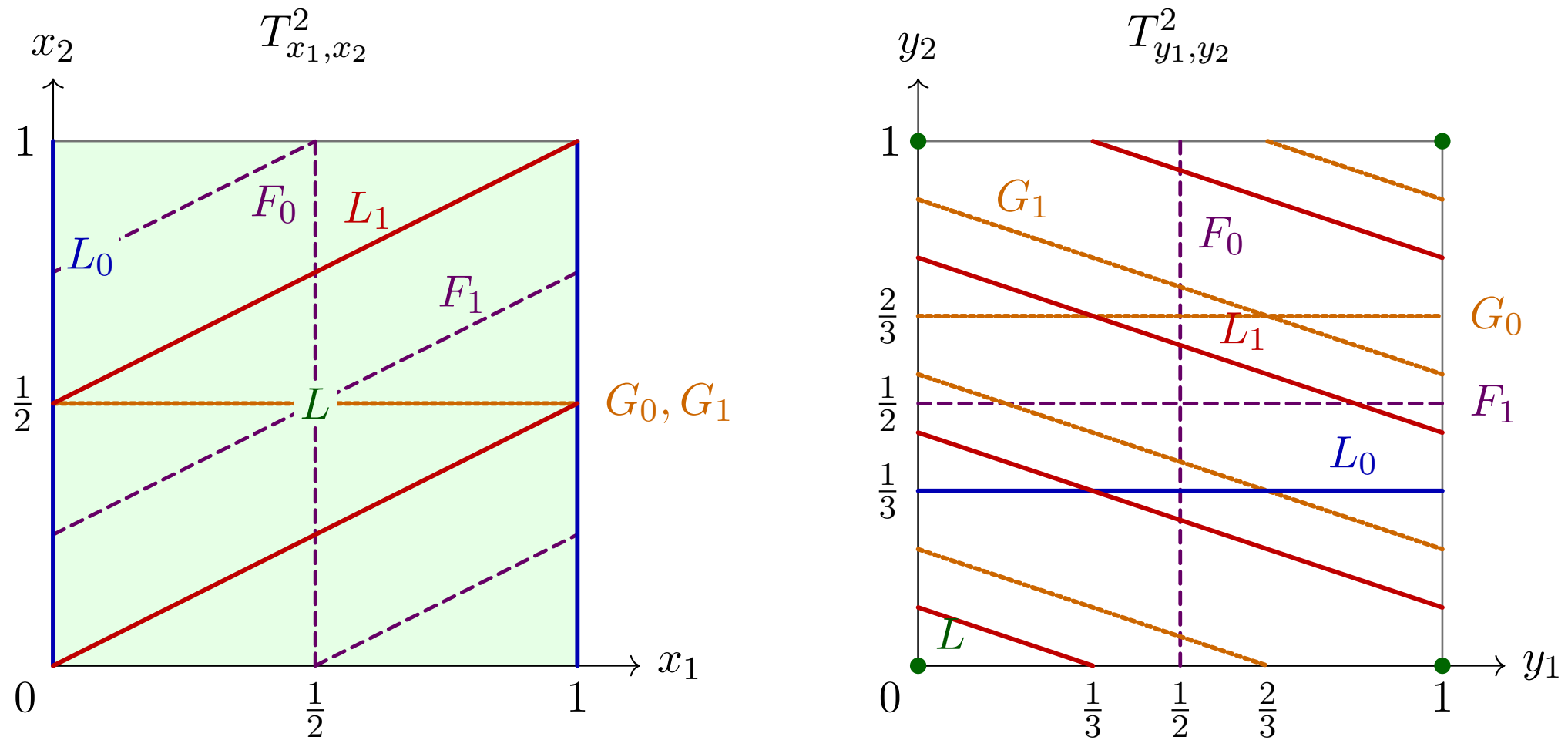
\begin{figure}[htbp]
\centering
\begin{tikzpicture}[
  x=4.5cm,y=4.5cm,font=\small,
  line cap=round,line join=round,
  zero/.style={blue!70!black,line width=1.1pt},
  one/.style={red!75!black,line width=1.1pt},
  sympF/.style={violet!80!black,line width=0.9pt,densely dashed},
  sympG/.style={orange!80!black,line width=1.1pt,densely dotted},
  boundary/.style={black!55,line width=0.5pt}
]

\begin{scope}
  \fill[green!10] (0,0) rectangle (1,1);
  \draw[boundary] (0,0) rectangle (1,1);

  \draw[->] (0,0) -- (1.12,0) node[right] {$x_1$};
  \draw[->] (0,0) -- (0,1.12) node[above] {$x_2$};

  \draw[sympF] ({1/2},0) -- ({1/2},1);
  \draw[sympF] ({1/2},0) -- (1,{1/4});
  \draw[sympF] (0,{1/4}) -- (1,{3/4});
  \draw[sympF] (0,{3/4}) -- ({1/2},1);
  \draw[sympG] (0,{1/2}) -- (1,{1/2});

  \draw[zero] (0,0) -- (0,1);
  \draw[zero] (1,0) -- (1,1);

  \draw[one] (0,0) -- (1,{1/2});
  \draw[one] (0,{1/2}) -- (1,1);

  \node[below left] at (0,0) {$0$};
  \node[below] at (1,0) {$1$};
  \node[left] at (0,1) {$1$};
  \node[left] at (0,{1/2}) {$\tfrac12$};
  \node[below] at ({1/2},0) {$\tfrac12$};

  \node[blue!70!black,right,fill=green!10,inner sep=1pt,outer sep=2pt]
    at (0,0.78) {$L_0$};
  \node[red!75!black,above] at (0.60,0.80) {$L_1$};
  \node[green!35!black,fill=green!10,inner sep=1.5pt]
    at (0.5,0.5) {$L$};

  \node[violet!80!black,left] at (0.5,0.88) {$F_0$};
  \node[violet!80!black,above] at (0.78,0.64) {$F_1$};
  \node[orange!80!black,right] at (1.02,0.5) {$G_0,G_1$};

  \node at (0.5,1.20) {$T^2_{x_1,x_2}$};
\end{scope}

\begin{scope}[shift={(1.65,0)}]
  \draw[boundary] (0,0) rectangle (1,1);

  \draw[->] (0,0) -- (1.12,0) node[right] {$y_1$};
  \draw[->] (0,0) -- (0,1.12) node[above] {$y_2$};

  \draw[sympF] ({1/2},0) -- ({1/2},1);
  \draw[sympF] (0,{1/2}) -- (1,{1/2});

  \draw[sympG] (0,{2/3}) -- (1,{2/3});

  \draw[sympG] (0,{2/9}) -- ({2/3},0);
  \draw[sympG] (0,{5/9}) -- (1,{2/9});
  \draw[sympG] (0,{8/9}) -- (1,{5/9});
  \draw[sympG] ({2/3},1) -- (1,{8/9});

  \draw[zero] (0,{1/3}) -- (1,{1/3});

  \draw[one] (0,{1/9}) -- ({1/3},0);
  \draw[one] (0,{4/9}) -- (1,{1/9});
  \draw[one] (0,{7/9}) -- (1,{4/9});
  \draw[one] ({1/3},1) -- (1,{7/9});

  \foreach \p in {(0,0),(1,0),(0,1),(1,1)}
    \fill[green!40!black] \p circle (2pt);

  \node[below left] at (0,0) {$0$};
  \node[below] at (1,0) {$1$};
  \node[left] at (0,1) {$1$};
  \node[below] at ({1/3},0) {$\tfrac13$};
  \node[left] at (0,{1/3}) {$\tfrac13$};
  \node[below] at ({1/2},0) {$\tfrac12$};
  \node[below] at ({2/3},0) {$\tfrac23$};
  \node[left] at (0,{1/2}) {$\tfrac12$};
  \node[left] at (0,{2/3}) {$\tfrac23$};

  \node[blue!70!black,above] at (0.83,{1/3}) {$L_0$};
  \node[red!75!black,above] at (0.62,0.57) {$L_1$};
  \node[green!35!black,above right] at (0,0) {$L$};

  \node[violet!80!black,right] at (0.5,0.82) {$F_0$};
  \node[violet!80!black,right] at (1.02,{1/2}) {$F_1$};
  \node[orange!80!black,right] at (1.02,{2/3}) {$G_0$};
  \node[orange!80!black,above] at (0.20,{8/9-0.20/3}) {$G_1$};

  \node at (0.5,1.20) {$T^2_{y_1,y_2}$};
\end{scope}
\end{tikzpicture}
\caption{
Coordinate projections of the three Lagrangian tori $L,L_0,L_1$ and the four symplectic tori $F_0,F_1,G_0,G_1$. The tori $L_0$ and $L_1$ are the products of the blue and red circles respectively. The torus $L$ is the product of the shaded $(x_1,x_2)$-torus with the origin of the $(y_1,y_2)$-torus. The symplectic tori $F_i$ and $G_i$ are the products of the correspondingly labelled purple dashed and orange dotted circles respectively. These four symplectic tori will be used in Step~3 to determine the periods of $\omega_0$ and $\omega_1$.}
\label{fig:hyperelliptic-three-tori}
\end{figure}

Write $a_1,b_1,a_2,b_2$ for the oriented coordinate directions corresponding to $x_1,y_1,x_2,y_2$, and set \(u=2a_1+a_2,\, v=3b_1-b_2\) to be the directions of the red circles in Figure \ref{fig:hyperelliptic-three-tori}. We use the affine translation framings, which are exactly the Lagrangian framings since nearby parallel tori are still Lagrangian. Perform the $(L,a_1,1)$ Luttinger surgery to obtain $(W,\Omega_W)$, using a neighborhood disjoint from $L_0\cup L_1$. Thus $W$ is diffeomorphic to the Kodaira--Thurston manifold by Example~\ref{example:torus-betti-change}. Then, let $(X_0,\omega_0)$ and $(X_1,\omega_1)$ be obtained from $(W,\Omega_W)$ by the Luttinger surgeries using \((L_0,b_1,1)\) and \((L_1,v,1)\) respectively.  Let $L_0'$ be the Lagrangian torus in $X_0$ corresponding to $L_0\subset W$ after the surgery. Hence, $(X_1,\omega_1)$ and $(X_0,\omega_0)$ are sequentially Luttinger related by performing one surgery on $L_0'$ followed by another one on $L_1$. The construction is summarized in the following diagram. 
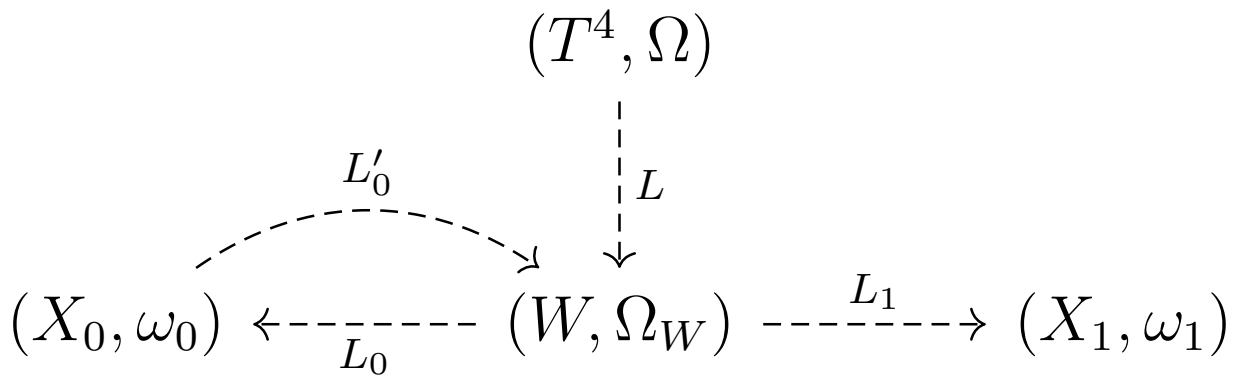
\begin{equation*}
\begin{tikzcd}[column sep=large, row sep=large]
&
(T^4,\Omega)
  \arrow[d, dashed, "L"]
&
\\
(X_0,\omega_0)
  \arrow[r, dashed, bend left=35, "L_0'"]
&
(W,\Omega_W)
  \arrow[l, dashed, "L_0"]
  \arrow[r, dashed, "L_1"]
&
(X_1,\omega_1)
\end{tikzcd}
\end{equation*}

\smallskip
\noindent {\bf Step 2. Identify the underlying manifolds.}
Consider the diffeomorphism $f$ of $T^4$ induced by the map in $\text{SL}(4,\ZZ)$ given by
\[F(x_1,y_1,x_2,y_2)=(x_1+2x_2,\ 3y_1+y_2,\ x_2,\ -y_1).\] It satisfies \( f(L)=L,\, f(L_0)=L_1,\)
and sends
\[a_1\longmapsto a_1,\qquad a_2\longmapsto u,\qquad b_1\longmapsto v,\qquad b_2\longmapsto b_1.\]
The map $f$ therefore naturally induces a diffeomorphism between $X_0$ and $X_1$. It therefore suffices to identify $X_0$ with $X$.

By Example~\ref{example:torus-betti-change}, we have \(W\cong S^1_{x_2}\times M_{D_{x_1}},\) where $M_{D_{x_1}}$ is the mapping torus of the Dehn twist along the $x_1$-circle with base $S^1_{y_2}$. The torus $L_0$ is the product of $S^1_{x_2}$ with the $y_1$-circle in the torus fiber over $y_2=1/3$. Thus, the second Luttinger surgery on $L_0$ inserts a  Dehn twist $D_{y_1}$ about the $y_1$-circle into the monodromy. Consequently, \(X_0\cong S^1_{x_2}\times M_A,\) where $M_A$ is the mapping torus with monodromy
\[A=D_{x_1}D_{y_1}=\begin{pmatrix}1&1\\0&1\end{pmatrix}\begin{pmatrix}1&0\\-1&1\end{pmatrix}=\begin{pmatrix}0&1\\-1&1\end{pmatrix}.\]
Reversing both the mapping torus base circle and the $S^1_{x_2}$-factor, $S^1_{x_2}\times M_A$ is oriented diffeomorphic to $S^1_{-x_2}\times M_{A^{-1}}$. The bundle data of $X_0$ in the notation of \cite{Geiges1992} are therefore
\[\Bigl\{\begin{pmatrix}1&-1\\1&0\end{pmatrix},I,(0,0)\Bigr\},\]
which is the case in \cite[Table~1(c)]{Geiges1992}, whose total space is a hyperelliptic surface. To identify its hyperelliptic type, note that since $A^3=-I,\,A^6=I$, the $6$-fold cover of the mapping torus base circle $S^1_{y_2}$ gives rise to a product of two tori $E\times F$ as the covering space of $S^1_{x_2}\times M_A$ with deck group $\ZZ_6$. Thus $X_0$ is of cyclic $\mathbb Z_6$-type.

\smallskip
\noindent {\bf Step 3. Determine the symplectic periods.}
In the original $T^4$, choose symplectic tori
\[F_0=\{x_1=y_1=1/2\}, \qquad G_0=\{x_2=1/2,\ y_2=2/3\}\]
They are disjoint from $L\cup L_0$. Their images under the diffeomorphism $f$ are
\[F_1=\{x_1-2x_2=1/2,\ y_2=1/2\},\qquad G_1=\{x_2=1/2,\ y_1+3y_2=2/3\},\]
which are disjoint from $L\cup L_1$; see Figure \ref{fig:hyperelliptic-three-tori}. The tori $F_0,G_0$ have exactly one positively transverse intersection point, and the same holds for their images under $f$. Their self-intersections are zero. Therefore, $\{[F_i],[G_i]\}$ forms the basis $\{s_1,s_2\}$ of the hyperbolic lattice $\Lambda=H_2(X_i;\mathbb Z)$ such that $s_1^2=s_2^2=0,\,s_1\cdot s_2=1$. Since the symplectic form $\omega_i$ agrees with $\Omega$ near the tori $F_i,G_i$, we can compute
\[\bigl(\omega_0(F_0),\omega_0(G_0)\bigr)=\bigl(\Omega(F_0),\Omega(G_0)\bigr)=(6,1),\quad\bigl(\omega_1(F_1),\omega_1(G_1)\bigr)=\bigl(\Omega(F_1),\Omega(G_1)\bigr)=(2,3).\]
Hence, $[\omega_0],[\omega_1]$ represent $6s_1+s_2,2s_1+3s_2$ in the hyperbolic lattice. In particular, \(\mathcal A(X_i,\omega_i)=(0,0,0,12)\) for $i=1,2$. The obstruction established in the discussion preceding the proposition therefore shows that $(X_0,\omega_0)$ and $(X_1,\omega_1)$ are not disjointly Luttinger related.
\end{proof}

\end{appendices}

\bibliographystyle{alpha}
\bibliography{mybib}

@article {HoLi2012,
    AUTHOR = {Ho, Chung-I and Li, Tian-Jun},
     TITLE = {Luttinger surgery and {K}odaira dimension},
   JOURNAL = {Asian J. Math.},
  FJOURNAL = {Asian Journal of Mathematics},
    VOLUME = {16},
      YEAR = {2012},
    NUMBER = {2},
     PAGES = {299--318},
      ISSN = {1093-6106,1945-0036},
   MRCLASS = {53D35 (57R17 57R65)},
  MRNUMBER = {2916366},
MRREVIEWER = {Federica\ Pasquotto},
       DOI = {10.4310/AJM.2012.v16.n2.a7},
       URL = {https://doi.org/10.4310/AJM.2012.v16.n2.a7},
}

@article {Luttinger1995,
    AUTHOR = {Luttinger, Karl Murad},
     TITLE = {Lagrangian tori in {${\bf R}^4$}},
   JOURNAL = {J. Differential Geom.},
  FJOURNAL = {Journal of Differential Geometry},
    VOLUME = {42},
      YEAR = {1995},
    NUMBER = {2},
     PAGES = {220--228},
      ISSN = {0022-040X,1945-743X},
   MRCLASS = {57R40 (57R95 58F05)},
  MRNUMBER = {1366546},
MRREVIEWER = {Karl\ Friedrich\ Siburg},
       URL = {http://projecteuclid.org/euclid.jdg/1214457232},
}

@incollection {Auroux2005OpenQuestions,
    AUTHOR = {Auroux, Denis},
     TITLE = {Some open questions about symplectic 4-manifolds, singular
              plane curves and braid group factorizations},
 BOOKTITLE = {European {C}ongress of {M}athematics},
     PAGES = {23--40},
 PUBLISHER = {Eur. Math. Soc., Z\"urich},
      YEAR = {2005},
      ISBN = {3-03719-009-4},
   MRCLASS = {53D35 (14H20 20F36 57M25 57N13 57R17)},
  MRNUMBER = {2185734},
MRREVIEWER = {Michael\ J.\ Usher},
}

@article{AdaloglouEvans2026,
  author  = {Nikolas Adaloglou and Jonathan David Evans},
  title   = {Lagrangian {K}lein bottles in ${S}^2 \times {S}^2$},
  JOURNAL = {Comment. Math. Helv.},
  FJOURNAL = {Commentarii Mathematici Helvetici},
  year    = {2026},
  doi     = {10.4171/CMH/620}
}

@article {AurouxDonaldsonKatzarkov2003,
    AUTHOR = {Auroux, D. and Donaldson, S. K. and Katzarkov, L.},
     TITLE = {Luttinger surgery along {L}agrangian tori and non-isotopy for
              singular symplectic plane curves},
   JOURNAL = {Math. Ann.},
  FJOURNAL = {Mathematische Annalen},
    VOLUME = {326},
      YEAR = {2003},
    NUMBER = {1},
     PAGES = {185--203},
      ISSN = {0025-5831,1432-1807},
   MRCLASS = {57R17 (53D35 57R52)},
  MRNUMBER = {1981618},
MRREVIEWER = {Ignasi\ Mundet-Riera},
       DOI = {10.1007/s00208-003-0418-9},
       URL = {https://doi.org/10.1007/s00208-003-0418-9},
}

@article {Nemirovski2009,
    AUTHOR = {Nemirovski\u{i}, S. Yu.},
     TITLE = {The homology class of a {L}agrangian {K}lein bottle},
   JOURNAL = {Izv. Ross. Akad. Nauk Ser. Mat.},
  FJOURNAL = {Izvestiya Rossiiskoi Akademii Nauk. Seriya Matematicheskaya},
    VOLUME = {73},
      YEAR = {2009},
    NUMBER = {4},
     PAGES = {37--48},
      ISSN = {1607-0046,2587-5906},
   MRCLASS = {53D12},
  MRNUMBER = {2583965},
MRREVIEWER = {Alex\ Degtyarev},
       DOI = {10.1070/IM2009v073n04ABEH002462},
       URL = {https://doi.org/10.1070/IM2009v073n04ABEH002462},
}

@book {McDuffSalamon2017,
    AUTHOR = {McDuff, Dusa and Salamon, Dietmar},
     TITLE = {Introduction to symplectic topology},
    SERIES = {Oxford Graduate Texts in Mathematics},
   EDITION = {Third},
 PUBLISHER = {Oxford University Press, Oxford},
      YEAR = {2017},
     PAGES = {xi+623},
      ISBN = {978-0-19-879490-5; 978-0-19-879489-9},
   MRCLASS = {53D35 (53D40 57R17 57R57 57R58)},
  MRNUMBER = {3674984},
MRREVIEWER = {Hansj\"org\ Geiges},
       DOI = {10.1093/oso/9780198794899.001.0001},
       URL = {https://doi.org/10.1093/oso/9780198794899.001.0001},
}

@article {Shevchishin2009,
    AUTHOR = {Shevchishin, V. V.},
     TITLE = {Lagrangian embeddings of the {K}lein bottle and the
              combinatorial properties of mapping class groups},
   JOURNAL = {Izv. Ross. Akad. Nauk Ser. Mat.},
  FJOURNAL = {Izvestiya Rossiiskoi Akademii Nauk. Seriya Matematicheskaya},
    VOLUME = {73},
      YEAR = {2009},
    NUMBER = {4},
     PAGES = {153--224},
      ISSN = {1607-0046,2587-5906},
   MRCLASS = {57R17 (57N05 57R40)},
  MRNUMBER = {2583968},
MRREVIEWER = {Mustafa\ Korkmaz},
       DOI = {10.1070/IM2009v073n04ABEH002465},
       URL = {https://doi.org/10.1070/IM2009v073n04ABEH002465},
}

@article {KK17,
    AUTHOR = {Karshon, Yael and Kessler, Liat},
     TITLE = {Distinguishing symplectic blowups of the complex projective
              plane},
   JOURNAL = {J. Symplectic Geom.},
  FJOURNAL = {The Journal of Symplectic Geometry},
    VOLUME = {15},
      YEAR = {2017},
    NUMBER = {4},
     PAGES = {1089--1128},
      ISSN = {1527-5256,1540-2347},
   MRCLASS = {53D05 (14E05)},
  MRNUMBER = {3734610},
MRREVIEWER = {Weiyi\ Zhang},
       DOI = {10.4310/JSG.2017.v15.n4.a5},
       URL = {https://doi.org/10.4310/JSG.2017.v15.n4.a5},
}

@article {EP94,
    AUTHOR = {Eliashberg, Yakov and Polterovich, Leonid},
     TITLE = {New applications of {L}uttinger's surgery},
   JOURNAL = {Comment. Math. Helv.},
  FJOURNAL = {Commentarii Mathematici Helvetici},
    VOLUME = {69},
      YEAR = {1994},
    NUMBER = {4},
     PAGES = {512--522},
      ISSN = {0010-2571,1420-8946},
   MRCLASS = {58F05 (57R40)},
  MRNUMBER = {1303225},
MRREVIEWER = {Jean-Claude\ Sikorav},
       DOI = {10.1007/BF02564502},
       URL = {https://doi.org/10.1007/BF02564502},
}

@article {Liu1996,
    AUTHOR = {Liu, Ai-Ko},
     TITLE = {Some new applications of general wall crossing formula,
              {G}ompf's conjecture and its applications},
   JOURNAL = {Math. Res. Lett.},
  FJOURNAL = {Mathematical Research Letters},
    VOLUME = {3},
      YEAR = {1996},
    NUMBER = {5},
     PAGES = {569--585},
      ISSN = {1073-2780},
   MRCLASS = {57R55 (57R57)},
  MRNUMBER = {},
MRREVIEWER = {Vicente\ Mu\~noz},
       DOI = {10.4310/MRL.1996.v3.n5.a1},
       URL = {https://doi-org.ezp2.lib.umn.edu/10.4310/MRL.1996.v3.n5.a1},
}

@article {OhtaOno1996,
    AUTHOR = {Ohta, Hiroshi and Ono, Kaoru},
     TITLE = {Notes on symplectic {$4$}-manifolds with {$b^+_2=1$}. {II}},
   JOURNAL = {Internat. J. Math.},
  FJOURNAL = {International Journal of Mathematics},
    VOLUME = {7},
      YEAR = {1996},
    NUMBER = {6},
     PAGES = {755--770},
      ISSN = {0129-167X,1793-6519},
   MRCLASS = {57R57 (57R15 58D10)},
  MRNUMBER = {},
MRREVIEWER = {Zhenbo\ Qin},
       DOI = {10.1142/S0129167X96000402},
       URL = {https://doi-org.ezp2.lib.umn.edu/10.1142/S0129167X96000402},
}

@article{Hamilton2012,
  author   = {Hamilton, M. J. D.},
  title    = {Representing homology classes by symplectic surfaces},
  journal  = {Math. Res. Lett.},
  volume   = {19},
  year     = {2012},
  number   = {5},
  pages    = {1021--1024},
  doi      = {10.4310/MRL.2012.v19.n5.a5},
  mrnumber = {2997055}
}

@article{iida2026,
      title={On {A}uroux's {L}uttinger-surgery uniqueness question}, 
      author={Nobuo Iida},
      year={2026},
      journal={arXiv:2609.29561},
      archivePrefix={arXiv},
      primaryClass={math.SG},
      url={https://arxiv.org/abs/2609.29561}, 
}

@incollection{Dolgachev2016,
  author    = {Dolgachev, Igor V.},
  title     = {A brief introduction to {E}nriques surfaces},
  booktitle = {Development of Moduli Theory---Kyoto 2013},
  series    = {Advanced Studies in Pure Mathematics},
  volume    = {69},
  pages     = {1--32},
  publisher = {Mathematical Society of Japan},
  address   = {Tokyo},
  year      = {2016},
  mrnumber  = {3586505}
}

@article{BaykurSunukjian,
    AUTHOR = {Baykur, R. \.{I}nan\c{c} and Sunukjian, Nathan},
     TITLE = {Round handles, logarithmic transforms and smooth
              4-manifolds},
   JOURNAL = {J. Topol.},
  FJOURNAL = {Journal of Topology},
    VOLUME = {6},
      YEAR = {2013},
    NUMBER = {1},
     PAGES = {49--63},
      ISSN = {1753-8416},
  MRNUMBER = {3029421},
       DOI = {10.1112/jtopol/jts027},
       URL = {https://doi.org/10.1112/jtopol/jts027},
}

@article {Welschinger,
    AUTHOR = {Welschinger, Jean-Yves},
     TITLE = {Effective classes and {L}agrangian tori in symplectic
              four-manifolds},
   JOURNAL = {J. Symplectic Geom.},
  FJOURNAL = {The Journal of Symplectic Geometry},
    VOLUME = {5},
      YEAR = {2007},
    NUMBER = {1},
     PAGES = {9--18},
      ISSN = {1527-5256,1540-2347},
   MRCLASS = {53D35 (53D12 53D45 57R17)},
  MRNUMBER = {2371182},
MRREVIEWER = {Tobias\ Ekholm},
       DOI = {10.4310/jsg.2007.v5.n1.a3},
       URL = {https://doi.org/10.4310/jsg.2007.v5.n1.a3},
}

@article {Smirnov,
    AUTHOR = {Smirnov, Gleb},
     TITLE = {On {L}agrangian tori in {K}3 surfaces},
   JOURNAL = {J. Symplectic Geom.},
  FJOURNAL = {The Journal of Symplectic Geometry},
    VOLUME = {23},
      YEAR = {2025},
    NUMBER = {4},
     PAGES = {835--851},
      ISSN = {1527-5256,1540-2347},
   MRCLASS = {53D12},
  MRNUMBER = {4965277},
MRREVIEWER = {Stefan\ Nemirovski},
       DOI = {10.4310/jsg.250820165748},
       URL = {https://doi.org/10.4310/jsg.250820165748},
}

@article{LalondeMcDuff1996Ruled,
  author  = {Lalonde, Fran{\c{c}}ois and McDuff, Dusa},
  title   = {The classification of ruled symplectic
             {$4$}-manifolds},
  journal = {Mathematical Research Letters},
  volume  = {3},
  number  = {6},
  pages   = {769--778},
  year    = {1996},
  doi     = {10.4310/MRL.1996.v3.n6.a5}
}

@book {GSbook,
    AUTHOR = {Gompf, Robert E. and Stipsicz, Andr\'as I.},
     TITLE = {{$4$}-manifolds and {K}irby calculus},
    SERIES = {Graduate Studies in Mathematics},
    VOLUME = {20},
 PUBLISHER = {American Mathematical Society, Providence, RI},
      YEAR = {1999},
     PAGES = {xvi+558},
      ISBN = {0-8218-0994-6},
   MRCLASS = {57N13 (14J80 32Q55 57-02 57R17 57R57 57R65)},
  MRNUMBER = {1707327},
MRREVIEWER = {Nikolai\ N.\ Saveliev},
       DOI = {10.1090/gsm/020},
       URL = {https://doi.org/10.1090/gsm/020},
}

@article {Li2006Quaternionic,
    AUTHOR = {Li, Tian-Jun},
     TITLE = {Quaternionic bundles and {B}etti numbers of symplectic
              4-manifolds with {K}odaira dimension zero},
   JOURNAL = {Int. Math. Res. Not.},
  FJOURNAL = {International Mathematics Research Notices},
      YEAR = {2006},
     PAGES = {Art. ID 37385, 28},
      ISSN = {1073-7928,1687-0247},
   MRCLASS = {57R57 (53D35 57R17)},
  MRNUMBER = {2264722},
MRREVIEWER = {Tobias\ Ekholm},
       DOI = {10.1155/IMRN/2006/37385},
       URL = {https://doi.org/10.1155/IMRN/2006/37385},
}

@article {Ning2025,
    AUTHOR = {Ning, Shengzhen},
     TITLE = {Comparing {K}\"ahler cone and symplectic cone of one-point
              blowup of {E}nriques surface},
   JOURNAL = {J. Lond. Math. Soc. (2)},
  FJOURNAL = {Journal of the London Mathematical Society. Second Series},
    VOLUME = {111},
      YEAR = {2025},
    NUMBER = {3},
     PAGES = {Paper No. e70117, 23},
      ISSN = {0024-6107,1469-7750},
   MRCLASS = {57K43 (14J28)},
  MRNUMBER = {4877602},
MRREVIEWER = {Zhenbo\ Qin},
       DOI = {10.1112/jlms.70117},
       URL = {https://doi-org.proxy-um.researchport.umd.edu/10.1112/jlms.70117},
}

@article {SFTcompactness,
    AUTHOR = {Bourgeois, F. and Eliashberg, Y. and Hofer, H. and Wysocki, K.
              and Zehnder, E.},
     TITLE = {Compactness results in symplectic field theory},
   JOURNAL = {Geom. Topol.},
  FJOURNAL = {Geometry and Topology},
    VOLUME = {7},
      YEAR = {2003},
     PAGES = {799--888},
      ISSN = {1465-3060,1364-0380},
   MRCLASS = {53D45 (53D35 53D40 57R17)},
  MRNUMBER = {2026549},
MRREVIEWER = {Kai\ Cieliebak},
       DOI = {10.2140/gt.2003.7.799},
       URL = {https://doi.org/10.2140/gt.2003.7.799},
}

@article {MorganSzabo,
    AUTHOR = {Morgan, John W. and Szab\'o, Zolt\'an},
     TITLE = {Homotopy {$K3$} surfaces and mod {$2$} {S}eiberg-{W}itten
              invariants},
   JOURNAL = {Math. Res. Lett.},
  FJOURNAL = {Mathematical Research Letters},
    VOLUME = {4},
      YEAR = {1997},
    NUMBER = {1},
     PAGES = {17--21},
      ISSN = {1073-2780},
   MRCLASS = {57R57 (57R55)},
  MRNUMBER = {1432806},
MRREVIEWER = {Paolo\ Lisca},
       DOI = {10.4310/MRL.1997.v4.n1.a2},
       URL = {https://doi.org/10.4310/MRL.1997.v4.n1.a2},
}

@article {Gompf1998,
    AUTHOR = {Gompf, Robert E.},
     TITLE = {Symplectically aspherical manifolds with nontrivial {$\pi_2$}},
   JOURNAL = {Math. Res. Lett.},
  FJOURNAL = {Mathematical Research Letters},
    VOLUME = {5},
      YEAR = {1998},
    NUMBER = {5},
     PAGES = {599--603},
      ISSN = {1073-2780},
   MRCLASS = {57R57 (57R15 57R55)},
  MRNUMBER = {1666848},
MRREVIEWER = {Vicente\ Mu\~noz},
       DOI = {10.4310/MRL.1998.v5.n5.a4},
       URL = {https://doi-org.proxy-um.researchport.umd.edu/10.4310/MRL.1998.v5.n5.a4},
}

@phdthesis{Pasquotto2004,
  author = {Pasquotto, Federica},
  title  = {On the Geography of Symplectic Manifolds},
  school = {Leiden University},
  year   = {2004}
}

@article {Auroux2000,
    AUTHOR = {Auroux, Denis},
     TITLE = {Symplectic 4-manifolds as branched coverings of {$\bold C\bold
              P^2$}},
   JOURNAL = {Invent. Math.},
  FJOURNAL = {Inventiones Mathematicae},
    VOLUME = {139},
      YEAR = {2000},
    NUMBER = {3},
     PAGES = {551--602},
      ISSN = {0020-9910,1432-1297},
   MRCLASS = {53D35 (57R17)},
  MRNUMBER = {1738061},
MRREVIEWER = {Ignasi\ Mundet-Riera},
       DOI = {10.1007/s002220050019},
       URL = {https://doi-org.proxy-um.researchport.umd.edu/10.1007/s002220050019},
}

@article {Rohlin71,
    AUTHOR = {Rohlin, V. A.},
     TITLE = {Two-dimensional submanifolds of four-dimensional manifolds},
   JOURNAL = {Funkcional. Anal. i Prilo\v zen.},
  FJOURNAL = {Akademija Nauk SSSR. Funkcional\cprime nyi Analiz i ego
              Prilo\v zenija},
    VOLUME = {5},
      YEAR = {1971},
    NUMBER = {1},
     PAGES = {48--60},
      ISSN = {0374-1990},
   MRCLASS = {57D40},
  MRNUMBER = {298684},
MRREVIEWER = {G.\ Lusztig},
}

@book {BHPV,
    AUTHOR = {Barth, Wolf P. and Hulek, Klaus and Peters, Chris A. M. and
              Van de Ven, Antonius},
     TITLE = {Compact complex surfaces},
    SERIES = {},
    VOLUME = {},
   EDITION = {Second},
 PUBLISHER = {Springer-Verlag, Berlin},
      YEAR = {2004},
     PAGES = {xii+436},
      ISBN = {3-540-00832-2},
   MRCLASS = {14Jxx (14-02 32-02 32J15 57R57)},
  MRNUMBER = {2030225},
MRREVIEWER = {I.\ Dolgachev},
       DOI = {10.1007/978-3-642-57739-0},
       URL = {https://doi-org.proxy-um.researchport.umd.edu/10.1007/978-3-642-57739-0},
}

@article {Suwa,
    AUTHOR = {Suwa, Tatsuo},
     TITLE = {On hyperelliptic surfaces},
   JOURNAL = {J. Fac. Sci. Univ. Tokyo Sect. I},
  FJOURNAL = {Journal of the Faculty of Science. University of Tokyo.
              Section I},
    VOLUME = {16},
      YEAR = {1969},
     PAGES = {469--476},
      ISSN = {0368-2269},
   MRCLASS = {32.47 (14.00)},
  MRNUMBER = {267133},
MRREVIEWER = {Manfred\ Schneider},
}

@article {Geiges1992,
    AUTHOR = {Geiges, Hansj\"org},
     TITLE = {Symplectic structures on {$T^2$}-bundles over {$T^2$}},
   JOURNAL = {Duke Math. J.},
  FJOURNAL = {Duke Mathematical Journal},
    VOLUME = {67},
      YEAR = {1992},
    NUMBER = {3},
     PAGES = {539--555},
      ISSN = {0012-7094,1547-7398},
   MRCLASS = {57R22 (53C15 57R15 58F05)},
  MRNUMBER = {1181312},
MRREVIEWER = {Marisa\ Fern\'andez},
       DOI = {10.1215/S0012-7094-92-06721-4},
       URL = {https://doi-org.proxy-um.researchport.umd.edu/10.1215/S0012-7094-92-06721-4},
}

@article {TaubesConstraints,
    AUTHOR = {Taubes, Clifford Henry},
     TITLE = {More constraints on symplectic forms from {S}eiberg-{W}itten
              invariants},
   JOURNAL = {Math. Res. Lett.},
  FJOURNAL = {Mathematical Research Letters},
    VOLUME = {2},
      YEAR = {1995},
    NUMBER = {1},
     PAGES = {9--13},
      ISSN = {1073-2780},
   MRCLASS = {57R57 (57N13 57R15 58D29)},
  MRNUMBER = {1312973},
MRREVIEWER = {Jim\ A.\ Bryan},
       DOI = {10.4310/MRL.1995.v2.n1.a2},
       URL = {https://doi-org.proxy-um.researchport.umd.edu/10.4310/MRL.1995.v2.n1.a2},
}

@incollection {KirbyTaylor1990,
    AUTHOR = {Kirby, R. C. and Taylor, L. R.},
     TITLE = {{${\rm Pin}$} structures on low-dimensional manifolds},
 BOOKTITLE = {Geometry of low-dimensional manifolds, 2 ({D}urham, 1989)},
    SERIES = {London Math. Soc. Lecture Note Ser.},
    VOLUME = {151},
     PAGES = {177--242},
 PUBLISHER = {Cambridge Univ. Press, Cambridge},
      YEAR = {1990},
      ISBN = {0-521-40001-5},
   MRCLASS = {57R15 (57M25 57R90)},
  MRNUMBER = {1171915},
MRREVIEWER = {Elmer\ G.\ Rees},
}

@incollection {Matsumoto,
    AUTHOR = {Matsumoto, Yukio},
     TITLE = {An elementary proof of {R}ochlin's signature theorem and its
              extension by {G}uillou and {M}arin},
 BOOKTITLE = {\`A{} la recherche de la topologie perdue},
    SERIES = {Progr. Math.},
    VOLUME = {62},
     PAGES = {119--139},
 PUBLISHER = {Birkh\"auser Boston, Boston, MA},
      YEAR = {1986},
      ISBN = {0-8176-3329-4},
   MRCLASS = {57N13},
  MRNUMBER = {900248},
       DOI = {10.1007/BF01585163},
       URL = {https://doi-org.proxy-um.researchport.umd.edu/10.1007/BF01585163},
}

@article {ConwayOrsonPowell2025,
    AUTHOR = {Conway, Anthony and Orson, Patrick and Powell, Mark},
     TITLE = {Unknotting nonorientable surfaces},
   JOURNAL = {J. Eur. Math. Soc. (JEMS)},
      YEAR = {2025},
      NOTE = {Published online first},
       DOI = {10.4171/JEMS/1638},
       URL = {https://doi.org/10.4171/JEMS/1638},
}

@article {Audin,
    AUTHOR = {Audin, Mich\`ele},
     TITLE = {Quelques remarques sur les surfaces lagrangiennes de
              {G}ivental},
   JOURNAL = {J. Geom. Phys.},
  FJOURNAL = {Journal of Geometry and Physics},
    VOLUME = {7},
      YEAR = {1990},
    NUMBER = {4},
     PAGES = {583--598},
      ISSN = {0393-0440},
   MRCLASS = {57R10 (53C15 58E05)},
  MRNUMBER = {1131914},
MRREVIEWER = {Louis\ Funar},
       DOI = {10.1016/0393-0440(90)90008-Q},
       URL = {https://doi.org/10.1016/0393-0440(90)90008-Q},
}

@article {LiWu2012,
    AUTHOR = {Li, Tian-Jun and Wu, Weiwei},
     TITLE = {Lagrangian spheres, symplectic surfaces and the symplectic
              mapping class group},
   JOURNAL = {Geom. Topol.},
  FJOURNAL = {Geometry \& Topology},
    VOLUME = {16},
      YEAR = {2012},
    NUMBER = {2},
     PAGES = {1121--1169},
      ISSN = {1465-3060,1364-0380},
   MRCLASS = {53D05 (53D12 53D42)},
  MRNUMBER = {2946805},
MRREVIEWER = {Jonathan\ David\ Evans},
       DOI = {10.2140/gt.2012.16.1121},
       URL = {https://doi-org.proxy-um.researchport.umd.edu/10.2140/gt.2012.16.1121},
}

@article {Furuta2001,
    AUTHOR = {Furuta, M.},
     TITLE = {Monopole equation and the {$\frac{11}8$}-conjecture},
   JOURNAL = {Math. Res. Lett.},
  FJOURNAL = {Mathematical Research Letters},
    VOLUME = {8},
      YEAR = {2001},
    NUMBER = {3},
     PAGES = {279--291},
      ISSN = {1073-2780},
   MRCLASS = {57R15 (53C27 57R57)},
  MRNUMBER = {1839478},
       DOI = {10.4310/MRL.2001.v8.n3.a5},
       URL = {https://doi-org.proxy-um.researchport.umd.edu/10.4310/MRL.2001.v8.n3.a5},
}

@article {DaiHoLi2019,
    AUTHOR = {Dai, Bo and Ho, Chung-I and Li, Tian-Jun},
     TITLE = {Nonorientable {L}agrangian surfaces in rational 4-manifolds},
   JOURNAL = {Algebr. Geom. Topol.},
  FJOURNAL = {Algebraic \& Geometric Topology},
    VOLUME = {19},
      YEAR = {2019},
    NUMBER = {6},
     PAGES = {2837--2854},
      ISSN = {1472-2747,1472-2739},
   MRCLASS = {53D12 (57Q35)},
  MRNUMBER = {4023330},
MRREVIEWER = {Roman\ Golovko},
       DOI = {10.2140/agt.2019.19.2837},
       URL = {https://doi-org.proxy-um.researchport.umd.edu/10.2140/agt.2019.19.2837},
}

@article {LiLi2002,
    AUTHOR = {Li, Bang-He and Li, Tian-Jun},
     TITLE = {Symplectic genus, minimal genus and diffeomorphisms},
   JOURNAL = {Asian J. Math.},
  FJOURNAL = {Asian Journal of Mathematics},
    VOLUME = {6},
      YEAR = {2002},
    NUMBER = {1},
     PAGES = {123--144},
      ISSN = {1093-6106,1945-0036},
   MRCLASS = {57R17 (57R57 57R95)},
  MRNUMBER = {1902650},
MRREVIEWER = {Ivan\ Smith},
       DOI = {10.4310/AJM.2002.v6.n1.a7},
       URL = {https://doi-org.proxy-um.researchport.umd.edu/10.4310/AJM.2002.v6.n1.a7},
}

@article {Wall1964,
    AUTHOR = {Wall, C. T. C.},
     TITLE = {Diffeomorphisms of {$4$}-manifolds},
   JOURNAL = {J. London Math. Soc.},
  FJOURNAL = {The Journal of the London Mathematical Society},
    VOLUME = {39},
      YEAR = {1964},
     PAGES = {131--140},
      ISSN = {0024-6107,1469-7750},
   MRCLASS = {57.20},
  MRNUMBER = {163323},
MRREVIEWER = {N.\ Kuiper},
       DOI = {10.1112/jlms/s1-39.1.131},
       URL = {https://doi-org.proxy-um.researchport.umd.edu/10.1112/jlms/s1-39.1.131},
}

@article {Traynor,
    AUTHOR = {Traynor, Lisa},
     TITLE = {Symplectic packing constructions},
   JOURNAL = {J. Differential Geom.},
  FJOURNAL = {Journal of Differential Geometry},
    VOLUME = {42},
      YEAR = {1995},
    NUMBER = {2},
     PAGES = {411--429},
      ISSN = {0022-040X,1945-743X},
   MRCLASS = {53C15 (57R40 58F05)},
  MRNUMBER = {1366550},
MRREVIEWER = {Hansj\"org\ Geiges},
       URL = {http://projecteuclid.org.proxy-um.researchport.umd.edu/euclid.jdg/1214457236},
}

@incollection {Symington,
    AUTHOR = {Symington, Margaret},
     TITLE = {Four dimensions from two in symplectic topology},
 BOOKTITLE = {Topology and geometry of manifolds ({A}thens, {GA}, 2001)},
    SERIES = {Proc. Sympos. Pure Math.},
    VOLUME = {71},
     PAGES = {153--208},
 PUBLISHER = {Amer. Math. Soc., Providence, RI},
      YEAR = {2003},
      ISBN = {0-8218-3507-6},
   MRCLASS = {53D35 (53D20 55R55 57R17)},
  MRNUMBER = {2024634},
MRREVIEWER = {Vicente\ Mu\~noz},
       DOI = {10.1090/pspum/071/2024634},
       URL = {https://doi-org.proxy-um.researchport.umd.edu/10.1090/pspum/071/2024634},
}

@article {HK19,
    AUTHOR = {Holm, Tara S. and Kessler, Liat},
     TITLE = {Circle actions on symplectic four-manifolds},
   JOURNAL = {Comm. Anal. Geom.},
  FJOURNAL = {Communications in Analysis and Geometry},
    VOLUME = {27},
      YEAR = {2019},
    NUMBER = {2},
     PAGES = {421--464},
      ISSN = {1019-8385,1944-9992},
   MRCLASS = {53D35 (53D45)},
  MRNUMBER = {4003013},
MRREVIEWER = {Sonja\ Hohloch},
       DOI = {10.4310/CAG.2019.v27.n2.a6},
       URL = {https://doi-org.proxy-um.researchport.umd.edu/10.4310/CAG.2019.v27.n2.a6},
}

@article {DorfmeisterLi2010,
    AUTHOR = {Dorfmeister, Josef G. and Li, Tian-Jun},
     TITLE = {The relative symplectic cone and {$T^2$}-fibrations},
   JOURNAL = {J. Symplectic Geom.},
  FJOURNAL = {The Journal of Symplectic Geometry},
    VOLUME = {8},
      YEAR = {2010},
    NUMBER = {1},
     PAGES = {1--35},
      ISSN = {1527-5256,1540-2347},
   MRCLASS = {53D35 (57R17 57R57)},
  MRNUMBER = {2609626},
MRREVIEWER = {Timothy\ Perutz},
       DOI = {10.4310/jsg.2010.v8.n1.a1},
       URL = {https://doi.org/10.4310/jsg.2010.v8.n1.a1},
}

@article {Usher2009,
    AUTHOR = {Usher, Michael},
     TITLE = {Kodaira dimension and symplectic sums},
   JOURNAL = {Comment. Math. Helv.},
  FJOURNAL = {Commentarii Mathematici Helvetici. A Journal of the Swiss
              Mathematical Society},
    VOLUME = {84},
      YEAR = {2009},
    NUMBER = {1},
     PAGES = {57--85},
      ISSN = {0010-2571,1420-8946},
   MRCLASS = {57R17 (53D35 57R57)},
  MRNUMBER = {2466075},
MRREVIEWER = {Timothy\ Perutz},
       DOI = {10.4171/CMH/152},
       URL = {https://doi.org/10.4171/CMH/152},
}
\bigskip
\noindent
\textsc{Department of Mathematics, Rutgers University, Piscataway, NJ, U.S.A.}
\medskip
\noindent
\\
\textit{Email address}: \texttt{xinle.dai@rutgers.edu}

\medskip
\noindent
\textsc{Department of Mathematics, University of Southern California, Los Angeles, CA, U.S.A.}
\medskip
\noindent
\\
\textit{Email address: }\texttt{boxihao@usc.edu}

\medskip
\noindent
\textsc{Department of Mathematics, University of Maryland, College Park, MD, U.S.A}
\\
\textit{Email address}: \texttt{szning@umd.edu}
\end{document}